\documentclass[a4paper,12pt, onecolumn, twoside]{book}

\usepackage[latin1]{inputenc} 
\usepackage[T1]{fontenc}
\usepackage[francais,frenchb]{babel}
\usepackage[fixlanguage]{babelbib}
\selectbiblanguage{french}
\usepackage[top=3cm, bottom=3cm, left=3cm, right=3cm]{geometry}
\usepackage{amsmath, amsfonts, amssymb}
\usepackage{mathrsfs}
\usepackage[mathscr]{eucal}
\usepackage{multicol}
\usepackage{graphicx}
\usepackage[thmmarks,amsmath]{ntheorem}
\usepackage{fancyhdr}
\usepackage{pstricks,pstricks-add,pst-plot}
\usepackage{cancel}
\usepackage{enumerate}
\usepackage{array}
\usepackage{tocbibind} 
\usepackage{url}
\usepackage{tikz-cd} 
\usepackage{calrsfs} 
\usepackage{xcolor}
\usepackage{fancybox}
\usepackage{tabvar}
\usepackage{eurosym}
\usepackage{listings} 
\usepackage{frcursive} 
\usepackage{calligra} 
\usepackage{caption} 
\usepackage{float} 
\usepackage{array,multirow,makecell} 
\usepackage{manfnt} 
\usepackage{lscape} 
\usepackage{pdfpages}  
\usepackage[new]{old-arrows} 
\usepackage{multirow} 
\DecimalMathComma

\renewcommand{\d}{\mathrm{d}}

\DeclareMathOperator{\pgcd}{pgcd}

\DeclareMathOperator{\vect}{Vect}

\newcommand{\id}{\mathrm{id}}
\DeclareMathOperator{\rg}{rg}

\DeclareMathOperator{\com}{com}

\renewcommand{\arcsin}{\mathrm{Arcsin}}

\renewcommand{\pgcd}{\mathrm{pgcd}}
\DeclareMathOperator{\card}{Card}

\DeclareMathOperator{\GL}{GL}

\DeclareMathOperator{\M}{M}

\DeclareMathOperator{\Gr}{Gr}
\DeclareMathOperator{\covol}{covol}
\DeclareMathOperator{\vol}{vol}

\DeclareMathOperator{\codim}{codim}

\newcommand{\Q}{\ensuremath{\mathbb{Q}}}
\newcommand{\R}{\ensuremath{\mathbb{R}}}

\newcommand{\N}{\ensuremath{\mathbb{N}}}
\newcommand{\Z}{\ensuremath{\mathbb{Z}}}

\newcommand{\A}{\ensuremath{\mathbb{A}}}

\renewcommand{\H}{\ensuremath{\mathbb{H}}}

\renewcommand{\P}{\ensuremath{\mathbb{P}}}

\newcommand{\1}{\ensuremath{\mathbf{1}}}

\renewcommand{\epsilon}{\varepsilon}
\renewcommand{\le}{\leqslant}
\renewcommand{\ge}{\geqslant}
\newcommand{\transp}{\,{}^t\!}

\newcommand{\fonction}[5]{\begin{array}{lrcl}
#1\colon & #2 & \longrightarrow & #3 \\
    & #4 & \longmapsto &  #5  \end{array}}

\newcommand{\abs}[1]{\left\vert #1 \right\vert}
\newcommand{\norme}[1]{\left\Vert #1 \right\Vert}

\newcommand{\muexp}[4]{\mathring\mu_{#1}(#2\vert #3)_{#4}}
\newcommand{\muexpA}[4]{\mu_{#1}(#2\vert #3)_{#4}}

\theorembodyfont{\upshape}     
\newtheorem{theoreme}{Th\'eor\`eme}[chapter]
\newtheorem*{theoremenonum}{Th\'eor\`eme}

\newtheorem{conjecture}[theoreme]{Conjecture}
\newtheorem{proposition}[theoreme]{Proposition}

\newtheorem{lemme}[theoreme]{Lemme}

\newtheorem{corollaire}[theoreme]{Corollaire}
\newtheorem*{corollairenonum}[theoreme]{Corollaire}
\newtheorem{notation}[theoreme]{Notation}

\theorembodyfont{\upshape} 
\newtheorem{exemple}[theoreme]{Exemple}

\newtheorem{remarque}[theoreme]{Remarque}

\newtheorem{definition}[theoreme]{D\'efinition}

\newtheorem{probleme}[theoreme]{Problème}

\theoremstyle{nonumberbreak} 
\theoremheaderfont{\itshape}
\theoremsymbol{\ensuremath{\square}} 

\newtheorem{preuve}{Preuve.}

\psset{algebraic=true}
\usepackage{xypic}

\title{Approximation rationnelle de sous-espaces vectoriels}
\author{Elio Joseph}
\date{}

\usepackage{float}

\begin{document}

\includepdf[pages=-]{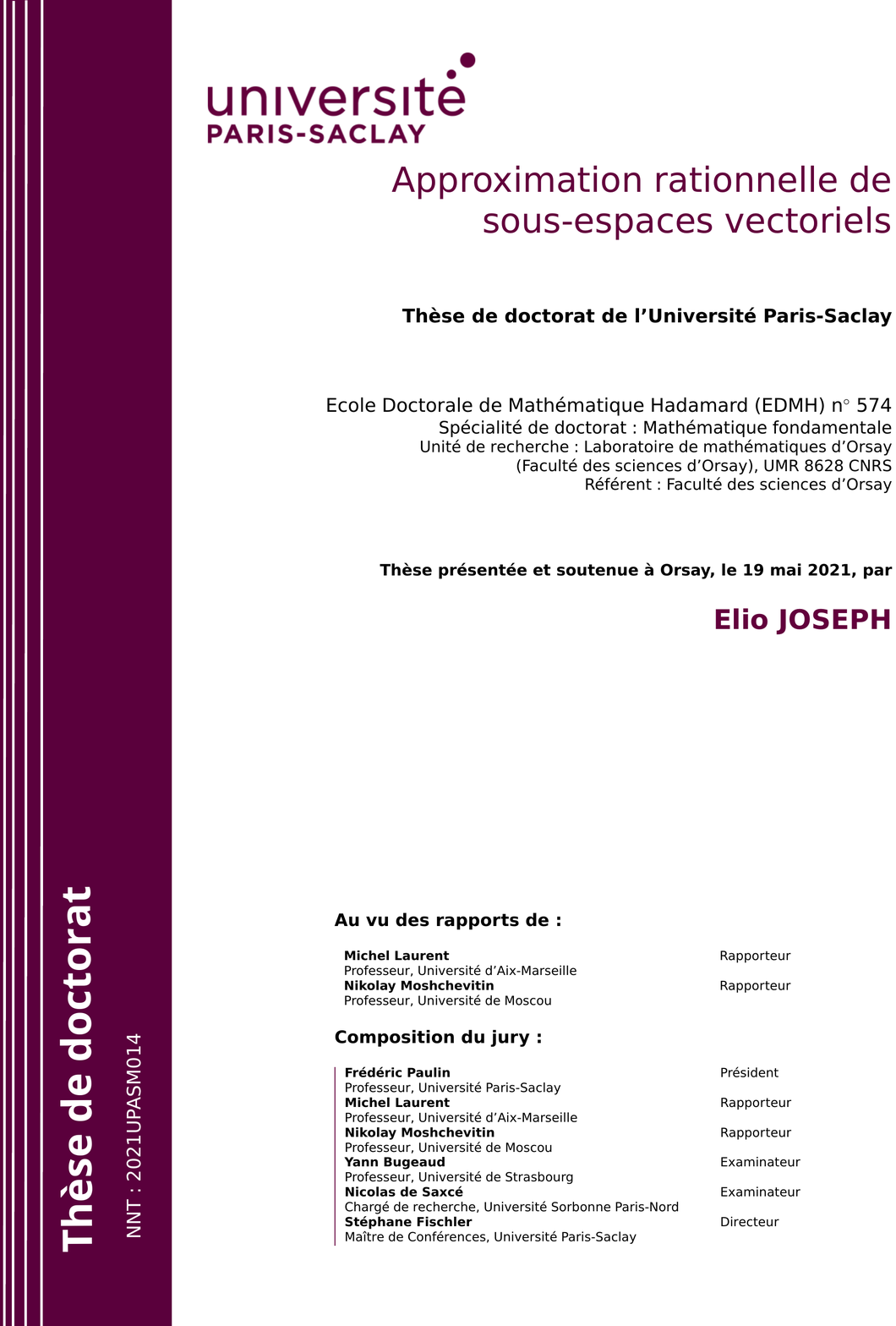}

\newpage
~
\pagestyle{empty}

\cleardoublepage

~
\vspace{3.5cm}
\begin{flushright}
\textit{À Michel Fournier.}
\end{flushright}

\chapter*{Remerciements}

En tout premier lieu, j'adresse de sincères remerciements à mon directeur de thèse, Stéphane Fischler. Merci pour la disponibilité constante, merci pour toutes les relectures aussi nombreuses qu'attentives, merci pour la confiance accordée en me laissant de l'autonomie dans le travail de recherche, et enfin merci pour tous les conseils avisés qui m'éclairèrent pendant ces trois années. \\

Je remercie chaleureusement les deux rapporteurs de cette thèse, Michel Laurent et Nikolay Moshchevitin, pour avoir généreusement donné de leur temps pour relire cette thèse. Merci pour votre intérêt et vos retours précieux. Un grand merci également à Yann Bugeaud, Frédéric Paulin et Nicolas de Saxcé pour avoir accepté de faire partie du jury, je vous en suis très reconnaissant. \\

Un grand merci à tous les chercheurs, enseignants et doctorants que j'ai pu croiser pendant mon doctorat, spécifiquement ceux que j'ai croisés à Cambridge, Istanbul et Orsay. Je remercie tout spécialement les doctorants du Laboratoire de Mathématiques d'Orsay. De peur d'en oublier, je ne prends pas le risque de les citer nommément, excepté mes valeureux co-bureau : Cyril Falcon, Hugo Federico, Guillaume Maillard et Adrien Béguinet. \\

Merci à Cyril Falcon pour avoir lancé avec moi le séminaire \emph{Explique-moi...}, merci à Nathalie Carrierre pour l'aide précieuse à l'organisation, et merci à Ella Blair et Adrien Béguinet pour l'avoir repris fin 2020. \\

J'adresse toute ma reconnaissance à Stéphane Fischler, Lucie Flammarion et Boris Joseph pour leurs relectures attentives de ce manuscrit. Toute erreur qui subsisterait reste évidemment de mon fait. \\

Merci à Yoann Pouligo, Louis Vialle, Hugo Guichaoua, Hugo Lefèvre, Morgan Catez, Stéphane Zlammansuck, Vincent Galbrun, Benoît Vacher, Anthony de Oliveira, Hippolyte Beaudroit, Matthieu Perrin, Maxime Caramona, Alexandre de Lemos, Donovan Dugeny, Benoît Sabourin, Florian Dangleant, Clément Labadie, Yoshi et les autres pour les drôles de moments pendant la préparation de cette thèse, ainsi que pour la convivialité pendant les corrections de copies. \\

Je suis reconnaissant aux professeurs qui auront, chacun à leur façon, contribué à la poursuite de mes études mathématiques. Pour n'en citer que quelques-uns, je remercie François Bertholon, Bruno Arsac, Etienne Fouvry, et bien sûr Olivier Fouquet. \\

C'est grâce à leur présence bienveillante et à leur complicité enjouée que la vie fut si joyeuse pendant ces trois dernières années. Je remercie tout particulièrement et par ordre alphabétique Alexandre Barrat pour m'avoir allégé d'une grande tâche, Charles Bignaud et Raphaël Huille pour la vie dans la forêt, Clément Jean et Lauren Oliel pour les apéros bios, Clément Walter pour les apéros-fléchettes, Cyril Falcon pour les nombreuses et diverses aventures, Camille Masson, Faustin Besiers, Joris Baraillon et Julien Crémy pour les vacances colorées, Florian Granger et Gédéon Chevallier pour les blagues et les grandes nouvelles, Hugo Federico pour les parties d'échecs, Lucie Flammarion pour toute la vie, Marguerite et Guillaume Matheron pour leur céleste présence, Olivier Rembliere pour tout le sport, Sandrine Gauthier pour les repas toujours joyeux, Thomas Gastellu pour les gardes et Victoria Soubeiran pour les nouvelles d'outre-mer. \\

Je remercie chaleureusement tous ceux qui ont émis l'idée de venir assister à ma soutenance ou qui ont proposé leur aide pour le pot. \\

Je remercie du fond du c\oe ur ma famille, pour leur soutien et leur présence attentive. Plus particulièrement ma s\oe ur Fanny et mes parents Anne et Boris. \\

Enfin, j'adresse non sans malice un merci tout particulier à ma chère fiancée et complice Lucie Flammarion. Cette thèse est évidemment trop courte pour te remercier entièrement, mais je souligne ici ton rayonnement quotidien, ton humour joyeux et tes conseils éclairés. Chère Lucie, 9\emph{0}4\emph{19}5\reflectbox{\emph{3}}10\emph{9}\reflectbox{10}\emph{14} !

\renewcommand{\headrulewidth}{0pt}
\pagestyle{fancy}

\tableofcontents{}
\chapter{Introduction}

La thèse qui suit s'inscrit dans le domaine de l'approximation diophantienne. On expose brièvement ce domaine dans la section \ref{section_approximation_diophantienne_classique_moairgnibva}, puis on introduit le problème de l'approximation diophantienne de sous-espaces vectoriels de $\R^n$ dans la section \ref{section_approx_sev_maobinaiovnoivb}. La section \ref{section_presentation_resultats_aboininbv} est consacrée à une présentation des résultats connus de ce domaine. Enfin, un plan de la thèse et des résultats principaux obtenus se trouve dans la section \ref{section_plan_de_la_these_abinrebinv}, une présentation explicite des résultats connus et obtenus jusqu'en dimension $6$ est donnée dans la section \ref{section_exemples_pour_n_inferieur_6_maobuefobv}, et une conjecture est formulée dans la section \ref{section_conjecture_oabeneboinv}.

\section{Approximation diophantienne rationnelle classique}\label{section_approximation_diophantienne_classique_moairgnibva}

On introduit ici le domaine de l'approximation diophantienne avec un double objectif. Le premier objectif est d'exposer de manière large l'endroit où se situe ce travail dans le monde mathématique, le second est de préparer l'introduction de notions qui serviront de façon centrale dans cette thèse. Ces notions pourront alors être comparées à celles vues dans cette section en approximation diophantienne rationnelle classique. \\

L'approximation diophantienne a pour but originel d'approcher les nombres réels par des nombres rationnels. Plus précisément, pour un réel $\xi$ donné, on cherche un rationnel $p/q$ tel que la quantité 

\[\abs{\xi-\frac pq}\]

soit la plus petite possible. Par densité de $\Q$ dans $\R$, ce problème n'a pas d'intérêt en l'état, cette quantité pouvant être rendue arbitrairement petite pour tout réel $\xi$. \\

On précise donc la question en demandant que le rationnel $p/q$ approchant le réel $\xi$ soit le moins \emph{compliqué} possible. Plus précisément, on veut lier la \emph{qualité} de l'approximation $\abs{\xi-p/q}$ à la \emph{complexité} du rationnel. Ainsi, s'autoriser des rationnels plus compliqués impose d'obtenir en retour une approximation plus précise. \\

Écrivons plus rigoureusement ceci. On donne une notion de \emph{complexité} d'un rationnel en regardant la taille de son dénominateur. On fixe alors un réel $\xi\in\R$, et on cherche les réels $\mu>0$ tels qu'il existe une infinité de rationnels $p/q$ vérifiant 

\begin{equation}\label{approx_dioph_classique_riogberg}
\abs{\xi-\frac pq}< \frac{1}{q^\mu}.
\vspace{2.6mm}
\end{equation}

Un nombre émerge naturellement de ce problème : l'\emph{exposant d'irrationalité} de $\xi$. Celui-ci est défini comme la borne supérieure de l'ensemble des réels $\mu>0$ vérifiant l'inégalité \eqref{approx_dioph_classique_riogberg} pour une infinité de nombres rationnels $p/q$. On notera par la suite $\mu(\xi)$ cette borne supérieure (éventuellement infinie). C'est ce nombre qui est au coeur de la théorie de l'approximation diophantienne. \\

De nombreux résultats sont connus sur cette quantité, mais des conjectures demeurent. Par exemple, le théorème de Dirichlet (19\up{ème} siècle) affirme que si $\xi$ est irrationnel, alors pour tout entier $Q\ge 1$, il existe deux entiers $p$ et $q$, avec $q\in\{1,\ldots,Q\}$, tels que

\[\abs{\xi-\frac pq}\le \frac 1{qQ}.\]

De ceci on déduit immédiatement que pour tout $\xi\notin\Q$, $\mu(\xi)\ge 2$. On peut aussi énoncer un théorème profond (théorème page 2 de \cite{roth55}) que K. F. Roth a démontré en 1955 en s'appuyant sur des travaux de Thue, Siegel et Dyson : si $\xi$ est un nombre algébrique de degré au moins $2$, alors $\mu(\xi)=2$. La réciproque étant fausse : on peut par exemple montrer que $\mu(e)=2$. \\ 

Enfin, on mentionne que c'est la théorie des fractions continues -- à laquelle une introduction claire peut être trouvée dans le livre \cite{hardy07} de G. H. Hardy et \hbox{E. M. Wright} -- qui permet de répondre à de nombreuses questions d'approximation diophantienne. Les fractions continues permettent notamment de montrer que

\[\Big\{\mu(\xi),\hspace{1mm}\xi\in\R\Big\}=\{1\}\cup[2,+\infty].\]

Les nombres $\xi$ tels que $\mu(\xi)=+\infty$ sont appelés \emph{nombres de Liouville}, et ce sont eux qui fournirent le premier nombre transcendant défini explicitement (\cite{liouville44} pages 910--911). \\

Le lecteur qui souhaite approfondir le sujet de l'approximation diophantienne rationnelle classique peut consulter \cite{niven63}, \cite{cassels57} ou encore \cite{schmidt80}, et \cite{ten96} au sujet de l'approximation diophantienne de matrices.

\section{Approximation diophantienne de sous-espaces vectoriels de $\R^n$}\label{section_approx_sev_maobinaiovnoivb}

Cette thèse suit l'idée donnée par W. M. Schmidt en 1967 dans un article fondateur \cite{schmidt67}. L'idée consiste à généraliser le problème de l'approximation diophantienne rationnelle classique qui a été introduit en section \ref{section_approximation_diophantienne_classique_moairgnibva}, pour créer une théorie de l'approximation diophantienne des sous-espaces vectoriels de $\R^n$. \\

Il est alors intéressant de raisonner de façon similaire à l'approximation diophantienne classique. Soient $n\ge 2$ et $d,e\in\{1,\ldots,n-1\}$ tels que $d+e\le n$. Soit $A$ un sous-espace vectoriel de $\R^n$ de dimension $d$. Le but est de chercher à \emph{approcher} l'espace $A$ par des sous-espaces \emph{rationnels} $B$ de dimension $e$ "pas trop" \emph{compliqués}. Pour formaliser ce problème, il reste donc à définir les termes \emph{rationnel}, \emph{approcher} et \emph{compliqué}. 

\begin{remarque}

Comme le fait Schmidt dans \cite{schmidt67} à partir de la section III, on supposera toujours dans cette thèse que $d+e\le n$, avec $d$ la dimension du sous-espace approché et $e$ la dimension des sous-espaces rationnels approchants. 

\end{remarque}

\begin{definition}\label{def_sev_rationnel_amebvomabsd}

Un sous-espace vectoriel $B$ de $\R^n$ est dit \emph{rationnel} s'il admet une base formée de vecteurs à coordonnées rationnelles. On note $\mathfrak R_n(e)$ l'ensemble des sous-espaces rationnels de $\R^n$ de dimension $e$.

\end{definition}

Comme dans le cas de l'approximation des irrationnels par les rationnels, il est utile de définir une notion d'\emph{irrationalité} pour le sous-espace $A$ qu'on cherche à approcher. Schmidt en donne une dans le corollaire du théorème 12 page 459 de \cite{schmidt67}.

\begin{definition}\label{def_sev_irrationnel_areiogahmvi}

Soient $d,e\in\{1,\ldots,n-1\}$ tels que $d+e\le n$ ; on pose $t=\min(d,e)$. Soit $j\in\{1,\ldots,t\}$. Un sous-espace $A$ de dimension $d$ de $\R^n$ est dit \emph{$(e,j)$-irrationnel} si pour tout sous-espace rationnel $B$ de dimension $e$, 

\[\dim(A\cap B)<j.\]

On note $\mathfrak I_n(d,e)_j$ l'ensemble des sous-espaces $(e,j)$-irrationnels de dimension $d$ de $\R^n$.

\end{definition}

\begin{remarque}\label{remarque_condition_dirr_abrihfbamvoa}

On utilisera souvent la condition "être $(e,1)$-irrationnel", précisons donc celle-ci. Un sous-espace $A$ est $(e,1)$-irrationnel si, et seulement si, il intersecte trivialement tous les sous-espace vectoriels rationnels de dimension $e$. Le fait que $A$ soit $(e,1)$-irrationnel est plus fort que le fait qu'il ne contienne pas de vecteur rationnel. Le sous-espace $A$ est $(e,1)$-irrationnel si, et seulement si,

\[\forall(\xi_1,\ldots,\xi_n)\in A\setminus\{0\},\quad \dim_\Q\vect_\Q(\xi_1,\ldots,\xi_n)\ge e+1.\]

\end{remarque}

On a commencé à formuler le problème d'approximation diophantienne de sous-espaces vectoriels en définissant les notions de rationalité et d'irrationalité pour ceux-ci. Il reste à formuler les notions de \emph{proximité} et de \emph{complexité} pour des sous-espaces. \\

Commençons par la \emph{complexité}, notion définie dans \cite{schmidt67} pages 432--433. Pour cela, on a besoin des coordonnées de Plücker (définition \ref{defcoordgrassmann_reaoimghermoh}), auxquelles une brève introduction est exposée dans la sous-section \ref{coord_Pluck_teoghzleufve}, basée sur le livre \cite{caldero15} de \hbox{P. Caldero} et \hbox{J. Germoni}. \\

Sauf mention explicite du contraire, ici et dans toute la suite, la norme $\norme\cdot$ désignera la norme euclidienne canonique. \\

\begin{definition}\label{def_hauteur_coord_plucker_reuighgzmne}

Soit $B\in\mathfrak R_n(e)$. Posons $N=\binom ne$ et notons $\Xi=(\xi_1,\ldots,\xi_N)$ le vecteur d'un représentant des coordonnées de Plücker de $B$. Comme $B$ est un sous-espace rationnel, on peut choisir $\Xi$ à coordonnées entières et premières entre elles. On définit alors la \emph{hauteur} de $B$ comme

\[H(B)=\norme{\Xi}=\sqrt{\sum_{i=1}^N \xi_i^2}.\]

\end{definition}

Ainsi, de même que la complexité d'un nombre rationnel se mesure grâce à la taille de son dénominateur sous forme irréductible, la complexité d'un sous-espace rationnel se mesure grâce à sa hauteur.

\begin{remarque}\label{def_plus_generale_de_la_hauteur_aobmirnfmeobnv}

Si $\Xi=(\xi_1,\ldots,\xi_N)$ est un représentant quelconque des coordonnées de Plücker, notons $\mathfrak a$ l'idéal fractionnaire de $\Z$ engendré par les $\xi_i$ : $\mathfrak a=\xi_1\Z+\cdots+\xi_N\Z$. L'équation (1) page 432 de \cite{schmidt67} permet de généraliser la définition \ref{def_hauteur_coord_plucker_reuighgzmne} : \hbox{$H(B)=\norme{\Xi}/N(\mathfrak a)$}.

\end{remarque}

Il reste enfin à définir une notion de \emph{proximité} entre deux sous-espaces de $\R^n$, pas nécessairement de même dimension. On pose de nouveau $t=\min(d,e)$. Cette notion de proximité vient encore de l'article \cite{schmidt67}, page 443, où Schmidt définit -- en s'appuyant notamment sur l'article \cite{seidel55} de J. J. Seidel -- $t$ angles entre les sous-espaces $A$ et $B$ de la façon suivante. \\

On munit $\R^n$ ainsi que la puissance extérieure $\Lambda^2(\R^n)$ de leurs normes euclidiennes canoniques respectives. On note pour $X,Y\in\R^n\setminus\{0\}$,

\[\psi(X,Y)=\sin\widehat{(X,Y)}=\frac{\norme{X\wedge Y}}{\norme X\cdot\norme Y}\]

où $\widehat{(X,Y)}$ désigne l'angle géométrique entre les vecteurs $X$ et $Y$, et $\wedge$ le produit extérieur sur $\R^n$ à valeurs dans $\Lambda^2(\R^n)$. \\

On peut alors définir 

\begin{equation}\label{def_du_psi_1_piaeofdbouefbv}
\psi_1(A,B)=\min_{\substack{X\in A\setminus\{0\} \\ Y\in B\setminus\{0\}}}\psi(X,Y)
\vspace{2.6mm}
\end{equation}

et on note $X_1,Y_1$ des vecteurs unitaires réalisant ce minimum. \\

On construit alors d'autres quantités $\psi_2(A,B),\ldots,\psi_t(A,B)$\label{page_def_psij_aomebfvoaub} par récurrence sur \hbox{$j\in\{1,\ldots,t\}$}.

Soit $j\in\{1,\ldots,t-1\}$, on suppose que $\psi_1(A,B),\ldots,\psi_j(A,B)$ ont été construites, associées à des couples de vecteurs $(X_1,Y_1),\ldots,(X_j,Y_j)\in A\times B$. On note $A_j$ le supplémentaire orthogonal de $\vect(X_1,\ldots,X_j)$ dans $A$ et $B_j$ le supplémentaire orthogonal de $\vect(Y_1,\ldots,Y_j)$ dans $B$. On définit de la même façon

\[\psi_{j+1}(A,B)=\min_{\substack{X\in A_j\setminus\{0\} \\ Y\in B_j\setminus\{0\}}}\psi(X,Y)\]

et on note $X_{j+1},Y_{j+1}$ des vecteurs unitaires réalisant ce minimum.

\begin{definition}\label{def_angles_canoniques_eormighemoidvn}

Pour $j\in\{1,\ldots,t\}$ on pose $\theta_j=\arcsin(\psi_j(A,B))$ et on appelle $\theta_j$ le \emph{$j$-ème angle canonique} entre $A$ et $B$.

\end{definition}

Ces angles sont canoniques dans le sens du théorème \ref{la_raison_pour_laquelle_les_angles_sont_canoniques} suivant, qui est le théorème 4 page 443 de l'article \cite{schmidt67}. Dans ce qui suit, le point $\cdot$ désigne le produit scalaire usuel sur $\R^n$.

\begin{theoreme}\label{la_raison_pour_laquelle_les_angles_sont_canoniques}

Il existe des bases orthonormales $(X_1,\ldots,X_d)$ et $(Y_1,\ldots,Y_e)$ de $A$ et $B$ respectivement, et il existe des réels $0\le \theta_t\le \cdots\le \theta_1\le 1$ tels que

\[\forall (i,j)\in\{1,\ldots,d\}\times\{1,\ldots,e\},\quad X_i\cdot Y_j=\delta_{i,j}\cos\theta_i\]

où $\delta$ est le symbole de Kronecker.

De plus, les nombres $\theta_1,\ldots,\theta_t$ sont indépendants des bases $(X_1,\ldots,X_d)$ et $(Y_1,\ldots,Y_e)$ vérifiant cette propriété, et sont invariants sous l'effet d'une isométrie appliquée simultanément à $A$ et à $B$. Enfin, les réels $\theta_j$ sont les angles canoniques de la définition \ref{def_angles_canoniques_eormighemoidvn}, et les vecteurs $X_1,\ldots,X_t,Y_1,\ldots,Y_t$ construits avant la définition \ref{def_angles_canoniques_eormighemoidvn} peuvent être complétés pour former de telles bases orthonormales.

\end{theoreme}

Comme $\psi_j(A,B)=\sin(\theta_j)$, ces bases orthonormales vérifient

\[\forall (i,j)\in\{1,\ldots,d\}\times\{1,\ldots,e\},\quad \psi(X_i,Y_j)=\begin{cases} 1&\text{ si $i\ne j$}\\\psi_j(A,B) &\text{ si $i=j$}.\end{cases}\]

\begin{remarque}

On utilisera à plusieurs reprises le fait que 

\[\dim(A\cap B)\ge j\iff \psi_j(A,B)=0.\]

Ceci permet de formuler une définition équivalente à la définition \ref{def_sev_irrationnel_areiogahmvi} :

\[A\in\mathfrak I_n(d,e)_j\iff \forall B\in\mathfrak R_n(e),\quad \psi_j(A,B)\ne 0.\]

\end{remarque}

Avec les définitions \ref{def_sev_rationnel_amebvomabsd}, \ref{def_sev_irrationnel_areiogahmvi}, \ref{def_hauteur_coord_plucker_reuighgzmne} et avec la définition des $\psi_j$ page \pageref{page_def_psij_aomebfvoaub}, il est désormais possible d'énoncer un problème analogue à l'inéquation \eqref{approx_dioph_classique_riogberg} en approximation diophantienne classique, mais cette fois-ci avec des sous-espaces vectoriels de $\R^n$. C'est ce qui est fait page 459 de \cite{schmidt67}. \\

Soient $n\ge 2$ et $d,e\in\{1,\ldots,n-1\}$ des entiers tels que $d+e\le n$, soit \hbox{$j\in\{1,\ldots,\min(d,e)\}$}. Soit $A\in\mathfrak I_n(d,e)_j$ un sous-espace $(e,j)$-irrationnel de $\R^n$. On cherche les réels $\beta>0$ tels qu'il existe une infinité de sous-espaces rationnels $B\in\mathfrak R_n(e)$ tels que

\begin{equation}\label{inequation_du_probleme_principal_oubrzgurb}
\psi_j(A,B)\le \frac{1}{H(B)^\beta}.
\vspace{2.6mm}
\end{equation}

On peut alors observer la similarité entre les inéquations \eqref{approx_dioph_classique_riogberg} et \eqref{inequation_du_probleme_principal_oubrzgurb}. \\

Définissons ensuite deux quantités importantes.

\begin{notation}

Soit $A\in\mathfrak I_n(d,e)_j$. Notons $\muexpA nAej$ la borne supérieure de l'ensemble des $\beta>0$ tels que l'inéquation \eqref{inequation_du_probleme_principal_oubrzgurb} soit vérifiée pour une infinité de $B\in\mathfrak R_n(e)$ ; si cet ensemble n'est pas majoré, on pose $\muexpA nAej=+\infty$.

Notons aussi

\[\muexp ndej=\inf_{A\in\mathfrak I_n(d,e)_j}\muexpA nAej.\]

\end{notation}

Ainsi, on peut formuler trois problèmes : 

\begin{probleme}\label{probleme_sev_particulier_nreourbgea}

Pour $A\in\mathfrak I_n(d,e)_j$, déterminer $\muexpA nAej$ en fonction de $A$ et $(n,e,j)$.

\end{probleme}

\begin{probleme}\label{probleme_principal_gnzroign}

Déterminer $\muexp ndej$ en fonction de $(n,d,e,j)$.

\end{probleme}

\begin{probleme}\label{probleme_spectre_oaerinemoivn}

Déterminer l'ensemble $\muexpA n{\mathfrak I_n(d,e)_j}ej$ en fonction de $(n,d,e,j)$, \emph{i.e.} l'ensemble des valeurs prises par $\muexpA nAej$ lorsque $A$ décrit $\mathfrak I_n(d,e)_j$.

\end{probleme}

Continuons la comparaison avec l'approximation diophantienne rationnelle classique :
\begin{enumerate}[$\bullet$]
	\item Répondre au problème \ref{probleme_sev_particulier_nreourbgea} dans le cas de l'approximation diophantienne classique, revient à déterminer l'exposant d'irrationalité d'un irrationnel donné. Ce problème est encore largement ouvert. À titre d'exemple, la recherche de l'exposant d'irrationalité de $\pi$ est toujours en cours (voir figure \ref{histoire_pi_abjpibmefpivn}), et de nombreux progrès sont faits régulièrement à ce sujet ; on conjecture que $\mu(\pi)=2$.
	
	\item Répondre au problème \ref{probleme_principal_gnzroign} dans le cas de l'approximation diophantienne classique, revient à déterminer le meilleur exposant d'approximation qui soit valide pour tous les irrationnels. Le théorème de Dirichlet (voir théorème \ref{approx_simul_dirichlet_aemofibbv}) répond à cette question : $2$ est le meilleur exposant d'approximation qui fonctionne simultanément pour tous les irrationnels, c'est-à-dire qu'on a $\mu(\xi)\ge 2$ pour tout $\xi\in\R\setminus\Q$ avec égalité pour au moins un $\xi$.
	
	\item Le problème \ref{probleme_spectre_oaerinemoivn} dans le cas de l'approximation diophantienne classique revient à déterminer l'ensemble $\{\mu(\xi),\ \xi\in\R\setminus\Q\}$. Grâce à la théorie des fractions continues (voir par exemple \cite{jarnik29}), on peut montrer que $\{\mu(\xi),\ \xi\in\R\setminus\Q\}=[2,+\infty]$.
\end{enumerate}

\begin{figure}[H]
\begin{center}
\begin{tabular}{c|c|c}
1953 & K. Mahler \cite{mahler53} & $\mu(\pi)\le 42$ \\
\hline
1974 & M. Mignotte \cite{mignotte74} & $\mu(\pi)\le 20.6$ \\
\hline 
1982 & G. V. Chudnovsky \cite{chudnovsky82} & $\mu(\pi)\le 18.9$ \\
\hline
1993 & M. Hata \cite{hata93} & $\mu(\pi)\le 13.398$ \\
\hline
2008 & V. K. Salikhov \cite{salikhov08} & $\mu(\pi)\le 7.6063$ \\
\hline
2020 & W. Zudilin et D. Zeilberger \cite{zeilberger20} & $\mu(\pi)\le 7.10321$
\end{tabular}
\end{center}
\caption{Majorations successives de l'exposant d'irrationalité de $\pi$}
\label{histoire_pi_abjpibmefpivn}
\end{figure}

Enfin, on peut conclure cette généralisation en remarquant que le problème de l'approximation d'irrationnels par des rationnels est contenu dans la généralisation donnée avec des sous-espaces vectoriels. En effet, approcher $\xi\in\R\setminus\Q$ revient dans le langage de l'approximation rationnelle des sous-espaces vectoriels à approcher 

\[A_\xi=\vect \left(\begin{pmatrix}1\\\xi\end{pmatrix}\right)\in\mathfrak I_2(1,1)_1\]

par des sous-espaces vectoriels rationnels $B\in\mathfrak R_2(1)$. On a 

\[\muexpA 2{A_\xi}11=\mu(\xi).\]

\section{Présentation des résultats en direction du problème \ref{probleme_principal_gnzroign}}\label{section_presentation_resultats_aboininbv}

La réponse idéale au problème \ref{probleme_principal_gnzroign} serait d'avoir une formule donnant explicitement la valeur de la quantité $\muexp ndej$ en fonction de $(n,d,e,j)$. Cependant cet objectif est encore loin d'être atteint. Cette section présente les différents résultats connus à ce jour apportant des réponses partielles à ce problème. Dans la sous-section \ref{sssection_Schmidt_originaux_moaibofdbv}, on énonce les résultats initiaux obtenus par Schmidt en 1967, et dans la sous-section \ref{cas_dim_1_ezrogbeorgb} on détaille le cas $\min(d,e)=1$ à travers les résultats de M. Laurent en 2009. Enfin, on présente dans la sous-section \ref{sssection_moshchevitin_gaoirhgpihvd} les résultats récents obtenus par \hbox{N. Moshchevitin} en 2020, et dans la sous-section \ref{sssection_de_Saxce_amoinbefoinvnie} ceux obtenus par \hbox{N. de Saxcé}, aussi en 2020.

\subsection{Les résultats originaux de Schmidt}\label{sssection_Schmidt_originaux_moaibofdbv}

Exposons ici les résultats généraux obtenus par Schmidt sur $\muexp ndej$ dans \cite{schmidt67}. \\

Dans toute la suite, on fixe $n\ge 2$ et $d,e\in\{1,\ldots,n-1\}$ tels que $d+e\le n$. On note $t=\min(d,e)$. \\

Le théorème 12 page 459 de l'article \cite{schmidt67} fournit une première minoration.

\begin{theoreme}[Schmidt, 1967]\label{premiere_borne_Schmidt_hgoembeer}

On a pour tout $j\in\{1,\ldots,t\}$,

\[\muexp ndej\ge\frac{d(n-j)}{j(n-d)(n-e)},\]

et même mieux dans le cas $j=1$ :

\[\muexp nde1\ge\frac{n(n-1)}{(n-d)(n-e)}.\]

\end{theoreme}

Schmidt améliore ce résultat sous une condition dans le théorème 13 page 460 de \cite{schmidt67} :

\begin{theoreme}[Schmidt, 1967]\label{deuxieme_borne_Schmidt_reuhfgezliuhl}

Soit $j\in\{1,\ldots,t\}$. Si

\begin{equation}\label{condition_du_th_13_Schmdit_oizhfgori}
j+n-t\ge j(j+n-d-e),
\vspace{2.6mm}
\end{equation}

alors

\[\muexp ndej\ge \frac{j+n-t}{j(j+n-d-e)}.\]

\end{theoreme}

Ce théorème est surtout utile lorsque $d+e$ est grand, ce qui permet à la condition \eqref{condition_du_th_13_Schmdit_oizhfgori} d'être vérifiée. En effet, si $d+e$ est petit devant $n$, seul le cas $j=1$ permet de satisfaire l'inégalité \eqref{condition_du_th_13_Schmdit_oizhfgori}. \\

Le théorème 15 page 462 de \cite{schmidt67} donne une réponse exacte au problème \ref{probleme_principal_gnzroign} dans certains cas particuliers :

\begin{theoreme}[Schmidt, 1967]\label{th_Schmidt_minoration_exp_oeurhgosvdn}

Si

\begin{equation}\label{condition_du_th_15_Schmdit_rgobaelrg}
n\ge t(t+n-d-e),
\vspace{2.6mm}
\end{equation}

alors

\[\muexp ndet= \frac{n}{t(t+n-d-e)}.\]

\end{theoreme}

On peut remarquer que la condition \eqref{condition_du_th_15_Schmdit_rgobaelrg} est toujours vérifiée dans le cas où $t=\min(d,e)=1$.

Ainsi, Schmidt a donné une réponse complète au problème \ref{probleme_principal_gnzroign} dans le cas où on cherche :
\begin{enumerate}[$\bullet$]
	\item soit à approcher une droite par des sous-espaces rationnels,
	\item soit à approcher un sous-espace par des droites rationnelles.
\end{enumerate}

Ce cas a continué à être étudié depuis, et une présentation en est faite dans la sous-section \ref{cas_dim_1_ezrogbeorgb} ci-dessous. \\

Schmidt fournit aussi page 465 de \cite{schmidt67}, sans aller jusqu'à l'optimalité, une majoration dans le cas général :

\begin{theoreme}[Schmidt, 1967]\label{borne_sup_de_Schmidt_roimghnev}

Soit $j\in\{1,\ldots,t\}$. On a

\[\muexp ndej\le \frac 1j\left\lceil\frac{e(n-e)+1}{n+1-d-e}\right\rceil.\]

\end{theoreme}

Ces différentes bornes permettent d'obtenir des encadrements de $\muexp ndej$ en fonction de $(n,d,e,j)$. À titre d'exemple et d'illustration, les différents encadrements connus de $\muexp ndej$ pour $n\in\{2,\ldots,6\}$ sont exposés dans la section \ref{section_exemples_pour_n_inferieur_6_maobuefobv} ci-dessous. \\ 

Le premier cas encore ouvert est le cas 

\[(n,d,e,j)=(4,2,2,1).\]

Les résultats de Schmidt ne donnent que l'encadrement

\[3\le \muexp 4221\le 5,\]

qui a été amélioré par Moshchevitin -- voir sous-section \ref{sssection_moshchevitin_gaoirhgpihvd} -- en

\[3\le \muexp 4221\le 4\]

et sera l'objet du théorème \ref{cas_total_4221_egorihgoegn}.

\subsection{Le cas $\min(d,e)=1$}\label{cas_dim_1_ezrogbeorgb}

A. Khintchine et A. V. Groshev (voir \cite{groshev38}) ont travaillé sur l'approximation d'un sous-espace vectoriel par une droite, mais d'après l'article \cite{laurent09} de Laurent, c'est surtout Schmidt qui fut le premier à étudier les exposants $\muexpA nLe1$ où $L$ est une droite rationnelle de $\R^n$. Schmidt adopte un point de vue plus géométrique que celui de Khintchine et Groshev. Grâce au théorème \ref{th_Schmidt_minoration_exp_oeurhgosvdn}, Schmidt a complètement résolu le problème consistant à déterminer $\muexp n1e1$ et $\muexp nd11$, \emph{i.e.} à approcher des droites ou à approcher par des droites. Cependant d'autres résultats furent apportés sur ces exposants, comme celui du théorème 2 page 3 de l'article \cite{laurent09} de Laurent, énoncé ci-dessous. 

\begin{theoreme}[Laurent, 2009]\label{MLaurent_transfert_eorighovfidn}

Soit $A$ une droite $(n-1,1)$-irrationnelle de $\R^n$. Alors on a un résultat de Going-up :

\[\forall e\in\{1,\ldots,n-2\},\quad \muexpA nA{e+1}1\ge \frac{(n-e)\muexpA nAe1}{n-e-1},\]

et un résultat de Going-down : 

\[\forall e\in\{2,\ldots,n-1\},\quad \muexpA nA{e-1}1\ge \frac {e\muexpA nAe1}{\muexpA nAe1 +e-1}.\]

\end{theoreme}

La première partie de ce théorème améliore le théorème 11 du Going-up de Schmidt (\cite{schmidt67} page 457) dans le cas particulier dans lequel se place Laurent. Dans la seconde, il retrouve la même minoration que le théorème \ref{goingdown_ogihoeihfgoniv}.

Laurent déduit du théorème \ref{MLaurent_transfert_eorighovfidn} le corollaire suivant (\cite{laurent09} page 4), qui a un triple intérêt :
\begin{enumerate}[$\bullet$]
	\item il retrouve la minoration de Schmidt du théorème \ref{deuxieme_borne_Schmidt_reuhfgezliuhl} dans le cas $d=1$,
	\item il donne l'exposant générique,
	\item il donne le spectre de l'exposant dans le cas particulier où $A$ est une droite \hbox{$(n-1,1)$}-irrationnelle, répondant au problème \ref{probleme_spectre_oaerinemoivn} dans ce cas.
\end{enumerate}

\begin{corollaire}[Laurent]

Soit $A$ une droite $(n-1,1)$-irrationnelle de $\R^n$. On a

\[\forall e\in\{1,\ldots,n-1\},\quad \muexpA {n}Ae1\ge \frac{n}{n-e}\]

avec égalité pour presque tout $A$ par rapport à la mesure de Lebesgue sur $\R^n$.

De plus, l'ensemble $\muexpA {n}{\mathfrak I_n(1,n-1)_1}e1$ est entièrement déterminé : pour tout $e\in\{1,\ldots,n-1\}$, on a

\[\Big\{\muexpA nAe1,\ A\in\mathfrak I_n(1,n-1)_1\Big\}=\left[\frac{n}{n-e},+\infty\right].\]

\end{corollaire}

Cet article \cite{laurent09} soulève aussi le problème de déterminer l'ensemble 

\[\Big\{(\muexpA {n}{A}11,\ldots,\muexpA {n}{A}{n-1}1),\ A\in\mathfrak I_n(1,n-1)_1\Big\},\]

problème qui a connu des avancées grâce à A. Marnat en 2018 (voir \cite{marnat18}). \\

Enfin, pour approfondir le cas $\min(d,e)=1$, on peut consulter à ce sujet \cite{bugeaud05}, \cite{bugeaud05b} ou \cite{bugeaud10}.

\subsection{Les résultats de Moshchevitin}\label{sssection_moshchevitin_gaoirhgpihvd}

Schmidt a montré que $\muexp 4221\in[3,5]$. Ce résultat a été amélioré en 2020 : Moshchevitin montre dans \cite{moshchevitin20} le théorème \ref{moshchevitin2020_ehgoierhf} ci-dessous. \\

Ici, Moshchevitin construit une mesure sur la grassmannienne $\Gr_{2,4}$ en la plongeant dans $\R^4$ grâce à un système de six cartes affines, puis en tirant en arrière la mesure de Lebesgue sur $\R^4$ ; c'est ce qui lui permet de parler de \emph{presque tout} sous-espace vectoriel.

\begin{theoreme}[Moshchevitin, 2020]\label{moshchevitin2020_ehgoierhf}

Soit $\omega\colon\R^+\to\R^+$ une application décroissante telle que

\[\sum_{j=1}^\infty j\omega(\sqrt j)<\infty.\]

Pour presque tout $A\in\mathfrak I_4(2,2)_1$, il existe une constante $c>0$ telle que

\[\forall B\in\mathfrak R_4(2),\quad \psi_1(A,B)\ge c\,\omega(H(B)).\]

\end{theoreme}

Du théorème \ref{moshchevitin2020_ehgoierhf} (appliqué avec $\omega(t)=t^{-4-\epsilon}$ pour $\epsilon>0$) découle notamment l'existence d'un sous-espace $A\in\mathfrak I_4(2,2)_1$ approché à l'exposant au plus $4$, d'où le corollaire suivant :

\begin{corollaire}[Moshchevitin, 2020]\label{coro_premier_th_Mosh_aomreigbovb}

On a

\[\muexp 4221\le 4,\]

et même pour presque tout $A\in\mathfrak I_4(2,2)_1$, 

\[\muexpA 4A21\le 4.\]

\end{corollaire}

Ainsi, Moshchevitin montre que l'exposant générique relatif au premier angle pour deux plans dans $\R^4$ se situe dans $[3,4]$. Il énonce de plus une généralisation du théorème \ref{moshchevitin2020_ehgoierhf}, en affirmant que la preuve est similaire à celle du théorème \ref{moshchevitin2020_ehgoierhf}.

\begin{theoreme}[Moshchevitin, 2020]\label{th_Moshchevitin_2020_general_ergoighvoifvn}

Soient $n=2d$, et $\omega\colon\R^+\to\R^+$ une application décroissante telle que

\[\sum_{j=1}^\infty j^{d-1}\omega(\sqrt j)<\infty.\]

Pour presque tout $A\in\mathfrak I_n(d,d)_1$, il existe une constante $c>0$ telle que

\[\forall B\in\mathfrak R_n(d),\quad \psi_1(A,B)\ge c\,\omega(H(B)).\]

\end{theoreme}

Comme pour le corollaire \ref{coro_premier_th_Mosh_aomreigbovb}, le théorème \ref{th_Moshchevitin_2020_general_ergoighvoifvn} donne un sous-espace $A\in\mathfrak I_n(d,d)_1$ approché à l'exposant au plus $2d$, ce qui implique le corollaire suivant :

\begin{corollaire}[Moshchevitin, 2020]\label{cor_Moshchevitin2020_2d_gezroguheoufb}

Soit $n=2d$, on a

\[\muexp ndd1\le 2d,\]

et même pour presque tout $A\in\mathfrak I_n(d,d)_1$, 

\[\muexpA nAd1\le 2d.\]

\end{corollaire}

Ceci améliore considérablement la majoration de Schmidt, car le théorème \ref{borne_sup_de_Schmidt_roimghnev} avec $n=2d$ donne

\[\muexp ndd1\le d^2+1.\]

À titre d'exemple, le corollaire \ref{cor_Moshchevitin2020_2d_gezroguheoufb} montre que

\[\muexp 6331\le 6,\]

ce qui améliore la majoration de Schmidt $\muexp 6331\le 10$.

\subsection{Les résultats de Saxcé}\label{sssection_de_Saxce_amoinbefoinvnie}

Soient $n\ge 2$ et $d\in\{1,\ldots,\lfloor n/2\rfloor\}$. D'après le théorème \ref{premiere_borne_Schmidt_hgoembeer}, on a

\begin{equation}\label{min_Schmidt_de_Saxce_amoerubobuivc}
\muexp nddd\ge \frac 1{n-d}.
\vspace{2.6mm}
\end{equation}

Dans \cite{saxce20}, Saxcé donne dans le théorème 9.3.2 page 105 l'exposant générique relatif au $d$-ième angle pour deux sous-espaces vectoriels de dimension $d$ de $\R^n$. Cet exposant générique est supérieur à la minoration \eqref{min_Schmidt_de_Saxce_amoerubobuivc} de Schmidt. La classe de mesure que Saxcé utilise est la classe de la mesure de Lebesgue sur les points réels de la grassmannienne $\Gr(d,n)$, c'est ce qui lui permet de parler de \emph{presque tout} $A\in\Gr(d,n)$.

\begin{theoreme}[Saxcé, 2020]\label{th_generique_de_Saxce_oaimrnfosinv}

Soient $A$ un sous-espace vectoriel de $\R^n$ de dimension $d$ et $\omega\colon\R^+\to\R^+$ une fonction décroissante. On considère pour $B\in\mathfrak R_n(d)$ l'inégalité 

\begin{equation}\label{th_principal_Saxce_baineroeinbv}
\psi_d(A,B)\le \omega(H(B)) H(B)^{-n/(d(n-d))}.
\vspace{2.6mm}
\end{equation}

\begin{enumerate}[$\bullet$]
	\item Si $\int_1^\infty t^{-1}\omega(t)^{d(n-d)}\d t=+\infty$, alors pour presque tout $A\in\Gr(d,n)$ l'inégalité \eqref{th_principal_Saxce_baineroeinbv} admet une infinité de solutions $B\in\mathfrak R_n(d)$.
	
	\item Si $\int_1^\infty t^{-1}\omega(t)^{d(n-d)}\d t<\infty$, alors pour presque tout $A\in\Gr(d,n)$ l'inégalité \eqref{th_principal_Saxce_baineroeinbv} n'admet qu'un nombre fini de solutions $B\in\mathfrak R_n(d)$.
\end{enumerate}

\end{theoreme}

Comme $\mathfrak I_n(d,d)_d=X(\R)\setminus X(\Q)$ avec $X=\Gr(d,n)$, on déduit du théorème \ref{th_generique_de_Saxce_oaimrnfosinv} le corollaire suivant.

\begin{corollaire}[Saxcé, 2020]\label{cor_de_Saxce_aoeiboifbvoiefbvoi}

Pour presque tout $A\in\mathfrak I_n(d,d)_d$, on a

\[\muexpA nAdd=\frac n{d(n-d)}.\]

\end{corollaire}

\begin{corollaire}[Saxcé, 2020]\label{cor_maj_de_Saxce_amovrvrhbeofin}

On a 

\[\muexp nddd\le \frac n{d(n-d)}.\]

\end{corollaire}

Saxcé estime possible (\cite{saxce20} page 106) qu'on ait égalité dans ce corollaire. \\

Enfin, Saxcé s'intéresse au cas des sous-espaces définis par des équations à coefficients algébriques.

\begin{theoreme}[Saxcé, 2020]

Soit $A\in\mathfrak I_n(d,d)_d$ un sous-espace vectoriel défini par des équations à coefficients algébriques. Alors

\begin{equation}\label{minoration_de_Saxce_algebrique_abemoribiv}
\muexpA nAdd\ge \frac n{d(n-d)}.
\vspace{2.6mm}
\end{equation}

\end{theoreme}

\begin{remarque}

Dans \cite{saxce20} page 107, Saxcé donne une condition nécessaire et suffisante pour que la minoration \eqref{minoration_de_Saxce_algebrique_abemoribiv} soit une égalité, ce qui apporte une réponse partielle au problème \ref{probleme_sev_particulier_nreourbgea}.

\end{remarque}

\section{Plan de la thèse et résultats principaux}\label{section_plan_de_la_these_abinrebinv}

On termine cette introduction par un bref résumé du contenu de chaque chapitre. On énonce aussi les théorèmes principaux démontrés dans cette thèse.

\subsection*{Chapitre \ref{chap_outils_oaeibov}}

De nombreux outils qui seront utiles pour la suite sont regroupés dans ce chapitre. On commence par développer dans la sous-section \ref{coord_Pluck_teoghzleufve} ce que sont les coordonnées de Plücker, qui permettent de travailler avec la hauteur d'un sous-espace rationnel. \\

On donne ensuite plusieurs outils sur la hauteur dans les sous-sections \ref{ss_section_det_gene_omaiosbv} et \ref{ss_section_hauteur_transfo_lineaire_aboimfbevbv}, et sur la proximité dans la section \ref{section_proximite_aobomebvob}. \\

Enfin, on énonce deux théorèmes d'approximation simultanée dans la section \ref{section_resultats_approx_simult_aoeibhaoimnvao}, et les théorèmes du Going-up et du Going-down dans la section \ref{section_goingupdown_aogieovb}.

\subsection*{Chapitre \ref{chap_cas_particuliers_oamofbvoaube}}

Différents cas particuliers sont étudiés dans ce chapitre. On commence par énoncer quelques outils dans la section \ref{section_outils_en_toute_dimension_amoimaoihogierbobi}, puis on traite deux cas particuliers dans $\R^4$ et $\R^5$, respectivement dans les sections \ref{section_R4_maoefvbn} et \ref{section_R5_gaomfnbamobo}. \\

Le cas le plus simple encore ouvert est la détermination de $\muexp 4221$, quantité encadrée entre $3$ et $5$ par Schmidt en 1967, puis entre $3$ et $4$ par Moshchevitin en 2020. Le théorème \ref{cas_total_4221_egorihgoegn} la détermine exactement, grâce à un plan explicite de $\mathfrak I_4(2,2)_1$ particulièrement mal approché par les plans rationnels. C'est ainsi qu'on obtient le :

\begin{theoremenonum}\textbf{\ref{cas_total_4221_egorihgoegn}}
On a

\[\muexp 4221=3.\]

\end{theoremenonum}

Par ailleurs, en considérant un élément explicite de $\mathfrak I_5(3,2)_1$ mal approché par les plans rationnels, on obtient une amélioration de l'encadrement connu de $\muexp 5321$ dû à Schmidt. Cette quantité est encadrée par $4$ et $7$ d'après les théorèmes \ref{deuxieme_borne_Schmidt_reuhfgezliuhl} et \ref{borne_sup_de_Schmidt_roimghnev} respectivement.

\begin{theoremenonum}\textbf{\ref{amelioration_cas_R5_eoufbemovb}}
On a

\[\muexp 5321\le6.\]

\end{theoremenonum}

Dans la section \ref{section_app_Moshchevitin_aoimbnoibav}, on déduit un nouveau majorant de $\muexp nd{d-1}1$ comme corollaire du théorème \ref{th_Moshchevitin_2020_general_ergoighvoifvn} : pour $d\ge 2$ un entier et $n=2d$, on a

\[\muexp nd{d-1}1\le \frac{2d^2}{d+1}.\]

\subsection*{Chapitre \ref{chap_somme_directes_vaoinvoin}}

Des résultats de reconstruction sur le comportement de la proximité et de la hauteur par somme directe sont obtenus dans la sous-section \ref{ss_section_resultats_reconstruction_agmoibeamoiv}. Grâce à ces résultats de reconstruction, on améliore asymptotiquement le minorant connu de $\muexp ndej$ dans le théorème \ref{premiereborneobtenue_amorimeovbn}. Dans le cas particulier où $d=e=j$, on obtient alors le corollaire suivant.

\begin{corollairenonum}\textbf{\ref{corollaire_aeroibfnoiefbivo}}
Soient $n\ge 4$ et $d\le n/2$ un entier. On a

\[\muexp nddd\ge\frac{2dn-d^2+d+2}{2d^2n-d^3+d^2}.\]

\end{corollairenonum}

En combinant ce corollaire avec le corollaire \ref{cor_maj_de_Saxce_amovrvrhbeofin} de Saxcé, on obtient pour $n\ge 2d$ l'encadrement

\[\frac{2dn-d^2+d+2}{2d^2n-d^3+d^2}\le \muexp nddd\le \frac n{d(n-d)},\]

ce qui permet d'en déduire le résultat suivant.

\begin{corollairenonum}\textbf{\ref{cor_limite_avec_maj_desaxce_amoberubfvs}}
On a 

\[\lim_{n\to+\infty} \muexp nddd = \frac 1d.\]

\end{corollairenonum}

\subsection*{Chapitre \ref{chap_inclusion_sev_rationnel_amoriegmoingoa}}

Dans ce chapitre, on démontre le théorème \ref{th_inclusion_sev_rationnel_apeivpinpiaenv}, qui permet de baisser la dimension de l'espace ambiant lorsque le sous-espace vectoriel approché est inclus dans un sous-espace rationnel :

\begin{theoremenonum}\textbf{\ref{th_inclusion_sev_rationnel_apeivpinpiaenv}}
Soient $n\ge 2$ et $k\in\{2,\ldots,n\}$. Soient $d,e\in\{1,\ldots,k-1\}$ tels que $d+e\le k$, et $j\in\{1,\ldots,\min(d,e)\}$. Soit $A$ un sous-espace vectoriel de $\R^n$ de dimension $d$ tel qu'il existe un sous-espace vectoriel rationnel $F\in\mathfrak R_n(k)$ vérifiant $A\subset F$.

Notons $\varphi$ un isomorphisme rationnel de $F$ dans $\R^k$ et $\tilde A=\varphi(A)$, qui est un sous-espace vectoriel de dimension $d$ de $\R^k$.

Supposons que pour tout sous-espace rationnel $B'$ de $F$ de dimension $e$, on a 

\begin{equation}\label{condition_dirr_sur_A_ameirnmoafvonsd}
\dim (A\cap B')<j.
\vspace{2.6mm}
\end{equation}

Alors $A\in\mathfrak I_n(d,e)_j$, $\tilde A\in\mathfrak I_k(d,e)_j$ et 

\[\muexpA nAej=\muexpA k{\tilde A}ej.\]

\end{theoremenonum}

On déduit ensuite de ce théorème l'amélioration de quelques majorations connues de $\muexp ndej$ dans la section \ref{section_app_du_th_inclusiion_sev_rationnel_vaoihbmoifv}. On obtient que 

\[\muexp 5221\le 3,\]

ce qui améliore le majorant connu $\muexp 5221\le 4$ donné par le théorème \ref{borne_sup_de_Schmidt_roimghnev}. On déduit aussi que pour tout $d\ge 2$ et pour tout $n\ge 2d$, on a

\[\muexp ndd1\le 2d\]

et 

\[\muexp nd{d-1}1\le \frac{2d^2}{d+1},\]

ce qui améliore des majorations dans des cas où $n$ est proche de $2d$ : $\muexp 7331\le 7$ devient $\muexp 7331\le 6$ et $\muexp 9431\le 7$ devient $\muexp 9431\le 32/5=6.4$ par exemple.

\subsection*{Chapitre \ref{chapitre_spectre_aoimhremvoin}}

L'objectif de ce chapitre est de montrer le résultat suivant sur le spectre de l'exposant relatif au $\ell$-ième angle entre deux sous-espaces vectoriels de dimension $\ell$ :

\begin{theoremenonum}\textbf{\ref{theoreme_spectre_amoeribnefaomnv}}
Soient $n\ge 2$ et $\ell\in\{1,\ldots,\lfloor n/2\rfloor\}$.

Alors 

\[\left[1+\frac 1{2\ell}+\sqrt{1+\frac{1}{4\ell^2}},+\infty\right]\subset\Big\{\muexpA nA\ell\ell,\ A\in\mathfrak I_n(\ell,\ell)_\ell\Big\}.\]

\end{theoremenonum}

Pour ce faire, on construit explicitement un sous-espace vectoriel d'exposant prescrit dans la section \ref{section_resultats_sev_exposant_donne_aefoubv}.

\section{Les résultats connus jusqu'en dimension $6$}\label{section_exemples_pour_n_inferieur_6_maobuefobv}

On expose ici les premières valeurs de $\muexp ndej$ pour $n\in\{2,\ldots,6\}$, ce qui permet d'avoir une idée de l'ordre de grandeur de $\muexp ndej$. Ces valeurs sont présentées dans les tableaux \ref{casdeR2_regoimehoh}, \ref{casdeR3_regoimehoh}, \ref{casdeR4_regoimehoh}, \ref{casdeR5_regoimehoh} et \ref{casdeR6_regoimehoh} ci-dessous. On se place dans les cas où $d+e\le n$, avec $d=\dim A$ et $e=\dim B$ ; les étoiles correspondent aux cas $d+e>n$. Pour $n\in\{2,\ldots,6\}$, les encadrements connus sont principalement dus à Schmidt (théorèmes \ref{premiere_borne_Schmidt_hgoembeer}, \ref{deuxieme_borne_Schmidt_reuhfgezliuhl}, \ref{th_Schmidt_minoration_exp_oeurhgosvdn} et \ref{borne_sup_de_Schmidt_roimghnev}). Les majorations $\muexp 4221\le 4$ et $\muexp 6331\le 6$ sont dues à Moshchevitin (corollaires \ref{coro_premier_th_Mosh_aomreigbovb} et \ref{cor_Moshchevitin2020_2d_gezroguheoufb} respectivement), le théorème \ref{borne_sup_de_Schmidt_roimghnev} de Schmidt donnant seulement $\muexp 4221\le 5$ et $\muexp 6331\le 10$. Les majorations $\muexp 5222\le 5/6$, $\muexp 6222\le 3/4$ et $\muexp 6333\le 2/3$ sont dues à Saxcé (corollaire \ref{cor_maj_de_Saxce_amovrvrhbeofin}), le théorème \ref{borne_sup_de_Schmidt_roimghnev} donnant seulement $\muexp 5222\le 2$, $\muexp 6222\le 3/2$ et $\muexp 6333\le 10/3$. \\

Les lignes en gras correspondent aux améliorations présentées dans cette thèse des encadrement connus jusqu'ici de $\muexp ndej$. \\

Remarquons que Schmidt émet l'hypothèse dans \cite{schmidt67} page 471 que $\muexp ndej$ est une fonction symétrique en $d$ et $e$, sans aller jusqu'à conjecturer ce résultat. 

\begin{figure}[H]
\begin{center}
\begin{tabular}{c||c}
$\R^2$ &$\dim A=1$ \\
\hline\hline
$\dim B=1$ &$\muexp 2111=2$
\end{tabular}
\end{center}
\caption{Le cas de $\R^2$}
\label{casdeR2_regoimehoh}
\end{figure}

\begin{figure}[H]
\begin{center}
\begin{tabular}{c||c|c}
$\R^3$ &$\dim A=1$&$\dim A=2$ \\
\hline\hline
$\dim B=1$&$\muexp 3111=\frac 32$&$\muexp 3211=3$ \\
\hline
$\dim B=2$&$\muexp 3121=3$&$\star$
\end{tabular}
\end{center}
\caption{Le cas de $\R^3$}
\label{casdeR3_regoimehoh}
\end{figure}

\begin{figure}[H]
\begin{center}
\begin{tabular}{c||c|c|c}
$\R^4$ &$\dim A=1$&$\dim A=2$&$\dim A=3$ \\
\hline\hline
$\dim B=1$&$\muexp 4111=\frac 43$&$\muexp 4211=2$&$\muexp 4311=4$ \\
\hline
\multirow{3}{*}{$\dim B=2$}&\multirow{3}{*}{$\muexp 4121=2$}&$3\le\muexp 4221 \le 4$ \\
&&$\mathbf{\muexp 4221=3}$&$\star$ \\
&&$\muexp 4222=1$ \\
\hline
$\dim B=3$&$\muexp 4131=4$&$\star$&$\star$
\end{tabular}
\end{center}
\caption{Le cas de $\R^4$}
\label{casdeR4_regoimehoh}
\end{figure}

Le résultat $\muexp 4221=3$ est celui du théorème \ref{cas_total_4221_egorihgoegn}.

\begin{figure}[H]
\begin{center}
\begin{tabular}{c||c|c|c|c}
$\R^5$ &$\dim A=1$&$\dim A=2$&$\dim A=3$&$\dim A=4$ \\
\hline\hline
$\dim B=1$&$\muexp 5111=\frac 54$&$\muexp 5211=\frac 53$&$\muexp 5311=\frac 52$&$\muexp 5411=5$ \\
\hline
\multirow{4}{*}{$\dim B=2$}&\multirow{4}{*}{$\muexp 5121=\frac 53$}&$\frac{20}9\le\muexp 5221 \le 4$&\multirow{1.59}{*}{$4\le\muexp 5321\le 7$}&\multirow{4}{*}{$\star$} \\
&&$\mathbf{\muexp 5221\le 3}$&\multirow{1.80}{*}{$\mathbf{\muexp 5321\le 6}$} \\
&&$\frac 13\le\muexp 5222\le 5/6$&\multirow{2.16}{*}{$\muexp 5322=\frac 54$} \\
&&$\mathbf{\frac 59\le\muexp 5222}$& \\
\hline
\multirow{2}{*}{$\dim B=3$}&\multirow{2}{*}{$\muexp 5131=\frac 52$}&$4\le\muexp 5231\le 7$&\multirow{2}{*}{$\star$}&\multirow{2}{*}{$\star$} \\
&&$\muexp 5232=\frac 54$& \\
\hline
$\dim B=4$&$\muexp 5141=5$&$\star$&$\star$&$\star$
\end{tabular}
\end{center}
\caption{Le cas de $\R^5$}
\label{casdeR5_regoimehoh}
\end{figure}

La majoration $\muexp 5221\le 3$ est celle de la proposition \ref{amelioration_R5_grace_inclusion_R4_aoibvoeubv} (proposition qui combine la proposition \ref{proposition_R4_mu_Axi_zgiozovdib} et le théorème \ref{th_inclusion_sev_rationnel_apeivpinpiaenv}), la majoration $\muexp 5321\le 6$ provient du théorème \ref{amelioration_cas_R5_eoufbemovb} et la minoration $\muexp 5222\ge 5/9$ se déduit du théorème \ref{premiereborneobtenue_amorimeovbn}.

\begin{figure}[H]\footnotesize
\begin{center}
\begin{tabular}{c||c|c|c|c|c}
$\R^6$ &$\dim A=1$&$\dim A=2$&$\dim A=3$&$\dim A=4$&$\dim A=5$ \\
\hline\hline
$\dim B=1$&$\muexp 6111=\frac 65$&$\muexp 6211=\frac 32$&$\muexp 6311=2$&$\muexp 6411=3$&$\muexp 6511=6$ \\
\hline
\multirow{3}{*}{$\dim B=2$}&\multirow{3}{*}{$\muexp 6121=\frac 32$}&$\frac{15}8\le\muexp 6221 \le 3$&$\frac 52\le\muexp 6321\le 5$&\multirow{1.7}{*}{$5\le\muexp 6421\le 9$}&\multirow{3}{*}{$\star$} \\
&&$\frac 14\le\muexp 6222\le \frac 34$&$\mathbf{\muexp 6321\le \frac 92}$&\multirow{1.9}{*}{$\muexp 6422=\frac 32$} \\
&&$\mathbf{\frac 6{11}\le\muexp 6222}$&$\muexp 6322=1$& \\
\hline
\multirow{4}{*}{$\dim B=3$}&\multirow{4}{*}{$\muexp 6131=2$}&\multirow{2.5}{*}{$\frac 52\le\muexp 6231\le 5$}&$4\le\muexp 6331\le 6$&\multirow{4}{*}{$\star$}&\multirow{4}{*}{$\star$} \\
&&&$\frac 54\le\muexp 6332\le 5$& \\
&&\multirow{1.6}{*}{$\muexp 6232=1$}&$\frac 13\le\muexp 6333\le\frac {2}3$& \\
&&&$\mathbf{\frac{16}{45}\le\muexp 6333}$& \\
\hline
\multirow{2}{*}{$\dim B=4$}&\multirow{2}{*}{$\muexp 6141=3$}&$5\le\muexp 6241\le 9$&\multirow{2}{*}{$\star$}&\multirow{2}{*}{$\star$}&\multirow{2}{*}{$\star$} \\
&&$\muexp 6242=\frac 32$&& \\
\hline
$\dim B=5$&$\muexp 6151=6$&$\star$&$\star$&$\star$&$\star$
\end{tabular}
\end{center}
\caption{Le cas de $\R^6$}
\label{casdeR6_regoimehoh}
\end{figure}

Les minorations $\muexp 6222\ge 6/11$ et $\muexp 6333\ge 16/45$ sont des cas particuliers du théorème \ref{premiereborneobtenue_amorimeovbn}. La majoration $\muexp 6321\le 9/2$ se déduit du théorème \ref{theoreme_avec_moshchevitin_goingup_aeronfvn} qui découle des résultats de Moshchevitin (théorème \ref{th_Moshchevitin_2020_general_ergoighvoifvn}).

\section{Une conjecture sur $\muexp ndee$}\label{section_conjecture_oabeneboinv}

Commençons par mentionner le théorème 8 page 449 de \cite{schmidt67} :

\begin{theoreme}\label{theoreme_Schmidt_B_1_B_u_baemotibonfiv}

Soient $n\ge 2$ et $d\in\{1,\ldots,n-1\}$, posons $t=\min(d,n-d)$ ; soit $A\in\mathfrak I_n(d,t)_t$. Il existe des sous-espaces $B_1\subset \cdots\subset B_t$ de $\mathfrak R_n(1),\ldots,\mathfrak R_n(t)$ respectivement, tels que pour tout $j\in\{1,\ldots,t\}$, on a 

\[H(B_1)\psi_j(A,B_j)\le \frac{c}{H(B_j)^{d/(j(n-d))}}\]

où $c>0$ ne dépend uniquement de $n$, $d$ et $e$.

\end{theoreme}

On peut alors remarquer -- comme le fait Schmidt en conclusion de \cite{schmidt67} -- qu'un sous-espace mal approché par les droites rationnelles est très bien approché par des sous-espaces rationnels de dimension supérieure. Explicitons ce que donne l'approche qu'il propose. \\

Soient $n\ge 2$, $d,e\in\{1,\ldots,n-1\}$ des entiers tels que $e\le d\le n-e$, et $A\in\mathfrak I_n(d,e)_e$. Supposons que pour tout $\epsilon>0$, il n'existe qu'un nombre fini droites rationnelles $B$ telles que

\[\psi_1(A,B)\le \frac {1}{H(B)^{n/(n-d)+\epsilon}}.\]

Ainsi, pour tout $\epsilon>0$ il existe une constante $c_1>0$ dépendant uniquement de $A$, de $n$ et de $\epsilon$ telle que

\[\forall B\in\mathfrak R_n(1),\quad \psi_1(A,B)\ge \frac{c_1}{H(B)^{n/(n-d)+\epsilon}}\]

soit

\[\forall B\in\mathfrak R_n(1),\quad \frac 1{H(B)}\le c_2\psi_1(A,B)^{(n-d)/(n+\epsilon(n-d))}\]

avec $c_2>0$ ne dépendant que de $A$, de $n$ et de $\epsilon$. Le théorème \ref{theoreme_Schmidt_B_1_B_u_baemotibonfiv} donne alors une infinité de sous-espaces $B_1\subset \cdots\subset B_e$ de $\mathfrak R_n(1),\ldots,\mathfrak R_n(e)$ respectivement, tels que 

\[\psi_e(A,B_e)\le \frac 1{H(B_1)}\cdot\frac{c}{H(B_e)^{d/(e(n-d))}}\le \frac{c_3 \psi_1(A,B_1)^{(n-d)/(n+\epsilon(n-d))}}{H(B_e)^{d/(e(n-d))}}\]

avec $c_3>0$ ne dépendant que de $A$, de $n$, de $e$ et de $\epsilon$. Or d'après le lemme \ref{lemme_transfert_psi1_psiell_baoeoearv}, on a

\[\psi_1(A,B_1)\le \psi_e(A,B_e),\]

donc

\[\psi_e(A,B_e)\le c_4 H(B_e)^{-\frac{d}{e(n-d)}\cdot \frac{n+\epsilon(n-d)}{d+\epsilon(n-d)}}\]

avec $c_4>0$ ne dépendant que de $A$, de $n$, de $e$ et de $\epsilon$. Finalement, en faisant tendre $\epsilon$ vers $0$, on obtient

\[\muexpA nAee\ge \frac n{e(n-d)}.\]

En comparant ce résultat avec le corollaire \ref{cor_maj_de_Saxce_amovrvrhbeofin} qui donne $\muexp nddd\le n/(d(n-d))$, on émet la conjecture suivante qui rejoint, lorsque $d=e$, la possibilité évoquée par Saxcé (\cite{saxce20} page 106) qu'on ait égalité dans le corollaire \ref{cor_maj_de_Saxce_amovrvrhbeofin}. De plus, d'après le théorème \ref{th_Schmidt_minoration_exp_oeurhgosvdn}, cette conjecture est vraie si $n\ge e(n-d)$.

\begin{conjecture}\label{conjecture_rmoibfaobva}

Soient $n\ge 2$, $d,e\in\{1,\ldots,n-1\}$ des entiers tels que $e\le d\le n-e$. On a 

\[\muexp ndee=\frac n{e(n-d)}.\]

\end{conjecture}

On remarque aussi que le corollaire \ref{cor_de_Saxce_aoeiboifbvoiefbvoi} donne que $\muexpA nAdd=n/(d(n-d))$ pour presque tout $A\in\mathfrak I_n(d,d)_d$. \\

Enfin, Schmidt émet l'hypothèse que $\muexp ndej$ soit une fonction symétrique de $d$ et $e$ (\cite{schmidt67} page 471). Dans ce cas, la conjecture \ref{conjecture_rmoibfaobva} s'étend au cas $e>d$ et devient

\[\muexp nde{\min(d,e)}=\frac n{\min(d,e)(n-\max(d,e))}.\]

\chapter{Des outils}\label{chap_outils_oaeibov}

Dans ce chapitre, après avoir développé les coordonnées de Plücker dans la sous-section \ref{coord_Pluck_teoghzleufve}, qui permettront d'énoncer des résultats utiles pour travailler avec la hauteur d'un sous-espace rationnel, il s'agira de donner quelques résultats sur la hauteur (sous-sections \ref{ss_section_det_gene_omaiosbv} et \ref{ss_section_hauteur_transfo_lineaire_aboimfbevbv}) et sur la proximité (section \ref{section_proximite_aobomebvob}). \\

On énonce aussi quelques théorèmes connus, d'approximation simultanée et de Going-up/Going-down dans les sections \ref{section_resultats_approx_simult_aoeibhaoimnvao} et \ref{section_goingupdown_aogieovb}. 

\section{Outils liés à la hauteur}

\subsection{Les coordonnées de Plücker}\label{coord_Pluck_teoghzleufve}

Les coordonnées de Plücker sont au c\oe ur de la définition \ref{def_hauteur_coord_plucker_reuighgzmne} de la hauteur d'un sous-espace rationnel. Cette sous-section rappelle quelques résultats sur ces coordonnées, en s'appuyant sur le chapitre I du livre \cite{caldero15} de Caldero et Germoni, dans lequel tous les résultats énoncés ici sont démontrés. \\

Soient $m,n,r\in\N^*$ tels que $r\le\min(m,n)$ ; on note $\Lambda(r,n)$ l'ensemble des parties à $r$ éléments de $\{1,\ldots,n\}$. Pour $I\in\Lambda(r,m)$, $J\in\Lambda(r,n)$ et $A$ une matrice de taille $m\times n$, on note $\Delta_{I,J}(A)$ le mineur correspondant aux lignes de $A$ indexées par $I$, et aux colonnes de $A$ indexées par $J$.

Si $I$ et $J$ sont deux parties de $\{1,\ldots,n\}$, on définit $\iota(I,J)$ le nombre d'inversions de $I$ à $J$ par

\[\iota(I,J)= \card(\{(i,j)\in I\times J,\ i>j\}).\]

On note $\bar I$ le complémentaire de $I$ dans $\{1,\ldots,n\}$, et $\ell(I)=\iota(I,\bar I)$. \\

Une formule qui sera utile par la suite est la suivante. Elle généralise le développement par rapport à une ligne ou une colonne du déterminant d'une matrice.

\begin{proposition}[développement de Laplace]\label{devlaplace_aomrgihmoaeiv}

Soient $A\in\M_n(\R)$ et $J\in\Lambda(r,n)$. Alors

\[\det(A)=\sum_{I\in\Lambda(r,n)}(-1)^{\ell(I)+\ell(J)}\Delta_{I,J}(A)\Delta_{\bar I,\bar J}(A).\]

\end{proposition}

\begin{exemple}\label{exemple_Laplace_24_bariomnbofeiv}

Comme le cas de $A\in\M_4(\R)$ sera utilisé par la suite, illustrons le développement de Laplace en calculant le déterminant de $A$ par rapport à ses deux premières colonnes ($J=\{1,2\}$). On a 

\[\ell(\{1,2\})=\iota(\{1,2\},\{3,4\})=\card(\{(i,j)\in \{1,2\}\times \{3,4\},\ i>j\})=0,\]

et de même

\[\ell(\{1,3\})=1,\ \ell(\{1,4\})=2,\ \ell(\{2,3\})=2,\ \ell(\{2,4\})=3\text{ et }\ell(\{3,4\})=4,\]

donc
\begin{align*}
\det(A) =
	&\phantom{+} \Delta_{\{1,2\},\{1,2\}}(A)\Delta_{\{3,4\},\{3,4\}}(A)-\Delta_{\{1,3\},\{1,2\}}(A)\Delta_{\{2,4\},\{3,4\}}(A) \\
	&+\Delta_{\{1,4\},\{1,2\}}(A)\Delta_{\{2,3\},\{3,4\}}(A)+\Delta_{\{2,3\},\{1,2\}}(A)\Delta_{\{1,4\},\{3,4\}}(A) \\
	&-\Delta_{\{2,4\},\{1,2\}}(A)\Delta_{\{1,3\},\{3,4\}}(A)+\Delta_{\{3,4\},\{1,2\}}(A)\Delta_{\{1,2\},\{3,4\}}(A).
\end{align*}
\end{exemple}

Soient $A\in\M_{n,m}(\R)$ et $r\le \min(n,m)$. On utilise l'ordre lexicographique sur $\Lambda(r,n)$ et $\Lambda(r,m)$, ce qui permet d'identifier $\Lambda(r,n)$ à $\{1,\ldots,N\}$ et $\Lambda^r(\R^n)$ à $\R^N$, où $N=\binom nr$. On définit la \emph{$r$-ième puissance extérieure} de $A$, comme la matrice dont les coefficients sont les mineurs extraits de taille $r\times r$, \emph{i.e.}

\[\Lambda^r(A)=(\Delta_{I,J}(A))_{I\in\Lambda(r,n),J\in\Lambda(r,m)}.\]

Ces quelques notations permettent d'arriver au théorème du plongement de Plücker. L'idée du plongement de Plücker est de voir les sous-espaces vectoriels de $\R^n$ de dimension $d$, comme des points d'un espace projectif. Autrement dit, de voir tous les sous-espaces vectoriels de $\R^n$ comme des droites d'un espace plus grand. \\

On note $\P(\R^n)$ l'ensemble des droites vectorielles de $\R^n$ et $\Gr_{r,n}$ l'ensemble des sous-espaces vectoriels de $\R^n$ de dimension $r$, qui s'appelle la grassmannienne. 

\begin{theoreme}[plongement de Plücker]\label{plongementdepluck_egomhgoev}

Soient $n\ge 2$ un entier et \hbox{$r\in\{1,\ldots,n-1\}$}. Soient $F$ un sous-espace vectoriel de dimension $r$ de $\R^n$ et $A_F\in\M_{n,r}(\R)$ la matrice d'une base de $F$ dans la base canonique de $\R^n$.

Alors la droite engendrée par $\Lambda^r(A_F)$ ne dépend que de $F$. On note $[\Lambda^r(A_F)]$ le point correspondant dans l'espace projectif, et $N=\binom nr$. On peut donc définir une application $\psi_{r,n}$ par

\[\begin{matrix} \psi_{r,n}\colon & \Gr_{r,n} &\longhookrightarrow &\P(\R^N) \\ &F &\longmapsto &[\Lambda^r(A_F)].\end{matrix}\]

De plus, cette application -- appelée \emph{plongement de Plücker de $\Gr_{r,n}$} -- est injective.

\end{theoreme}

Arrive alors la définition clef de ce paragraphe.

\begin{definition}\label{defcoordgrassmann_reaoimghermoh}

Soit $F$ un sous-espace vectoriel de $\R^n$ de dimension $r$ ; on note $A_F$ la matrice d'une base de $F$ dans la base canonique. Les $\binom nr$ coordonnées du vecteur $\Lambda^r(A_F)$ sont appelées les \emph{coordonnées de Plücker} de $F$.

\end{definition}

Dans toute la suite, les coordonnées de Plücker (qui sont les mineurs de taille maximale de $A_F$) seront donc ordonnées par l'ordre lexicographique.

\begin{remarque}

Le théorème \ref{plongementdepluck_egomhgoev} assure que deux sous-espaces vectoriels différents auront des coordonnées de Plücker associées non proportionnelles, et que les coordonnées de Plücker d'un sous-espace donné sont uniques à multiplication par un scalaire près. Ce sont ces deux points cruciaux qui donnent tout leur intérêt aux coordonnées de Plücker.

\end{remarque}

Le théorème \ref{plongementdepluck_egomhgoev} construit un plongement $\psi_{r,n}$ de $\Gr_{r,n}$ dans $\P(\R^N)$, mais une question reste en suspens : quels sont les points de $\P(\R^N)$ qui ont un antécédent par $\psi_{r,n}$ ? Autrement dit, quels points de l'espace projectif sur $\R^N$ correspondent effectivement à un sous-espace vectoriel ? Cette question est résolue dans le théorème suivant :

\begin{theoreme}[relations de Plücker]\label{relation_de_Plucker_oerighbfduvb}

Soient $n\ge 2$ et $r\in\{1,\ldots,n-1\}$ ; posons $N=\binom nr$. Soit $\psi_{r,n}\colon\Gr_{r,n}\hookrightarrow\P(\R^N)$ le plongement de Plücker de $\Gr_{r,n}$.

Alors l'image de $\psi_{r,n}$ est l'ensemble des points $v=(v_K)_{K\in\Lambda(r,n)}$ de $\P(\R^N)$ vérifiant le système d'équations suivant :

\[\sum_{\substack{1\le k\le r+1\\j_k\notin I}} (-1)^{k+\iota(I,\{j_k\})} v_{I\cup\{j_k\}} v_{J\setminus\{j_k\}}=0\]

pour tout $I\in\Lambda(r-1,n)$ et pour tout $J=\{j_1<\cdots<j_{r+1}\}\in\Lambda(r+1,n)$.

Ces équations sont appelées les \emph{relations de Plücker}.

\end{theoreme}

Maintenant que les coordonnées de Plücker ont été définies, on donne un corollaire du développement de Laplace (proposition \ref{devlaplace_aomrgihmoaeiv}) en termes de coordonnées de Plücker.

\begin{corollaire}[développement de Laplace]\label{corollaire_dev_Laplace_amiovomvbnainv}

Soient $n\ge 2$, \hbox{$a,b\in\{1,\ldots,n-1\}$} tels que $a+b=n$, $M_A\in\M_{n,a}(\R)$ et $M_B\in\M_{n,b}(\R)$ deux matrices de rangs $a$ et $b$ respectivement. Notons $M$ la matrice carrée définie comme suit, par blocs :

\[M=\begin{pmatrix} M_A&M_B\end{pmatrix}\in\M_n(\R).\]

Notons $A$ et $B$ les sous-espaces vectoriels de $\R^n$ engendrés par les colonnes de $M_A$ et $M_B$ respectivement. Posons $N=\binom na=\binom nb$, et notons $(\zeta_1,\ldots,\zeta_N)$ les coordonnées de Plücker de $A$ associées à $M_A$ et classées par ordre lexicographique, et $(\eta_1,\ldots,\eta_N)$ celles de $B$ associées à $M_B$ et classées par ordre lexicographique.

Alors il existe une fonction $\epsilon$ à valeurs dans $\{\pm 1\}$ telle que 

\[\det(M)=\sum_{i=1}^N \epsilon(i)\zeta_i\eta_{N+1-i}.\]

\end{corollaire}

On finit avec un lemme permettant de montrer qu'on peut trouver un représentant des coordonnées de Plücker de n'importe quel sous-espace rationnel, tel que celles-ci soient entières et premières entre elles. On inclut une preuve de ce lemme, faute d'avoir trouvé une référence pour celui-ci.

\begin{lemme}\label{lemme_pour_coordPluck_premsentreelles_egouzbelvc}

Soient $n\ge 2$ un entier, $e\in\{1,\ldots,n\}$ et $B\in\mathfrak R_n(e)$. Il existe une base $(X_1,\ldots,X_e)$ de $B\cap\Z^n$ telle qu'en notant $\eta=(\eta_1,\ldots,\eta_N)$, où \hbox{$N=\binom ne$}, les coordonnées de Plücker associées à $(X_1,\ldots,X_e)$ et ordonnées selon l'ordre lexicographique, on ait $\eta\in\Z^N$ et 

\[\pgcd(\eta_1,\ldots,\eta_N)=1.\]

\end{lemme}

\begin{preuve}
Comme $B$ est un sous-espace rationnel, $B\cap\Z^n$ est un sous-$\Z$-module du $\Z$-module libre $\Z^n$. D'après le théorème de la base adaptée, il existe une base $(X_1,\ldots,X_n)$ de $\Z^n$ et des entiers $d_1,\ldots,d_e\ge 1$ tels que $(d_1X_1,\ldots,d_eX_e)$ soit une base de $B\cap\Z^n$. \\

Soit $i\in\{1,\ldots,e\}$. Comme $d_iX_i\in B\cap\Z^n$ et que $X_i\in\Z^n$, on a $X_i\in B\cap\Z^n$, d'où $d_i=1$. \\

Finalement, $(X_1,\ldots,X_e)$ est une base de $B\cap\Z^n$.\\

Notons $M$ la matrice de $\M_n(\Z)$ dont les colonnes sont respectivement $X_1,\ldots,X_n$. Notons $M_1$ la matrice de $\M_{n,e}(\Z)$ formée des $e$ premières colonnes de $M$ et $M_2$ la matrice de $M_{n,n-e}(\Z)$ formée des $n-e$ dernières colonnes de $M$. \\

Remarquons que les mineurs de taille $e\times e$ de $M_1$ donnent le représentant \hbox{$(\eta_1,\ldots,\eta_N)\in\Z^N$}, des coordonnées de Plücker de $B$ associées à sa base $(X_1,\ldots,X_e)$. Notons $\delta_1,\ldots,\delta_N$ les mineurs de taille $(n-e)\times (n-e)$ de $M_2$ dans l'ordre lexicographique. \\

En appliquant la proposition \ref{devlaplace_aomrgihmoaeiv} du développement de Laplace, on obtient

\[\abs{\det M}=\abs{\sum_{i=1}^N \epsilon(i)\eta_i\delta_{N+1-i}}\]

où $\epsilon$ est une fonction à valeurs dans $\{\pm1\}$. Or

\[\abs{\det M}=\covol(\Z^n)=1,\]

donc

\[\abs{\sum_{i=1}^N \epsilon(i)\eta_i\delta_{N+1-i}}=1,\]

ce qui est une relation de Bézout généralisée. Finalement,

\[\pgcd(\eta_1,\ldots,\eta_N)=1.\] 
\end{preuve}

\subsection{Déterminant généralisé}\label{ss_section_det_gene_omaiosbv}

Donnons ici plusieurs résultats montrés par Schmidt \cite{schmidt67} sur la complexité et la proximité des sous-espaces vectoriels de $\R^n$, qui ont en commun la notion suivante qui n'est rien d'autre que la racine carrée du déterminant de Gram. \\

\begin{definition}\label{def_det_gene_namobfd}

Soit $(X_1,\ldots,X_\ell)$ une famille de vecteurs de $\R^n$. Notons \hbox{$\M\in\M_{n,\ell}(\R)$} la matrice dont la $j$-ème colonne est $X_j$ pour tout $j\in\{1,\ldots,\ell\}$. On définit le \emph{déterminant généralisé} de la famille $(X_1,\ldots,X_\ell)$ comme

\[D(X_1,\ldots,X_\ell)=\sqrt{\det(\transp MM)}=\sqrt{\mathrm{Gram}(X_1,\ldots,X_\ell)}.\]

\end{definition}

Le déterminant généralisé possède les propriétés suivantes, données dans \cite{schmidt67} page 434.

\begin{proposition}\label{propriete_det_gene_rqgeuhlgzuh}

Soit $X_1,\ldots,X_\ell$ une famille de vecteurs de $\R^n$. On a

\begin{enumerate}[$(i)$]

	\item $D(X_1,\ldots,X_\ell)\ge 0$ ;
	
	\item $D(X_1,\ldots,X_\ell)=0$ si, et seulement si, $(X_1,\ldots,X_\ell)$ est liée ;
	
	\item $D(X_1,\ldots,X_\ell)=D(X_{\sigma(1)},\ldots,X_{\sigma(\ell)})$ pour toute permutation $\sigma\in\mathfrak S_\ell$ ;
	
	\item $D(\rho(X_1),\ldots,\rho(X_\ell))=D(X_1,\ldots,X_\ell)$ si $\rho$ est une transformation orthogonale ;
	
	\item $D(X_1,\ldots,tX_k,\ldots,X_\ell)=\abs tD(X_1,\ldots,X_\ell)$ pour tout $t\in\R$ et pour tout $k\in\{1,\ldots,\ell\}$ ;
	
	\item $D(X_1,\ldots,X_\ell)=D(X_1,\ldots,X_k+cX_m,\ldots,X_\ell)$ pour tout $c\in\R$ et pour tout $m\ne k$.
	
\end{enumerate}

\end{proposition}

Schmidt définit dans l'équation (6) page 446 de son article une quantité qui regroupe tous les angles canoniques entre deux sous-espaces $A$ et $B$ tels que \hbox{$\min(\dim A,\dim B)=t$} : 

\begin{equation}\label{eqdef_varphi_armibaoibv}
\varphi(A,B)=\prod_{i=1}^t\psi_i(A,B).
\vspace{2.6mm}
\end{equation}

Il démontre page 446 de \cite{schmidt67} le résultat suivant :

\begin{proposition}\label{lien_angletotal_detgene_orzuglziuhg}

Soient $(X_1,\ldots,X_d)$ une base de $A$ et $(Y_1,\ldots,Y_e)$ une base de $B$, alors

\[\varphi(A,B)=\frac{D(X_1,\ldots,X_d,Y_1,\ldots,Y_e)}{D(X_1,\ldots,X_d)D(Y_1,\ldots,Y_e)}.\]

\end{proposition}

Ce résultat montre que la proximité de $A$ et $B$ s'exprime en fonction de déterminants généralisés. C'est en fait surtout celui du numérateur qui importe ici : on va voir (remarque \ref{la_hauteur_comme_det_gene_aruihglnomfehdv} ci-dessous) que le dénominateur est le produit des hauteurs de $A$ et $B$, si les bases sont bien choisies. \\

En effet, cette notion conduit à une définition équivalente de la hauteur d'un sous-espace rationnel (proposition \ref{defequiv_hauteur_aveclecovol_eohuflbg}). Pour relier déterminant généralisé et hauteur, fixons quelques notions de géométrie.

\begin{definition}

On appelle \emph{réseau} de $\R^n$ tout $\Z$-module engendré par des vecteurs $\R$-linéairement indépendants dans $\R^n$.

Soit $\Gamma$ un réseau de $\R^n$ ; on appelle \emph{rang} de $\Gamma$, noté $\rg(\Gamma)$, le rang du $\Z$-module associé. Une \emph{base} de $\Gamma$ est une base de ce $\Z$-module.

\end{definition}

\begin{definition}

Soient $\Gamma$ un réseau de $\R^n$ de rang $e\in\{1,\ldots,n\}$ et $(x_1,\ldots,x_e)$ une base de $\Gamma$. Notons $X\in\M_{n,e}(\R)$ la matrice dont les colonnes sont respectivement $x_1,\ldots,x_e$. On définit le \emph{covolume} de $\Gamma$, noté $\covol (\Gamma)$, comme 

\[\covol (\Gamma)=\abs{\det X}.\]

\end{definition}

Géométriquement, le covolume d'un réseau correspond à la mesure d'un parallélépipède minimal, \emph{i.e.} de mesure minimale, de ce réseau (pour la mesure de Lebesgue sur $\vect(\Gamma)$). Autrement dit, c'est le volume d'une maille de ce réseau, et celui du quotient $\vect(\Gamma)/\Gamma$. Cette notion conduit à la proposition suivante (théorème 1 page 435 de \cite{schmidt67}).

\begin{proposition}\label{defequiv_hauteur_aveclecovol_eohuflbg}

Soit $B$ un sous-espace rationnel de $\R^n$. On a

\[H(B)=\covol(B\cap \Z^n).\]

\end{proposition}

Ainsi, la définition \ref{def_hauteur_coord_plucker_reuighgzmne} et la proposition \ref{defequiv_hauteur_aveclecovol_eohuflbg} fournissent deux points de vue équivalents sur la hauteur d'un sous-espace vectoriel. La définition \ref{def_hauteur_coord_plucker_reuighgzmne} donne une vision algébrique de la hauteur, quand la proposition \ref{defequiv_hauteur_aveclecovol_eohuflbg} en propose une plus géométrique. \\

Pour relier ceci au déterminant généralisé, on peut remarquer (voir page 298 de \cite{fischer14}) que si $X_1,\ldots,X_\ell$ sont des vecteurs de $\R^n$, $D(X_1,\ldots,X_\ell)$ est le $\ell$-volume du parallélotope défini par $X_1,\ldots,X_\ell$ (\emph{i.e.} le volume de ce parallélotope vu dans l'espace euclidien $\vect(X_1,\ldots,X_\ell)$ s'il est de dimension $\ell$, $0$ sinon). On en déduit la remarque suivante :

\begin{remarque}\label{la_hauteur_comme_det_gene_aruihglnomfehdv}

Si $B\in\mathfrak R_n(e)$ et que $(X_1,\ldots,X_e)$ est une $\Z$-base de $B\cap\Z^n$, alors

\[H(B)=D(X_1,\ldots,X_e).\]

\end{remarque}

\subsection{Comportement de la hauteur vis-à-vis d'une transformation linéaire}\label{ss_section_hauteur_transfo_lineaire_aboimfbevbv}

Cette sous-section s'inspire du résultat énoncé dans l'équation (6) page 433 de \cite{schmidt67}, sur le comportement de la hauteur vis-à-vis d'une transformation linéaire inversible. On étend ici ce résultat à certains cas où le morphisme n'est pas inversible. La preuve est très similaire à celle page 433 de \cite{schmidt67}. \\

Ici et dans toute la suite, $\norme\cdot$ désignera la norme euclidienne canonique.

\begin{proposition}\label{inegalitesurlahauteur_aroibvaoivbnoi}

Soient $n\ge 3$ et $e,p\in\{1,\ldots,n\}$ ; soient $B\in\mathfrak R_e(n)$ et $F$ deux sous-espaces vectoriels rationnels de $\R^n$ tels que $B\subset F$ ; soit $\varphi\colon F\to \R^p$ un morphisme linéaire rationnel tel que 

\[\dim \varphi (B)=\dim B.\]

Alors il existe une constante $c(\varphi)>0$, indépendante de $B$, telle que

\[H(\varphi(B))\le c(\varphi)H(B).\]

\end{proposition}

\begin{preuve}
On prolonge $\varphi$ en un endomorphisme rationnel de $\R^n$ en étendant son espace d'arrivée de $\R^p$ à $\R^n$ et en posant $\varphi(x)=x$ pour tout $x\in F^\perp$. \\

On commence par démontrer le résultat pour un sous-espace rationnel $L$ de dimension $1$. \\

Soit $\xi=(\xi_1,\ldots,\xi_n)\in\Z^n$ tel que $\pgcd(\xi_1,\ldots,\xi_n)=1$ et

\[L=\vect(\xi).\]

On a alors $\varphi(L)=\vect(\varphi(\xi))$, et il existe $c_1(\varphi)>0$ indépendante de $\xi$ telle que

\[\norme {\varphi(\xi)}\le c_1(\varphi) \norme \xi.\]

Notons $(\zeta_1,\ldots,\zeta_n)\in\Q^n$ les coordonnées de $\varphi(\xi)$.

Comme $\varphi\in\M_n(\Q)$, il existe $k\in\Z$ tel que $k\varphi\in\M_n(\Z)$. Alors $k\varphi(\xi)\in\Z^n$, donc

\[\forall i\in\{1,\ldots,n\},\quad k\zeta_i\in\Z.\]

En notant $\mathfrak b$ l'idéal fractionnaire engendré par les $\zeta_i$, on a alors

\[k\mathfrak b=k(\zeta_1\Z+\cdots+\zeta_n\Z)\subset \Z,\]

d'où $kN(\mathfrak b)=N(k\mathfrak b)\ge 1$. Ceci donne en utilisant la remarque \ref{def_plus_generale_de_la_hauteur_aobmirnfmeobnv} que

\[H(\varphi(L))=N(\mathfrak b)^{-1} \norme{\varphi(\xi)}\le k c_1(\varphi) \norme{\xi}=c(\varphi) \norme{\xi}=c(\varphi)H(L)\]

en posant $c(\varphi)=kc_1(\varphi)$. \\

Montrons maintenant le résultat pour un sous-espace rationnel $B$ de dimension $e$. On pose $N=\binom ne$. \\

On note $B^*$ la droite rationnelle de $\R^N$ engendrée par les coordonnées de Plücker de $B$. Comme par hypothèse $\dim\varphi(B)=\dim B$, la droite rationnelle $\varphi(B)^*$ engendrée par les coordonnées de Plücker de $\varphi(B)$ appartient elle aussi à $\R^N$. \\

D'après la définition \ref{def_hauteur_coord_plucker_reuighgzmne} de la hauteur, on a

\[H(B)=H(B^*)\quad \text{ et } \quad H(\varphi(B))=H(\varphi(B)^*).\]

Par ailleurs, si $S\in\M_n(\Q)$ est la matrice de $\varphi$ dans la base canonique de $\R^n$, alors $\Lambda^e(S)\in\M_N(\Q)$ est la matrice de $\varphi^{(e)}$, la $e$-ième puissance composée de $\varphi$, dans la base canonique de $\R^N$. D'après \cite{schmidt67} page 433, on a alors 

\[\varphi(B)^*=\varphi^{(e)}(B^*),\]

donc

\[H(\varphi(B)^*)=H(\varphi^{(e)}(B^*)).\]

On est alors ramené au cas de la dimension $1$ dans $\R^N$, et il est possible d'appliquer le début de la preuve pour conclure. \\

On remarque que la constante $c(\varphi)$ ne dépend pas de $e$ quitte à considérer

\[c(\varphi)=\max_{1\le e\le n} {c}^{(e)}(\varphi).\]
\end{preuve}

\begin{remarque}

La proposition \ref{inegalitesurlahauteur_aroibvaoivbnoi} s'applique à tout $B\in\mathfrak R_n(e)$ si $\varphi$ est un automorphisme de $\R^n$. 

\end{remarque}

\section{Sur la proximité}\label{section_proximite_aobomebvob}

\begin{figure}[H]
\begin{center}
\includegraphics[scale=0.65]{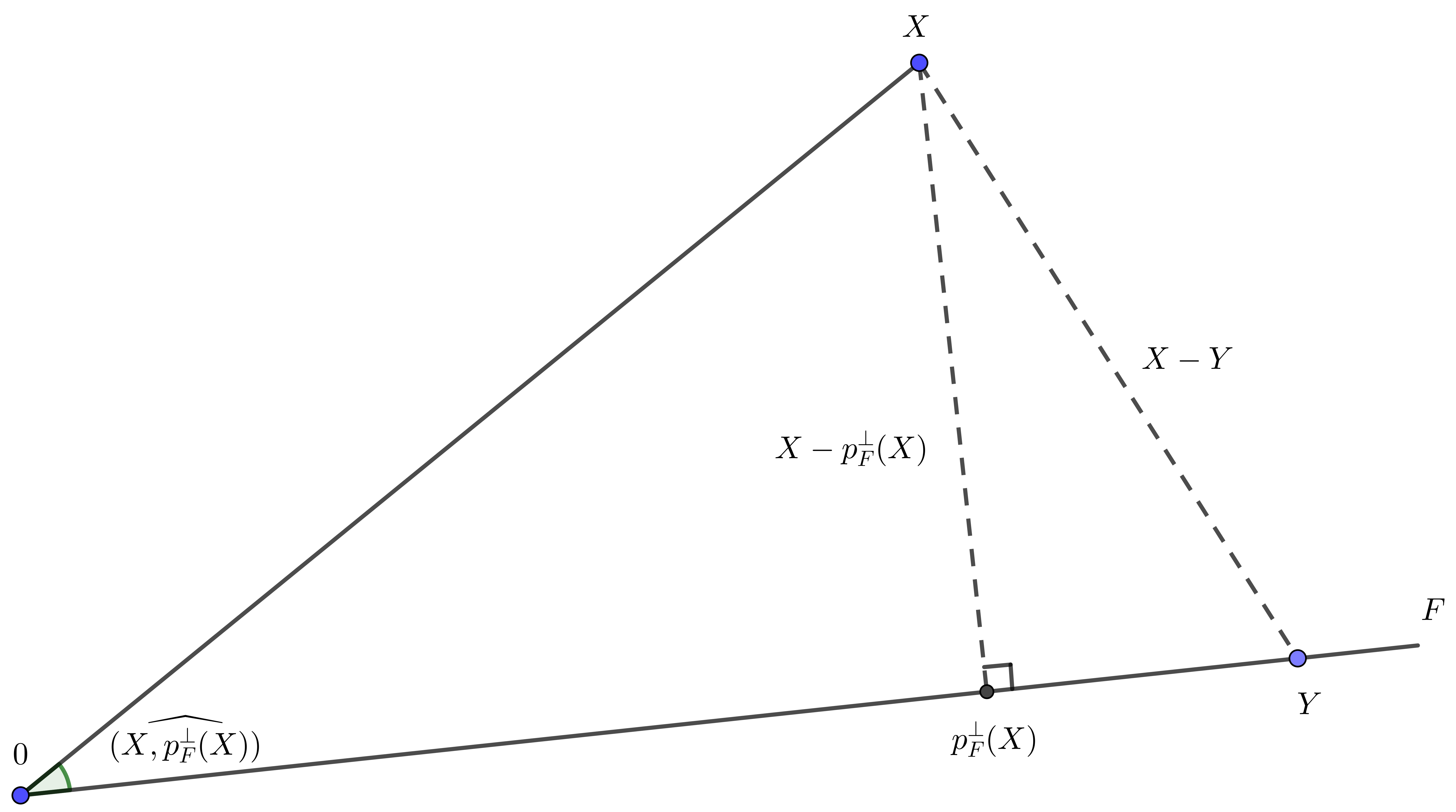}
\caption{Illustration des lemmes \ref{lemme_proximite_proj_ortho_apamoebnamon} et \ref{lemme_maj_psiXY_amoignaoivnovibs}}
\label{lemmes_geom_dessin_1}
\end{center}
\end{figure}

Dans cette section sont exposés plusieurs outils sur les angles entre deux sous-espaces vectoriels. Un premier lemme permet de calculer la proximité entre un vecteur et son projeté orthogonal sur un sous-espace vectoriel. Ce lemme est illustré sur la figure \ref{lemmes_geom_dessin_1} ci-dessous.

\begin{lemme}\label{lemme_proximite_proj_ortho_apamoebnamon}

Soient $F$ un sous-espace vectoriel de $\R^n$, $p_F^\perp$ la projection orthogonale sur $F$ et $X\in\R^n\setminus F^\perp$. Alors 

\[\psi(X,p_F^\perp(X))=\frac{\norme{X-p_F^\perp(X)}}{\norme X}.\]

\end{lemme}

\begin{preuve}
On utilise la définition géométrique du sinus dans le triangle rectangle formé par les vecteurs $X$ et $p_F^\perp(X)$.
\end{preuve}

Le lemme suivant permet de majorer l'angle entre deux vecteurs en fonction de leurs normes. Il est aussi illustré sur la figure \ref{lemmes_geom_dessin_1}.

\begin{lemme}\label{lemme_maj_psiXY_amoignaoivnovibs}

Soient $X$ et $Y$ deux vecteurs non nuls. On a 

\[\psi(X,Y)\le\frac{\norme{X-Y}}{\norme X}.\]

\end{lemme}

\begin{preuve}
On applique le lemme \ref{lemme_proximite_proj_ortho_apamoebnamon} avec $F=\vect(Y)$. On note $p_{\vect(Y)}^\perp$ la projection orthogonale sur $\vect(Y)$. D'après le théorème de Pythagore, on a

\[\norme{X-Y}^2=\norme{X-p_{\vect(Y)}^\perp(X)}^2+\norme{Y-p_F^\perp(X)}^2,\]

donc

\[\norme{X-p_{\vect(Y)}^\perp(X)}\le \norme{X-Y}\]

et le résultat découle du lemme \ref{lemme_proximite_proj_ortho_apamoebnamon}.
\end{preuve}

\begin{figure}[H]
\begin{center}
\includegraphics[scale=0.27]{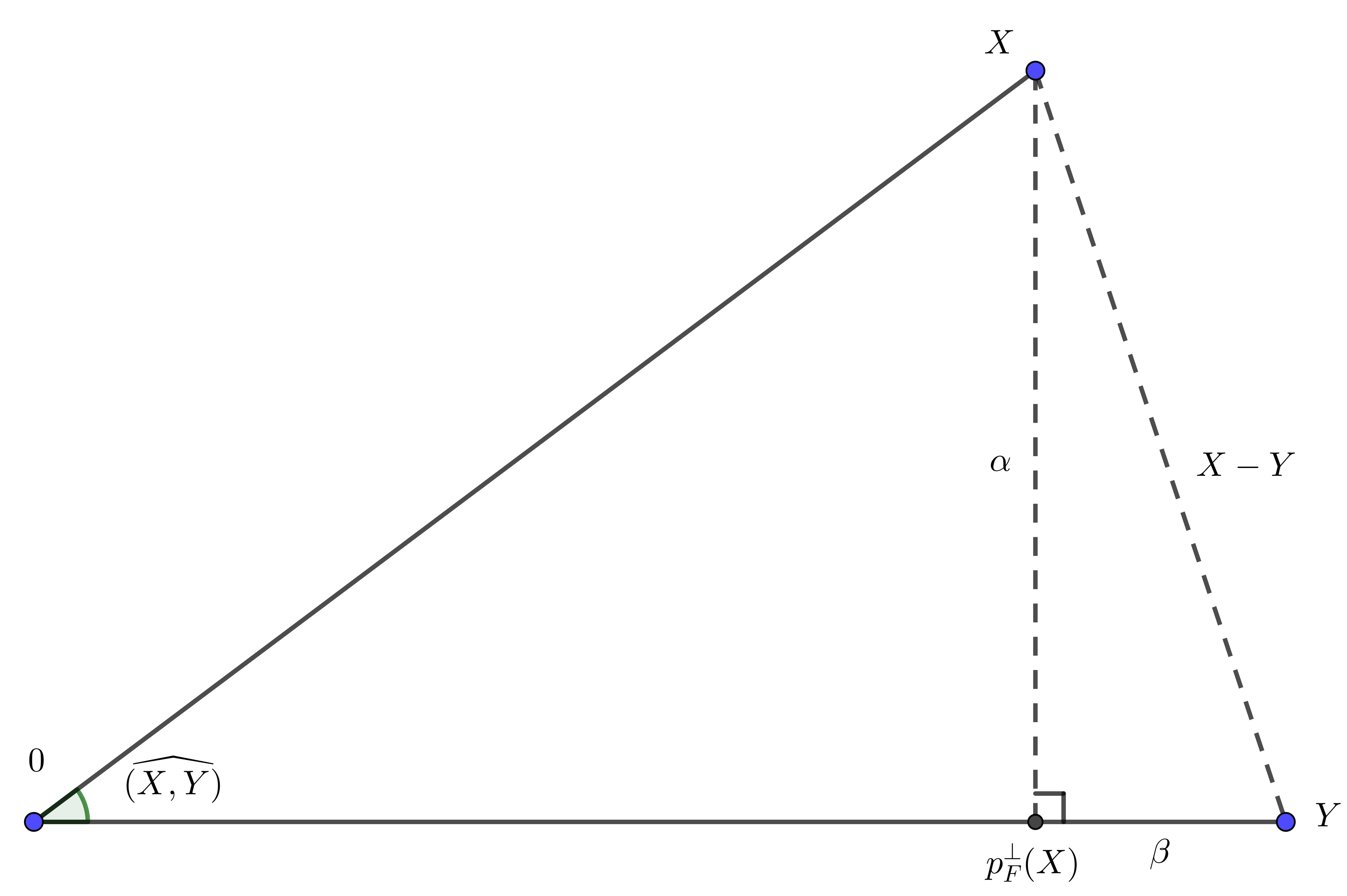}
\caption{Illustration du lemme \ref{lemme_minoration_psiXY_amoifbvoisbdvoibs}}
\label{lemmes_geom_dessin_2}
\end{center}
\end{figure}

À l'inverse du lemme \ref{lemme_maj_psiXY_amoignaoivnovibs}, le prochain lemme permet de minorer l'angle entre deux vecteurs en fonction de leurs normes. Celui-ci est illustré sur la figure \ref{lemmes_geom_dessin_2}. 

\begin{lemme}\label{lemme_minoration_psiXY_amoifbvoisbdvoibs}

Soient $X$ et $Y$ deux vecteurs unitaires tels que $X\cdot Y\ge 0$. On a

\[\psi(X,Y)\ge \frac{\sqrt 2}{2} \norme{X-Y}.\]

\end{lemme}

\begin{preuve}
On note $p_{\vect(Y)}^\perp$ la projection orthogonale sur $\vect(Y)$, ainsi que

\[\alpha=\norme{X-p_{\vect(Y)}^\perp(X)}\quad\text{ et }\quad\beta=\norme{Y-p_{\vect(Y)}^\perp(X)}.\]

D'après le théorème de Pythagore, on a 

\[\norme{X-Y}^2=\alpha^2+\beta^2.\]

Or $X$ est unitaire, donc d'après le lemme \ref{lemme_proximite_proj_ortho_apamoebnamon}, on a

\[\psi(X,Y)=\psi(X,p_{\vect(Y)}^\perp(X))=\norme{X-p_{\vect(Y)}^\perp(X)}=\alpha.\]

De plus $X\cdot Y\ge 0$, donc en appliquant à nouveau le théorème de Pythagore, on a

\[1=\norme{X}^2=(1-\beta)^2+\alpha^2.\]

Il existe donc $\theta\in[0,\pi/2]$ tel que $1-\beta=\cos\theta$ et $\alpha=\sin\theta$. Or

\[\forall \theta\in[0,\pi/2],\quad 1-\cos\theta\le \sin\theta,\]

donc

\[\beta\le \alpha=\norme{X-p_{\vect(Y)}^\perp(X)}.\]

Finalement,

\[\norme{X-Y}^2\le 2\alpha^2=2\psi(X,Y)^2.\]
\end{preuve}

\begin{figure}[H]
\begin{center}
\includegraphics[scale=0.77]{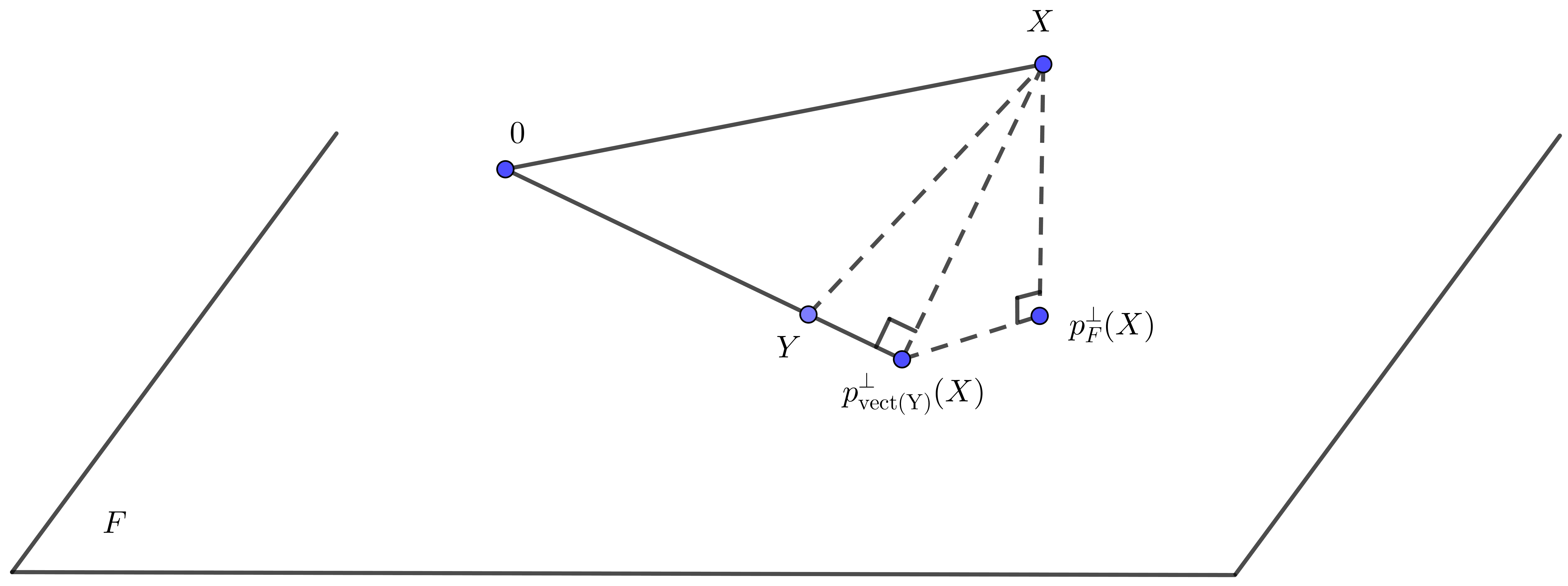}
\caption{Illustration du lemme \ref{lemme_angle_projete_ortho_amotinfqovn}}
\label{lemme_geom_dessin_3}
\end{center}
\end{figure}

Le lemme ci-dessous relie l'angle entre un vecteur et un sous-espace vectoriel, à l'angle entre le vecteur et son projeté orthogonal sur le sous-espace. Il est illustré sur la figure \ref{lemme_geom_dessin_3} ci-dessous.

\begin{lemme}\label{lemme_angle_projete_ortho_amotinfqovn}

Soient $X\in\R^n\setminus\{0\}$ et $F$ un sous-espace vectoriel de $\R^n$ tel que $\dim F\in\{1,\ldots,n-1\}$. Notons $p_F^\perp$ la projection orthogonale sur $F$. On a

\[\psi_1(\vect(X),F)=\psi(X,p_F^\perp(X)).\] 

\end{lemme}

\begin{preuve}
Sans perte de généralité, on peut supposer que $X$ est unitaire ; soit $Y\in F\setminus\{0\}$. Si $X\in F$ le résultat est trivial. Sinon, on note $p_{\vect(Y)}^\perp$ la projection orthogonale sur $\vect(Y)$, le triangle dont les sommets sont les vecteurs $X$, $p_{\vect(Y)}^\perp(X)$ et $p_F^\perp(X)$ est rectangle en $p_F^\perp(X)$. Ainsi, 

\[\norme{X-p_F^\perp(X)}\le \norme{X-p_{\vect (Y)}^\perp(X)}.\]

D'après le lemme \ref{lemme_proximite_proj_ortho_apamoebnamon} et comme $X$ est un vecteur unitaire, on a

\[\psi(X,p_F^\perp(X))=\norme{X-p_F^\perp(X)}\]

donc toujours d'après le lemme \ref{lemme_proximite_proj_ortho_apamoebnamon} 

\[\psi(X,Y)=\psi(X,p_{\vect(Y)}^\perp(X))=\norme{X-p_{\vect (Y)}^\perp(X)},\] 

donc

\[\psi(X,p_F^\perp(X))\le \psi(X,Y).\]

D'après l'équation \eqref{def_du_psi_1_piaeofdbouefbv} définissant $\psi_1$, on a donc

\[\psi_1(\vect(X),F)=\min_{Z\in F\setminus\{0\}}\psi(X,Z)=\psi(X,p_F^\perp(X)).\]
\end{preuve}

Les deux résultats suivants correspondent respectivement au lemme 12 et à son corollaire dans l'article \cite{schmidt67} page 444.

\begin{proposition}[Schmidt]\label{prop_proximite_aveclesmin_maoierngmoinv}

Soient $A$ et $B$ deux sous-espaces vectoriels de $\R^n$ de dimensions respectives $d$ et $e$. Alors pour tout $j\in\{1,\ldots,\min(d,e)\}$, $\psi_j(A,B)$ est le plus petit réel $\psi\ge 0$ tel qu'il existe un sous-espace $A'\subset A$ de dimension $j$ tel que

\[\forall X\in A'\setminus\{0\},\quad \exists Y\in B\setminus\{0\},\quad \psi(X,Y)\le \psi.\]

\end{proposition}

\begin{corollaire}[Schmidt]\label{proximite_et_inclusions_aeoibnfvoidsn}

Soient $A'\subset A$ et $B'\subset B$ quatre sous-espaces vectoriels de $\R^n$. Alors

\[\forall j\in\{1,\ldots,\min(A',B')\},\quad \psi_j(A,B)\le \psi_j(A',B').\]

\end{corollaire}

On déduit de la proposition \ref{prop_proximite_aveclesmin_maoierngmoinv} le lemme suivant. 

\begin{lemme}\label{lemme_transfert_psi1_psiell_baoeoearv}

Soient $A$ et $B$ deux sous-espaces vectoriels non triviaux de $\R^n$ tels que \hbox{$\dim A\le \dim B$}. Alors

\[\forall X\in A\setminus\{0\},\quad \psi_1(\vect(X),B)\le\psi_{\dim A}(A,B).\]

\end{lemme}

\begin{preuve}
Soit $X\in A\setminus\{0\}$. On a
\begin{align*}
\psi_1(\vect(X),B)
	&=\min_{Y\in B\setminus\{0\}}\psi(X,Y) \\
	&\le \max_{Z\in A\setminus\{0\}}\min_{Y\in B\setminus\{0\}}\psi(Z,Y) \\
	&=\min\{\varphi,\ \forall Z\in A\setminus\{0\},\quad \exists Y\in B\setminus\{0\},\quad \psi(Z,Y)\le \varphi\} \\
	&=\psi_{\dim A}(A,B)
\end{align*}
d'après la proposition \ref{prop_proximite_aveclesmin_maoierngmoinv}.
\end{preuve}

Mentionnons une inégalité triangulaire entre les angles (\cite{schmidt67}, équation (3) page 446).

\begin{proposition}[Schmidt]\label{inegalite_triang_angles_aeoifgnivn}

Pour tous $X,Y,Z\in\R^n$ non nuls, on a

\[\psi(X,Z)\le \psi(X,Y)+\psi(Y,Z).\]

\end{proposition}

Finalement, énonçons une proposition sur la proximité entre deux sous-espaces transformés par une transformation inversible. Celle-ci regroupe la fin du théorème \ref{la_raison_pour_laquelle_les_angles_sont_canoniques} et le lemme 13 page 446 de l'article \cite{schmidt67}.

\begin{proposition}[Schmidt]\label{proximite_transfo_rationnelle_inv_moiefhvoizcnzou}

Soient $A$ et $B$ deux sous-espaces vectoriels de $\R^n$ de dimensions respectives $d$ et $e$, et $\varphi\in\GL_n(\R)$. Il existe une constante $c_\varphi>0$ (indépendante de $A$ et de $B$) telle que

\[\forall j\in\{1,\ldots,\min(d,e)\},\quad \psi_j(\varphi(A),\varphi(B))\le c_\varphi\psi_j(A,B).\]

De plus, si $\varphi$ est une isométrie, alors

\[\forall j\in\{1,\ldots,\min(d,e)\},\quad \psi_j(\varphi(A),\varphi(B))=\psi_j(A,B).\]

\end{proposition}

\section{Deux résultats d'approximation simultanée}\label{section_resultats_approx_simult_aoeibhaoimnvao}

Dans cette section, on énonce deux théorèmes d'approximation simultanée qui serviront dans les chapitres \ref{chap_cas_particuliers_oamofbvoaube} et \ref{chap_somme_directes_vaoinvoin}. \\

Commençons par le théorème d'approximation simultanée de Dirichlet (voir \cite{hardy07} théorème 200 page 216). 

\begin{theoreme}[Dirichlet]\label{approx_simul_dirichlet_aemofibbv}

Soient $d\ge 1$ entier et $x\in\R^d\setminus\Q^d$. Il existe une infinité de couples $(p,q)\in\Z^d\times\N^*$ tels que $\pgcd(p_1,\ldots,p_d,q)=1$ et

\[\norme{x-\frac pq}_\infty\le \frac{1}{q^{1+1/d}}.\]

\end{theoreme}

Mentionnons aussi le corollaire du théorème 2 de l'article \cite{schmidt70} qui permet de minorer les valeurs entières d'une forme linéaire à coefficients algébriques.

\begin{proposition}[Schmidt, 1970]\label{prop_approx_formelin_nombrealg_neorifmnbv}

Soient $\alpha_1,\ldots,\alpha_k$ des nombres algébriques tels que $1,\alpha_1,\ldots,\alpha_k$ soient linéairement indépendants sur $\Q$. Alors pour tout $\epsilon>0$, il existe une constante $c>0$, telle que pour tous entiers $q_1,\ldots,q_k,p$ tels que $\max(\abs{q_1},\ldots,\abs{q_k})>0$, on ait

\[\abs{p+\sum_{i=1}^k q_i\alpha_i}\ge c\max(\abs{q_1},\ldots,\abs{q_k})^{-k-\epsilon}.\]

\end{proposition}

\section{Théorèmes de Going-up et Going-down}\label{section_goingupdown_aogieovb}

Ces théorèmes sont démontrés dans \cite{schmidt67} (théorèmes 9 et 10 page 453, en faisant attention, car il y a un "$+1$" au lieu d'un "$-1$" dans le théorème du Going-up, et un signe "$-$" oublié dans le théorème du Going-down), et permettent de déduire des résultats en changeant la dimension du sous-espace vectoriel approchant. 

\begin{theoreme}[Going-up, Schmidt, 1967]\label{goingup_eorihgzmefvnfon}

Soient $d,e\in\N^*$ tels que $d+e<n$ ; posons $t=\min(d,e)$. Soient $A$ un sous-espace de $\R^n$ de dimension $d$ et $B\in\mathfrak R_n(e)$. Soit $H\ge 1$ tel que $H(B)\le H$, et tel qu'il existe des $x_i,y_i\in\R$ tels que 

\[\forall i\in\{1,\ldots,t\},\quad H(B)^{x_i}\psi_i(A,B)\le c_1H^{-y_i}\]

avec $c_1>0$.

Alors il existe une constante $c_2>0$ dépendant uniquement de $n$ et de $e$, et une constante $c_3>0$ dépendant uniquement de $n$, $e$, $x_i$ et $y_i$, telles qu'en posant 

\[H'=c_2H^{(n-e-1)/(n-e)},\]

il existe $C\in\mathfrak R_n(e+1)$ tel que

\[\begin{cases} C\supset B \\H(C)\le H'\\\forall i\in\{1,\ldots,t\},\quad H(C)^{x_i(n-e)/(n-e-1)}\psi_i(A,C)\le c_1c_3H'^{-y_i(n-e)/(n-e-1)}.\end{cases}\]

\end{theoreme}

\begin{theoreme}[Going-down, Schmidt, 1967]\label{goingdown_ogihoeihfgoniv}

Soient $d,e\in\N^*$ tels que $d+e\le n$. Soient $A$ un sous-espace de $\R^n$ de dimension $d$ et $B\in\mathfrak R_n(e)$. Soit $H\ge 1$ tel que $H(B)\le H$. Posons $f=\min(d,e-1)$. Soient $h\in\{1,\ldots,f\}$, $c_1\ge 1$ et 

\[y_1\ge \cdots\ge y_h\ge \frac 1h\]

tels que

\[\forall i\in\{1,\ldots,h\},\quad H(B)\psi_i(A,B)\le c_1H^{-y_i+1}.\]

Alors il existe des constantes $c_2,\ldots,c_5$, qui ne dépendent que de $n,d,e,y_1,\ldots,y_h$ mais pas de $A$, $B$ ou $H$, telles que les propriétés suivantes soient vérifiées.

Posons $y=y_1+\cdots+y_h$ et supposons que

\[\forall i\in\{1,\ldots,h\},\quad y_i'=\frac{y_i e}{y+e-1}\ge 1.\]

Alors il existe $C\in\mathfrak R_n(e-1)$ tel qu'en posant

\[H'=c_2H^{(e+y-1)/e},\]

on ait

\[\begin{cases} C\subset B \\H(C)\le H'\\\forall i\in\{1,\ldots,t\},\quad \psi_i(A,C)\le c_3H(C)^{-y_i'}.\end{cases}\]

De plus, si

\[\forall i\in\{1,\ldots,h\},\quad \psi_i(A,B)=0,\]

alors on pose $y_0'=e/h$, et pour tout $H'\ge c_4H>0$, il existe $C\in\mathfrak R_n(e-1)$ tel que

\[\begin{cases} C\subset B \\H(C)\le H'\\\forall i\in\{1,\ldots,t\},\quad \psi_i(A,C)\le c_5H^{y_0'}H(C)^{-y_0'}.\end{cases}\]

\end{theoreme}

\chapter{Différents cas particuliers}\label{chap_cas_particuliers_oamofbvoaube}

Ce chapitre est découpé selon la dimension de l'espace ambiant dans lequel sont traités les cas particuliers abordés ici. \\

On commence à la section \ref{section_outils_en_toute_dimension_amoimaoihogierbobi} par quelques lemmes qui serviront en toute dimension ; puis on étudie le cas de $\R^4$ en montrant que $\muexp 4221=3$ (section \ref{section_R4_maoefvbn}) et celui de $\R^5$ en majorant $\muexp 5321$ par $6$ (section \ref{section_R5_gaomfnbamobo}). Après quelques commentaires (section \ref{section_commentaires_généraux_amirbnamoi}) on majore $\muexp {2d}d{d-1}1$ dans $\R^{2d}$ (section \ref{section_app_Moshchevitin_aoimbnoibav}).

\section{Outils en toute dimension}\label{section_outils_en_toute_dimension_amoimaoihogierbobi}

Il est possible de faire le lien entre proximité et hauteur grâce au lemme suivant, dans lequel 

\begin{equation}\label{eqphi_mpairhbnobiv}
\varphi(A,B)=\prod_{j=1}^{\min(\dim A,\dim B)} \psi_j(A,B)
\vspace{2.6mm}
\end{equation}

a été défini dans la sous-section \ref{ss_section_det_gene_omaiosbv} (équation \eqref{eqdef_varphi_armibaoibv}). 

\begin{lemme}\label{lien_proximite_hauteur_zozofvbcz}

Soient $n\ge 2$, $A$ un sous-espace vectoriel de $\R^n$ de dimension $d$ et \hbox{$B\in\mathfrak R_n(e)$} un sous-espace vectoriel rationnel de dimension $e$. On suppose que $d+e=n$. Notons $(X_1,\ldots,X_d)$ une base de $A$ et $(Y_1,\ldots,Y_e)$ une base de $B\cap \Z^n$, ainsi que $M\in\M_n(\R)$ la matrice dont les colonnes sont respectivement $X_1,\ldots,X_d,Y_1,\ldots,Y_e$. 

Alors il existe une constante $c>0$ ne dépendant que de $(X_1,\ldots,X_d)$ telle que

\[\varphi(A,B)=c\ \frac{\abs{\det M}}{H(B)}.\]

\end{lemme}

\begin{preuve}
D'après la proposition \ref{lien_angletotal_detgene_orzuglziuhg}, on a

\[\varphi(A,B)=\frac{D(X_1,\ldots,X_d,Y_1,\ldots,Y_e)}{D(Y_1,\ldots,Y_e)D(X_1,\ldots,X_d)}.\]

Comme $(Y_1,\ldots,Y_e)$ est une base de $B\cap \Z^n$, la remarque \ref{la_hauteur_comme_det_gene_aruihglnomfehdv} découlant du point de vue géométrique sur la hauteur permet d'obtenir

\[\varphi(A,B)=D(X_1,\ldots,X_d,Y_1,\ldots,Y_e)\frac{c}{H(B)}\]

où $c=D(X_1,\ldots,X_d)^{-1}>0$ est une constante ne dépendant que de la base $(X_1,\ldots,X_d)$. \\

De plus, la matrice $M$ est une matrice carrée, donc la définition \ref{def_det_gene_namobfd} donne

\[D(X_1,\ldots,X_d,Y_1,\ldots,Y_e)^2=\det(\transp MM)=\det (M)^2.\]

Ainsi, comme $D(X_1,\ldots,X_d,Y_1,\ldots,Y_e)\ge 0$ (proposition \ref{propriete_det_gene_rqgeuhlgzuh} (i)), on a

\[\varphi(A,B)=c\ \frac{\abs{\det M}}{H(B)}.\]
\end{preuve}

Donnons un autre lemme, qui permet de minorer les angles canoniques en fonction de $\varphi$.

\begin{lemme}\label{minoration_psiphi_eomivocvbn}

Soient $n\ge 2$, et $A$ et $B$ deux sous-espaces vectoriels de $\R^n$ de dimensions respectives $d$ et $e$. Alors

\[\forall j\in\{1,\ldots,\min(d,e)\},\quad \psi_j(A,B)\ge \varphi(A,B)^{1/j}.\]

\end{lemme}

\begin{preuve}
On pose $t=\min(d,e)$ ; soit $j\in\{1,\ldots,t\}$. D'après la définition \ref{def_angles_canoniques_eormighemoidvn} et le théorème \ref{la_raison_pour_laquelle_les_angles_sont_canoniques}, on a 

\[\psi_1(A,B)\le\cdots\le \psi_t(A,B)\le 1.\]

Ainsi, le produit \eqref{eqphi_mpairhbnobiv} peut se découper :

\[\varphi(A,B)=\left(\prod_{i=1}^{j}\underbrace{\psi_i(A,B)}_{\le\psi_j(A,B)}\right)\times\left(\prod_{i=j+1}^t\underbrace{\psi_i(A,B)}_{\le 1}\right)\le \psi_j(A,B)^{j}\]

ce qui conclut la preuve du lemme \ref{minoration_psiphi_eomivocvbn}.
\end{preuve}

\section{Dans $\R^4$}\label{section_R4_maoefvbn}

Comme mentionné en introduction de cette thèse, le cas le plus simple encore ouvert est le cas

\[(n,d,e,j)=(4,2,2,1).\]

En combinant les résultats de Schmidt et Moshchevitin, on a uniquement l'encadrement $3\le \muexp 4221\le 4$. Schmidt mentionne d'ailleurs explicitement ce cas en conclusion de \cite{schmidt67} page 471, et c'est de plus le seul cas encore ouvert dans $\R^4$. \\

La question de déterminer $\muexp 4221$ se trouve complètement résolue grâce au théorème \ref{cas_total_4221_egorihgoegn} :

\begin{theoreme}\label{cas_total_4221_egorihgoegn}

On a

\[\muexp 4221=3.\]

\end{theoreme}

Ainsi, ce théorème clôture le problème \ref{probleme_principal_gnzroign} dans le cas $n=4$.

\subsection{Construction de plans mal approchés de $\mathfrak I_4(2,2)_1$}\label{resultats_dans_R4_jaomfibenovn}

Pour démontrer le théorème \ref{cas_total_4221_egorihgoegn}, on construit explicitement dans cette sous-section des plans de $\R^4$ particulièrement mal approchés par les plans rationnels de $\R^4$. \\

Une famille de sous-espaces vectoriels particuliers sera centrale ici : définissons pour $\xi\in]0,\sqrt 7[$, le sous-espace vectoriel $A_\xi$ de $\R^4$ de dimension $2$ engendré par

\begin{equation}\label{def_Yi_R4_efulkhebovb}
X^{(1)}_\xi=\begin{pmatrix} 0 \\1 \\\xi \\\sqrt{7-\xi^2}\end{pmatrix}\quad \text{ et }\quad X^{(2)}_\xi=\quad\begin{pmatrix} 1 \\0\\-\sqrt{7-\xi^2}\\\xi\end{pmatrix},
\vspace{2.6mm}
\end{equation}

de sorte à faire apparaître les trois nombres $1$, $\xi$ et $\sqrt{7-\xi^2}$ dans les coordonnées de Plücker de $A_\xi$. \\

Différents résultats sont obtenus sur ces sous-espaces, dont les démonstrations se trouvent dans la sous-section \ref{preuves_R4_oiagrhgmaeriofvn}. \\

On montre que pour tout $\xi\in]0,\sqrt 7[$ algébrique de degré au moins égal à $3$, le plan $A_\xi$ :
\begin{enumerate}[$\bullet$]

	\item vérifie la condition de $(2,1)$-irrationalité,
	
	\item est mal approché par les plans rationnels de $\R^4$. \\
	
\end{enumerate} 

\begin{lemme}\label{condition_dirrationalite_4221_regouzhelfgu}

Pour tout $\xi\in]0,\sqrt 7[$ algébrique de degré au moins égal à $3$, on a

\begin{equation}\label{condition_dirr_R4_zodibvdoivn}
A_\xi\in\mathfrak I_4(2,2)_1.
\vspace{2.6mm}
\end{equation}

De plus, pour tous les nombres algébriques $\xi$ vérifiant \eqref{condition_dirr_R4_zodibvdoivn}, on a

\begin{equation}\label{minoration_hauteur_lemme_R4_zidovdbon}
\forall \epsilon>0,\quad \exists c>0,\quad \forall B\in\mathfrak R_4(2),\quad \varphi(A_\xi,B)\ge\frac c{H(B)^{3+\epsilon}}.
\vspace{2.6mm}
\end{equation}

\end{lemme}

\emph{A fortiori}, la condition de $(2,2)$-irrationalité étant plus faible, le lemme \ref{condition_dirrationalite_4221_regouzhelfgu} montre aussi que pour tout $\xi\in]0,\sqrt 7[$ algébrique de degré au moins égal à $3$, on a \hbox{$A_\xi\in\mathfrak I_4(2,2)_2$}. De plus, le lemme \ref{minoration_psiphi_eomivocvbn} montre que le lemme \ref{condition_dirrationalite_4221_regouzhelfgu} est aussi vrai en remplaçant $\varphi(A_\xi,B)$ par $\psi_1(A_\xi,B)$. \\

De ceci découle alors la proposition centrale suivante qui permet de démontrer le théorème \ref{cas_total_4221_egorihgoegn} et ainsi d'apporter une réponse au problème \ref{probleme_principal_gnzroign} pour $n=4$. Cette résolution est de plus explicite. 

\begin{proposition}\label{proposition_R4_mu_Axi_zgiozovdib}

Pour tout nombre algébrique $\xi\in]0,\sqrt 7[$ de degré au moins égal à $3$, on a

\[\muexpA 4{A_\xi}21=3.\]

\end{proposition}

\begin{remarque}

La conclusion de la proposition \ref{proposition_R4_mu_Axi_zgiozovdib} est encore vraie pour tout nombre algébrique $\xi\in]0,\sqrt 7[$ de degré $2$ vérifiant

\[\dim_\Q\vect_\Q(1,\xi,\sqrt{7-\xi^2})=3.\]

En particulier pour $\xi=\sqrt 2$, on a bien 

\[\dim_\Q\vect_\Q(1,\xi,\sqrt{7-\xi^2})=\dim_\Q\vect_\Q(1,\sqrt 2,\sqrt{5})=3,\]

donc 

\[\muexpA 4{A_{\sqrt{2}}}21=3.\]

\end{remarque}

\begin{remarque}

On peut construire une infinité non dénombrable de sous-espaces $A_\xi$ vérifiant \hbox{$\muexpA 4{A_\xi}21=3$} grâce à un théorème de D. Y. Kleinbock et G. A. Margulis. En effet, pour Lebesgue-presque tout $\xi\in]0,\sqrt 7[$, on a $\dim_\Q\vect_\Q(1,\xi,\sqrt{7-\xi^2})=3$. On remplace alors l'utilisation page \pageref{app_pro_2_28_aepihboeifveifbn} de la proposition \ref{prop_approx_formelin_nombrealg_neorifmnbv} dans la preuve du lemme \ref{condition_dirrationalite_4221_regouzhelfgu} par le théorème A page 3 de l'article \cite{kleinbock98}. Celui-ci s'applique car la variété de $\R^3$

\[M=\Big\{(1,\xi,\sqrt{7-\xi^2}),\ \xi\in]0,\sqrt 7[\Big\}\]

est non dégénérée (\emph{i.e.} aucun de ses points n'a de voisinage inclus dans un hyperplan affine de $\R^3$), et l'inégalité \eqref{minorationKB_R4_orubvdobv} est obtenue cette fois-ci pour Lebesgue-presque tout $\xi\in]0,\sqrt 7[$. 

\end{remarque}

\subsection{Les preuves}\label{preuves_R4_oiagrhgmaeriofvn}

Cette sous-section apporte les preuves des résultats énoncés sans démonstration dans la sous-section \ref{resultats_dans_R4_jaomfibenovn} précédente. Notons déjà que le théorème \ref{cas_total_4221_egorihgoegn} suit de la définition de $\mathring\mu$ et de la proposition \ref{proposition_R4_mu_Axi_zgiozovdib}. Par ailleurs, compte tenu du lemme \ref{minoration_psiphi_eomivocvbn} (appliqué avec $j=1$), le lemme \ref{condition_dirrationalite_4221_regouzhelfgu} montre que $\muexpA 4{A_\xi}21\le 3$. Or le théorème \ref{premiere_borne_Schmidt_hgoembeer} donne $\muexpA 4{A_\xi}21\ge \muexp 4221\ge 3$, ce qui démontre la proposition \ref{proposition_R4_mu_Axi_zgiozovdib}. \\

Il ne reste donc qu'à démontrer le lemme \ref{condition_dirrationalite_4221_regouzhelfgu}, qui montre simultanément la condition de $(2,1)$-irrationalité et établit l'inégalité cruciale pour la proposition \ref{proposition_R4_mu_Axi_zgiozovdib}. \\

\begin{preuve}
Soit $B\in\mathfrak R_4(2)$ et $(Y_1,Y_2)$ une base de $B$ donnée par le lemme \ref{lemme_pour_coordPluck_premsentreelles_egouzbelvc}. \\

Comme $\binom 42=6$, notons $(\eta_1,\ldots,\eta_6)$ les coordonnées de Plücker de $B$ associées à la base $(Y_1,Y_2)$. D'après le lemme \ref{lemme_pour_coordPluck_premsentreelles_egouzbelvc}, on a $(\eta_1,\ldots,\eta_6)\in\Z^6$ et 

\begin{equation}\label{premiersentreeux_R4_idhvoevbd}
\pgcd(\eta_1,\ldots,\eta_6)=1.
\vspace{2.6mm}
\end{equation}

De plus, ce vecteur vérifie la relation de Plücker donnée par le théorème \ref{relation_de_Plucker_oerighbfduvb} :

\begin{equation}\label{relationPluck_R4_oehvodbc}
\eta_1\eta_6-\eta_2\eta_5+\eta_3\eta_4=0.
\vspace{2.6mm}
\end{equation}

Explicitons d'où vient la relation \eqref{relationPluck_R4_oehvodbc}. On reprend les notations du théorème \ref{relation_de_Plucker_oerighbfduvb} et de la sous-section \ref{coord_Pluck_teoghzleufve}. On note $(\eta_K)_{K\in\Lambda(2,4)}$ les coordonnées de Plücker de $B$. On a d'après le théorème \ref{relation_de_Plucker_oerighbfduvb} que pour tout $I\in\Lambda(1,4)$ et pour tout \hbox{$J=\{j_1<j_2<j_3\}\in\Lambda(3,4)$} :

\[\sum_{\substack{1\le k\le 3\\j_k\notin I}} (-1)^{k+\iota(I,\{j_k\})} \eta_{I\cup\{j_k\}}\eta_{J\setminus\{j_k\}}=0.\]

En particulier, pour $I=\{1\}$ et $J=\{2,3,4\}$, on obtient 

\[(-1)^{1+\iota(\{1\},\{2\})}\eta_{\{1,2\}}\eta_{\{3,4\}}+(-1)^{2+\iota(\{1\},\{3\})}\eta_{\{1,3\}}\eta_{\{2,4\}}+(-1)^{3+\iota(\{1\},\{4\})}\eta_{\{1,4\}}\eta_{\{2,3\}}=0,\]

soit

\[-\eta_{\{1,2\}}\eta_{\{3,4\}}+\eta_{\{1,3\}}\eta_{\{2,4\}}-\eta_{\{1,4\}}\eta_{\{2,3\}}=0,\]

ce qui donne bien l'équation \eqref{relationPluck_R4_oehvodbc} en faisant correspondre les $\eta_{\{i,j\}}$ et les $\eta_k$ grâce à l'ordre lexicographique sur $\Lambda(2,4)$ :

\begin{equation}\label{correspondance_lexico_moaiboafdvoibf}
(\eta_{\{1,2\}},\eta_{\{1,3\}},\eta_{\{1,4\}},\eta_{\{2,3\}},\eta_{\{2,4\}},\eta_{\{3,4\}})=(\eta_1,\eta_2,\eta_3,\eta_4,\eta_5,\eta_6).
\vspace{2.6mm}
\end{equation}

Notons $M_\xi$ la matrice de $\M_4(\R)$ dont les colonnes sont respectivement $X_\xi^{(1)},X_\xi^{(2)},Y_1,Y_2$. Remarquons que

\[A_\xi\cap B=\{0\}\iff \det M_\xi\ne 0.\]

On calcule alors le déterminant de $M_\xi$ par un développement de Laplace par rapport aux deux premières colonnes (proposition \ref{devlaplace_aomrgihmoaeiv} et exemple \ref{exemple_Laplace_24_bariomnbofeiv} avec la correspondance \eqref{correspondance_lexico_moaiboafdvoibf}), et on trouve :
\begin{align}
\det M_\xi
	&=-\eta_6+\eta_5\xi-\eta_4\sqrt{7-\xi^2}-\eta_3\sqrt{7-\xi^2}-\eta_2\xi+7\eta_1\nonumber \\
	&=-\eta_6+7\eta_1+(\eta_5-\eta_2)\xi+(-\eta_3-\eta_4)\sqrt{7-\xi^2}\label{valeur_detM_R4_zghodidvc}.
\end{align}
Supposons par l'absurde que $\det M_\xi=0$, on a donc

\begin{equation}\label{coord_Pluck_egales0_R4_ziogozivnb}
-\eta_6+7\eta_1+(\eta_5-\eta_2)\xi+(-\eta_3-\eta_4)\sqrt{7-\xi^2}=0. 
\vspace{2.6mm}
\end{equation}

Or pour tout nombre algébrique $\xi\in]0,\sqrt 7[$ de degré au moins égal à $3$, on a

\begin{equation}\label{indep_alg_R4_oersgbozdvbn}
\dim_\Q\vect_\Q(1,\xi,\sqrt{7-\xi^2})=3.
\vspace{2.6mm}
\end{equation}

En effet, supposons par l'absurde qu'il existe $(x,y,z)\in\Q^3\setminus\{(0,0,0)\}$ tel que

\begin{equation}\label{liesurQR4_eroginvoidn}
x+y\xi+z\sqrt{7-\xi^2}=0.
\vspace{2.6mm}
\end{equation}

Comme $\xi\notin\Q$, on peut supposer que $z\ne 0$. En posant $x'=x/z$ et $y'=y/z$, l'équation \eqref{liesurQR4_eroginvoidn} se réécrit

\[(y'^2+1)\xi^2+2x'y'\xi+x'^2-7=0.\]

Or $\xi$ n'est pas un nombre algébrique de degré inférieur ou égal à $2$, et donc ne peut annuler le polynôme

\[(y'^2+1)X^2+2x'y'X+x'^2-7\]

de degré $2$ à coefficients rationnels, ce qui établit l'indépendance linéaire sur $\Q$ de $1,\xi,\sqrt{7-\xi^2}$ pour tout nombre algébrique $\xi\in]0,\sqrt 7[$ de degré au moins égal à $3$. \\

Comme les $\eta_i$ sont entiers, l'équation \eqref{coord_Pluck_egales0_R4_ziogozivnb} sur les coordonnées de Plücker de $B$ combinée avec l'indépendance linéaire \eqref{indep_alg_R4_oersgbozdvbn}, donne que pour tout \hbox{$\xi\in]0,\sqrt 7[$} algébrique de degré au moins $3$, on a

\begin{equation}\label{cases_R4_peariohoeiv}
\begin{cases}\eta_6=7\eta_1\\\eta_5=\eta_2\\\eta_4=-\eta_3.\end{cases}
\vspace{2.6mm}
\end{equation}

Ainsi, la relation \eqref{relationPluck_R4_oehvodbc} de Plücker devient

\[\eta_2^2+\eta_3^2=7\eta_1^2.\]

On regarde cette équation modulo $4$ :

\[\eta_2^2+\eta_3^2\equiv 3\eta_1^2\pmod 4.\]

Or un carré est toujours congru à $0$ ou $1$ modulo $4$, donc $4\mid \eta_1^2$, donc $2\mid \eta_1$. De plus, $4\mid \eta_2^2$ et $4\mid\eta_3^2$, donc $2\mid\eta_2$ et $2\mid\eta_3$. Avec les équations \eqref{cases_R4_peariohoeiv}, on en déduit que $2$ divise tous les $\eta_i$, ce qui est absurde car $\pgcd(\eta_1,\ldots,\eta_6)=1$ d'après \eqref{premiersentreeux_R4_idhvoevbd}. \\

On a donc $\det M_\xi\ne 0$, et ainsi $A_\xi\cap B=\{0\}$ pour tout nombre algébrique $\xi\in]0,\sqrt 7[$ de degré au moins égal à $3$, ce qui établit la condition de $(2,1)$-irrationalité. \\

Il reste à établir l'inégalité \eqref{minoration_hauteur_lemme_R4_zidovdbon}. Pour cela, les deux points de vue sur la hauteur seront utilisés de façon combinée grâce au lemme \ref{lien_proximite_hauteur_zozofvbcz}. \\

Comme la base $(Y_1,Y_2)$ de $B$ donnée par le lemme \ref{lemme_pour_coordPluck_premsentreelles_egouzbelvc} est aussi une $\Z$-base de $B\cap \Z^4$, on peut utiliser le lemme \ref{lien_proximite_hauteur_zozofvbcz} qui donne

\begin{equation}\label{egalite_phi_R4_peiohbeofivn}
\varphi(A_\xi,B)=\abs{\det(M_\xi)}\frac{c_1}{H(B)}
\vspace{2.6mm}
\end{equation}

pour une certaine constante $c_1>0$ ne dépendant que de $A_\xi$. \\

Or d'après l'équation \eqref{valeur_detM_R4_zghodidvc}, on a

\[\abs{\det(M_\xi)}=\abs{-\eta_6+7\eta_1+(\eta_5-\eta_2)\xi+(-\eta_3-\eta_4)\sqrt{7-\xi^2}},\]

et il reste donc à minorer cette quantité en fonction de la hauteur de $B$. \\

Comme les coordonnées de Plücker $\eta=(\eta_1,\ldots,\eta_6)$ de $B$ sont entières et premières entre elles dans leur ensemble (équation \eqref{premiersentreeux_R4_idhvoevbd}), on peut utiliser la définition \ref{def_hauteur_coord_plucker_reuighgzmne} de la hauteur : 

\begin{equation}\label{R4_hauteur_fonction_de_eta_iahrnbodianv}
H(B)^2=\sum_{i=1}^6 \eta_i^2=\norme{\eta}^2,
\vspace{2.6mm}
\end{equation}

où $\norme\cdot$ est la norme euclidienne canonique sur $\R^6$, ce qui va permettre de minorer $\abs{\det(M_\xi)}$ en fonction de $H(B)$. \\

Soit $\epsilon>0$. Supposons désormais que $\xi$ soit un nombre algébrique tel que le triplet de nombres algébriques $(1,\xi,\sqrt{7-\xi^2})$ soit $\Q$-libre (\emph{i.e.} tel que la relation \eqref{indep_alg_R4_oersgbozdvbn} soit vérifiée). \\

Alors la proposition \ref{prop_approx_formelin_nombrealg_neorifmnbv}\label{app_pro_2_28_aepihboeifveifbn} s'applique, et donne une constante $c_2>0$ ne dépendant que de $A_\xi$ et de $\epsilon$, telle que pour tout $q=(a,b,c)\in\Z^3\setminus\{(0,0,0)\}$ :

\begin{equation}\label{minorationKB_R4_orubvdobv}
\abs{a\sqrt{7-\xi^2}+b\xi+c}\ge c_2\norme{q}^{-2-\epsilon}.
\vspace{2.6mm}
\end{equation}

On remarque que pour $q=(-\eta_3-\eta_4,\eta_5-\eta_2,-\eta_6+7\eta_1)$, on a $q\ne (0,0,0)$ sinon le système \eqref{cases_R4_peariohoeiv} serait vérifié, et on a vu que cela était impossible. De plus,
\begin{align*}
\norme q^2
	&=\abs{\eta_3^2+\eta_4^2+2\eta_3\eta_4+\eta_5^2+\eta_2^2-2\eta_2\eta_5+\eta_6^2+49\eta_1^2-14\eta_1\eta_6} \\
	&\le 49(\eta_1^2+\cdots+\eta_6^2)+18\cdot \max_{1\le i\le 6}(\abs{\eta_i})^2 \\
	&\le 67\norme{\eta}^2.
\end{align*}
La minoration \eqref{minorationKB_R4_orubvdobv} donne donc

\[\abs{\det(M_\xi)}\ge c_3\norme{\eta}^{-2-\epsilon}\]

pour une certaine constante $c_3>0$ qui ne dépend que de $A_\xi$ et de $\epsilon$. \\

En combinant cette minoration avec \eqref{egalite_phi_R4_peiohbeofivn} et \eqref{R4_hauteur_fonction_de_eta_iahrnbodianv}, on obtient une constante $c_4>0$ (qui ne dépend que de $A_\xi$ et de $\epsilon$) telle que

\[\varphi(A_\xi,B)\ge \frac{c_4}{H(B)^{3+\epsilon}},\]

ce qui termine la preuve du lemme \ref{condition_dirrationalite_4221_regouzhelfgu}.
\end{preuve}

\begin{remarque}

En appliquant le lemme \ref{minoration_psiphi_eomivocvbn} avec $j=2$, le lemme \ref{condition_dirrationalite_4221_regouzhelfgu} permet de montrer de la même façon que 

\[\muexp 4222\le \frac 32,\]

mais le théorème \ref{borne_sup_de_Schmidt_roimghnev} donne déjà que

\[\muexp 4222= 1.\]

\end{remarque}

\section{Dans $\R^5$}\label{section_R5_gaomfnbamobo}

Maintenant que le cas de $\R^4$ est complètement résolu, on s'intéresse au prochain cas encore ouvert : 

\[(n,d,e,j)=(5,3,2,1).\]

D'après les théorèmes \ref{deuxieme_borne_Schmidt_reuhfgezliuhl} et \ref{borne_sup_de_Schmidt_roimghnev}, on sait déjà que 

\[4\le\muexp 5321\le 7.\]

Sans déterminer complètement la valeur de $\muexp 5321$, on améliore cet encadrement dans le théorème suivant :

\begin{theoreme}\label{amelioration_cas_R5_eoufbemovb}

On a

\[\muexp 5321\le 6.\]

\subsection{Construction de sous-espaces mal approchés de $\mathfrak I_5(3,2)_1$}\label{ss_section_resultats_R5_amobibfaonvona}

Comme dans la section \ref{section_R4_maoefvbn} sur $\R^4$, construisons ici explicitement des sous-espaces de $\R^5$ de dimension $3$ vérifiant simultanément les deux conditions suivantes :
\begin{enumerate}[$\bullet$]
	
	\item être $(2,1)$-irrationnels,
	
	\item être mal approchés par les plans rationnels de $\R^5$. \\
	
\end{enumerate}

Soit $\zeta_3$ un nombre réel algébrique tel que

\[[\Q(\zeta_3):\Q]\ge 33.\]

Considérons les quatre nombres algébriques suivants :

\small\[\zeta_1=-\frac{112 \, \zeta_{3}^{4} - 196 \, \zeta_{3}^{3} - {\left(42 \, \sqrt{2} \zeta_{3}^{3} - 17 \, \sqrt{2} \zeta_{3}^{2} + 13 \, \sqrt{2} \zeta_{3}\right)} \sqrt{4 \, \zeta_{3} - 5} \sqrt{\zeta_{3} - 1} + 88 \, \zeta_{3}^{2} - 30 \, \zeta_{3} + 6}{4 \, {\left(10 \, \zeta_{3}^{4} - 7 \, \zeta_{3}^{3} - {\left(4 \, \sqrt{2} \zeta_{3}^{3} + 3 \, \sqrt{2} \zeta_{3}^{2} + \sqrt{2}\right)} \sqrt{4 \, \zeta_{3} - 5} \sqrt{\zeta_{3} - 1} - 10 \, \zeta_{3}^{2} + 5 \, \zeta_{3} - 2\right)}}
,\]

\[\hspace{-1cm}\zeta_2=-\frac{52 \, \zeta_{3}^{4} - 154 \, \zeta_{3}^{3} - {\left(18 \, \sqrt{2} \zeta_{3}^{3} - 35 \, \sqrt{2} \zeta_{3}^{2} + 13 \, \sqrt{2} \zeta_{3} - 6 \, \sqrt{2}\right)} \sqrt{4 \, \zeta_{3} - 5} \sqrt{\zeta_{3} - 1} + 148 \, \zeta_{3}^{2} - 60 \, \zeta_{3} + 18}{4 \, {\left(10 \, \zeta_{3}^{4} - 7 \, \zeta_{3}^{3} - {\left(4 \, \sqrt{2} \zeta_{3}^{3} + 3 \, \sqrt{2} \zeta_{3}^{2} + \sqrt{2}\right)} \sqrt{4 \, \zeta_{3} - 5} \sqrt{\zeta_{3} - 1} - 10 \, \zeta_{3}^{2} + 5 \, \zeta_{3} - 2\right)}},\]

\normalsize\[\zeta_4=-\frac{\sqrt{2} \sqrt{4 \, \zeta_{3} - 5} \sqrt{\zeta_{3} - 1} \zeta_{3}^{2} - 6 \, \zeta_{3}^{3} + 3 \, \zeta_{3}^{2} + 3 \, \zeta_{3}}{2 \, {\left(\zeta_{3}^{2} - 1\right)}},\]

\[\zeta_5=-\frac{\sqrt{2} \sqrt{4 \, \zeta_{3} - 5} \sqrt{\zeta_{3} - 1} \zeta_{3} - 3 \, \zeta_{3}^{2} + 3 \, \zeta_{3}}{2 \, {\left(\zeta_{3}^{2} - 1\right)}},\]

en supposant $\zeta_3\ge 5/4$ de sorte à ce que toutes les racines soient bien définies, et $\zeta_3$ suffisamment grand de sorte à ce que les dénominateurs ne soient pas nuls (ceux-ci ne s'annulent qu'en un nombre fini de valeurs). En étudiant les polynômes qui définissent les dénominateurs, on pourrait montrer que $\zeta_3\ge5/4$ suffit pour que les dénominateurs soient non nuls, mais cela ne sera pas utile ici. \\ 

On pose alors

\begin{equation}\label{relation_xi_zeta_R5_apihvdiifebnv}
\begin{cases}
\xi_1=1 \\
\xi_2=\zeta_2+\zeta_5 \\
\xi_3=-\zeta_1 \\
\xi_4=1+\zeta_1+\zeta_5 \\
\xi_5=\zeta_2 \\
\xi_6=2\zeta_2-\zeta_5 \\
\xi_7=-\zeta_3 \\
\xi_8=\zeta_3 \\
\xi_9=\zeta_4 \\
\xi_{10}=\zeta_5.
\end{cases}
\vspace{2.6mm}
\end{equation}

et

\begin{equation}\label{def_xi_R5_zovvozdibi}
\xi=(\xi_1,\ldots,\xi_{10}).
\vspace{2.6mm}
\end{equation}

Tout ceci permet de définir la famille de sous-espaces annoncée grâce au lemme suivant (qui sera démontré, ainsi que les lemmes \ref{dimensionrationnelle_amorigomegin} et \ref{lemme_crucial_R5_moqefivpzdin} et la proposition \ref{proposition_R5_mu_Axi_ozgbzeomfibnv}, dans la sous-section \ref{ss_section_preuves_R5_pbaifbapin}).

\begin{lemme}\label{existence_sev_R5_efubvodubcd}

Il existe un sous-espace vectoriel $A_\xi$ de $\R^5$ de dimension $3$ dont le vecteur $\xi$ est un représentant des coordonnées de Plücker.

\end{lemme}

Maintenant que les sous-espaces $A_\xi$ sont définis, il reste à montrer, comme dans la section précédente, qu'ils vérifient la condition de $(2,1)$-irrationalité et qu'ils sont mal approchés par les plans rationnels de $\R^5$. \\

Pour cela, on montre d'abord le lemme technique suivant. On rappelle que \hbox{$\zeta_3\ge 5/4$} a été défini comme un nombre réel suffisamment grand, et algébrique de degré supérieur à $33$, ce qu'on utilise ici.

\begin{lemme}\label{dimensionrationnelle_amorigomegin}

On a

\[\dim_{\Q}\vect_{\Q}(1,\zeta_1,\zeta_2,\zeta_3,\zeta_4,\zeta_5)=6.\]

\end{lemme}

Avec ce bagage, on peut alors énoncer des résultats qui ressemblent à ceux de la section précédente dans $\R^4$.

\begin{lemme}\label{lemme_crucial_R5_moqefivpzdin}

Pour tout $\xi$ défini précédemment en \eqref{def_xi_R5_zovvozdibi}, le sous-espace $A_\xi$ associé à $\xi$ par le lemme \ref{existence_sev_R5_efubvodubcd} vérifie

\begin{equation}\label{condition_irr_R5_eofivbemob}
A_\xi\in\mathfrak I_5(3,2)_1.
\vspace{2.6mm}
\end{equation}

De plus,

\[\forall \epsilon>0,\quad \exists c>0,\quad \forall B\in\mathfrak R_5(2),\quad \varphi(A_\xi,B)\ge\frac c{H(B)^{6+\epsilon}}.\]

\end{lemme}

\emph{A fortiori}, la condition de $(2,2)$-irrationalité étant plus faible, le lemme \ref{lemme_crucial_R5_moqefivpzdin} montre aussi que pour tout $\xi$, $A_\xi\in\mathfrak I_5(3,2)_2$. De plus, le lemme \ref{minoration_psiphi_eomivocvbn} (appliqué avec $j=1$) montre que le lemme \ref{lemme_crucial_R5_moqefivpzdin} est aussi vrai en remplaçant $\varphi(A_\xi,B)$ par $\psi_1(A_\xi,B)$. \\ 

La proposition suivante permet de démontrer le théorème \ref{amelioration_cas_R5_eoufbemovb}, et ainsi de s'approcher un peu plus d'une réponse au problème \ref{probleme_principal_gnzroign} dans $\R^5$.

\begin{proposition}\label{proposition_R5_mu_Axi_ozgbzeomfibnv}

Pour tout $\xi$ défini en \eqref{def_xi_R5_zovvozdibi}, le sous-espace $A_\xi$ associé à $\xi$ par le lemme \ref{existence_sev_R5_efubvodubcd} vérifie

\[\muexpA 5{A_\xi}21\le 6.\]

\end{proposition}

\end{theoreme}

\subsection{Les preuves}\label{ss_section_preuves_R5_pbaifbapin}

On démontre ici les résultats énoncés sans démonstration dans la sous-section \ref{ss_section_resultats_R5_amobibfaonvona}. Comme dans $\R^4$ (voir sous-section \ref{preuves_R4_oiagrhgmaeriofvn}), on déduit immédiatement la proposition \ref{proposition_R5_mu_Axi_ozgbzeomfibnv} et le théorème \ref{amelioration_cas_R5_eoufbemovb} du lemme \ref{lemme_crucial_R5_moqefivpzdin}. \\ 

Commençons par démontrer le lemme \ref{existence_sev_R5_efubvodubcd} montrant que le sous-espace $A_\xi$ est bien défini.

\begin{preuve}
On remarque que $\binom 5 3=10$. Il existe un sous-espace ayant $\xi$ pour coordonnées de Plücker si, et seulement si, les coordonnées de $\xi$ vérifient les relations de Plücker (données par le théorème \ref{relation_de_Plucker_oerighbfduvb}) pour un sous-espace de dimension $3$ de $\R^5$, \emph{i.e.}

\begin{equation}\label{relationPlucker5_mroegomgme}
\begin{cases}
	\xi_2\xi_5=\xi_3\xi_4+\xi_1\xi_6 \\
	\xi_2\xi_8=\xi_3\xi_7+\xi_1\xi_9 \\
	\xi_4\xi_8=\xi_5\xi_7+\xi_1\xi_{10} \\
	\xi_4\xi_9=\xi_6\xi_7+\xi_2\xi_{10} \\
	\xi_5\xi_9=\xi_6\xi_8+\xi_3\xi_{10}.
\end{cases}
\vspace{2.6mm}
\end{equation}

Explicitons d'où viennent les relations \eqref{relationPlucker5_mroegomgme}\label{explications_relations_plucker_R5_baonoevbovbn}. On reprend les notations du théorème \ref{relation_de_Plucker_oerighbfduvb} et de la sous-section \ref{coord_Pluck_teoghzleufve}. On note $(\xi_K)_{K\in\Lambda(3,5)}$ les coordonnées de Plücker de $B$. On a d'après le théorème \ref{relation_de_Plucker_oerighbfduvb} que pour tout $I\in\Lambda(2,5)$ et pour tout \hbox{$J=\{j_1<j_2<j_3<j_4\}\in\Lambda(4,5)$} : 

\[\sum_{\substack{1\le k\le 4\\j_k\notin I}} (-1)^{k+\iota(I,\{j_k\})} \xi_{I\cup\{j_k\}}\xi_{J\setminus\{j_k\}}=0.\]

Par exemple, pour $I=\{1,2\}$ et $J=\{2,3,4,5\}$, on obtient 

\[(-1)^{2+\iota(\{1,2\},\{3\})}\xi_{\{1,2,3\}}\xi_{\{2,4,5\}}+(-1)^{3+\iota(\{1,2\},\{4\})}\xi_{\{1,2,4\}}\xi_{\{2,3,5\}}+(-1)^{4+\iota(\{1,2\},\{5\})}\xi_{\{1,2,5\}}\xi_{\{2,3,4\}}=0,\]

soit

\[\xi_{\{1,2,3\}}\xi_{\{2,4,5\}}-\xi_{\{1,2,4\}}\xi_{\{2,3,5\}}+\xi_{\{1,2,5\}}\xi_{\{2,3,4\}}=0,\]

ce qui donne bien la deuxième équation de \eqref{relationPlucker5_mroegomgme} en faisant correspondre les $\xi_{\{i,j,k\}}$ et les $\xi_\ell$ avec l'ordre lexicographique sur $\Lambda(3,5)$ :
\begin{align*}
(\xi_{\{1,2,3\}},\xi_{\{1,2,4\}},\xi_{\{1,2,5\}},\xi_{\{1,3,4\}},\xi_{\{1,3,5\}},\xi_{\{1,4,5\}},\xi_{\{2,3,4\}},\xi_{\{2,3,5\}},\xi_{\{2,4,5\}},\xi_{\{3,4,5\}}) \\
=(\xi_1,\xi_2,\xi_3,\xi_4,\xi_5,\xi_6,\xi_7,\xi_8,\xi_9,\xi_{10}).
\end{align*}
Les quatre autres équations de \eqref{relationPlucker5_mroegomgme} s'obtiennent de façon similaire. \\

On vérifie que $\xi$ vérifie le système \eqref{relationPlucker5_mroegomgme} grâce au logiciel de calcul formel SageMath, avec les commandes suivantes (les variables \texttt{tmp\_zi} stockent les $\zeta_i$ comme définis en début de section) :

\begin{verbatim}
E_1=xi_2*xi_5-xi_3*xi_4-xi_1*xi_6;
E_1=calc(E_1.subs(xi_1=1).subs(xi_2=zeta_2+zeta_5).subs(xi_3=-zeta_1)
.subs(xi_4=1+zeta_1+zeta_5).subs(xi_5=zeta_2).subs(xi_6=2*zeta_2-zeta_5)
.subs(xi_7=-zeta_3).subs(xi_8=zeta_3).subs(xi_9=zeta_4).subs(xi_10=zeta_5)
.subs(zeta_1=tmp_z1).subs(zeta_2=tmp_z2).subs(zeta_4=tmp_z4));

E_2=xi_2*xi_8-xi_3*xi_7-xi_1*xi_9;
E_2=calc(E_2.subs(xi_1=1).subs(xi_2=zeta_2+zeta_5).subs(xi_3=-zeta_1)
.subs(xi_4=1+zeta_1+zeta_5).subs(xi_5=zeta_2).subs(xi_6=2*zeta_2-zeta_5).
subs(xi_7=-zeta_3).subs(xi_8=zeta_3).subs(xi_9=zeta_4).subs(xi_10=zeta_5)
.subs(zeta_1=tmp_z1).subs(zeta_2=tmp_z2).subs(zeta_4=tmp_z4));

E_3=xi_4*xi_8-xi_5*xi_7-xi_1*xi_10;
E_3=calc(E_3.subs(xi_1=1).subs(xi_2=zeta_2+zeta_5).subs(xi_3=-zeta_1)
.subs(xi_4=1+zeta_1+zeta_5).subs(xi_5=zeta_2).subs(xi_6=2*zeta_2-zeta_5)
.subs(xi_7=-zeta_3).subs(xi_8=zeta_3).subs(xi_9=zeta_4).subs(xi_10=zeta_5)
.subs(zeta_1=tmp_z1).subs(zeta_2=tmp_z2).subs(zeta_4=tmp_z4));

E_4=xi_4*xi_9-xi_6*xi_7-xi_2*xi_10;
E_4=calc(E_4.subs(xi_1=1).subs(xi_2=zeta_2+zeta_5).subs(xi_3=-zeta_1)
.subs(xi_4=1+zeta_1+zeta_5).subs(xi_5=zeta_2).subs(xi_6=2*zeta_2-zeta_5)
.subs(xi_7=-zeta_3).subs(xi_8=zeta_3).subs(xi_9=zeta_4).subs(xi_10=zeta_5)
.subs(zeta_1=tmp_z1).subs(zeta_2=tmp_z2).subs(zeta_4=tmp_z4));

E_5=xi_5*xi_9-xi_6*xi_8-xi_3*xi_10;
E_5=calc(E_5.subs(xi_1=1).subs(xi_2=zeta_2+zeta_5).subs(xi_3=-zeta_1)
.subs(xi_4=1+zeta_1+zeta_5).subs(xi_5=zeta_2).subs(xi_6=2*zeta_2-zeta_5)
.subs(xi_7=-zeta_3).subs(xi_8=zeta_3).subs(xi_9=zeta_4).subs(xi_10=zeta_5)
.subs(zeta_1=tmp_z1).subs(zeta_2=tmp_z2).subs(zeta_4=tmp_z4));

[E_1,E_2,E_3,E_4,E_5]
\end{verbatim}

qui renvoie 

\[\mathtt{[0,0,0,0,0]}\]

comme espéré.
\end{preuve}

Les sous-espaces $A_\xi$ étant bien définis, on démontre le lemme technique \ref{dimensionrationnelle_amorigomegin}.

\begin{preuve}
Soit $(a_0,\ldots,a_5)\in\Q^6$ tel que

\[\sum_{i=0}^5 a_i\zeta_i=0\]

où on a posé $\zeta_0:=1$. On veut montrer que les $a_i$ sont tous nuls. Pour cela, travaillons l'expression grâce au logiciel SageMath. Après avoir écrit une fonction \texttt{calc} pour aider SageMath à simplifier les expressions, on chasse les dénominateurs avec les lignes de calculs suivantes :

\begin{verbatim}
def calc(X):
    return (expand(numerator(expand(X)))/expand(denominator(expand(X))))
    					.canonicalize_radical()

D=calc(denominator(zeta_1)*denominator(zeta_2)*denominator(zeta_3)
					*denominator(zeta_4)*denominator(zeta_5));
D_1=calc(denominator(zeta_2)*denominator(zeta_3)*denominator(zeta_4)
					*denominator(zeta_5));
D_2=calc(denominator(zeta_1)*denominator(zeta_3)*denominator(zeta_4)
					*denominator(zeta_5));
D_3=calc(denominator(zeta_1)*denominator(zeta_2)*denominator(zeta_4)
					*denominator(zeta_5));
D_4=calc(denominator(zeta_1)*denominator(zeta_2)*denominator(zeta_3)
					*denominator(zeta_5));
D_5=calc(denominator(zeta_1)*denominator(zeta_2)*denominator(zeta_3)
					*denominator(zeta_4));

E=calc(a_0*D+a_1*(zeta_1)*D_1+a_2*(zeta_2)*D_2+a_3*(zeta_3)*D_3+a_4
*(zeta_4)*D_4+a_5*(zeta_5)*D_5);
E;
\end{verbatim}

Ceci donne 

\begin{equation}\label{sum_ai_zeta_i_agbeoibiobeiozbfc}
0=\sum_{i=0}^5 a_i\zeta_i=\frac{P_1(\zeta_3)+P_2(\zeta_3)\sqrt{P_3(\zeta_3)}}{10\zeta_3^3+7\zeta_3-2-(4\zeta_3^2-\zeta_3+1)\sqrt{P_3(\zeta_3)}}
\vspace{2.6mm}
\end{equation}

avec $P_1,P_2,P_3\in\Q[X]$ :
\begin{multline*}
P_1=77440 \, a_{3} X^{16} + 64 \, {\left(1210 \, a_{0} - 3127 \, a_{3}\right)} X^{15} - 32 \, {\left(6254 \, a_{0} + 3263 \, a_{3} - 5342 \, a_{4}\right)} X^{14}\\
- 16 \, {\left(6526 \, a_{0} - 40479 \, a_{3} + 35080 \, a_{4} - 7054 \, a_{5}\right)} X^{13} \\
 + 16 \, {\left(40479 \, a_{0} - 17471 \, a_{3} + 22640 \, a_{4} - 25699 \, a_{5}\right)} X^{12} \\
- 8 \, {\left(34942 \, a_{0} - 308 \, a_{1} - 137 \, a_{2} + 76708 \, a_{3} - 77312 \, a_{4} - 38023 \, a_{5}\right)} X^{11} \\
- 4 \, {\left(153416 \, a_{0} + 1942 \, a_{1} + 1246 \, a_{2} - 155896 \, a_{3} + 246351 \, a_{4} - 70711 \, a_{5}\right)} X^{10} \\
+ 2 \, {\left(311792 \, a_{0} + 2308 \, a_{1} + 3337 \, a_{2} - 8368 \, a_{3} + 194376 \, a_{4} - 377760 \, a_{5}\right)} X^{9} \\
 - 2 \, {\left(8368 \, a_{0} - 4403 \, a_{1} - 812 \, a_{2} + 158704 \, a_{3} - 101094 \, a_{4} - 256674 \, a_{5}\right)} X^{8} \\
 - 2 \, {\left(158704 \, a_{0} + 6646 \, a_{1} + 5797 \, a_{2} - 129152 \, a_{3} + 189576 \, a_{4} + 17808 \, a_{5}\right)} X^{7} \\
 + 2 \, {\left(129152 \, a_{0} + 2335 \, a_{1} + 4786 \, a_{2} - 27088 \, a_{3} + 122274 \, a_{4} - 121002 \, a_{5}\right)} X^{6} \\
 - 2 \, {\left(27088 \, a_{0} - 1556 \, a_{1} + 337 \, a_{2} + 25192 \, a_{3} + 37696 \, a_{4} - 122400 \, a_{5}\right)} X^{5} \\
 - 2 \, {\left(25192 \, a_{0} + 2349 \, a_{1} + 2088 \, a_{2} - 24168 \, a_{3} - 4210 \, a_{4} + 65986 \, a_{5}\right)} X^{4} \\
 + 2 \, {\left(24168 \, a_{0} + 1434 \, a_{1} + 1977 \, a_{2} - 11984 \, a_{3} + 3880 \, a_{4} + 24652 \, a_{5}\right)} X^{3} \\
 - 2 \, {\left(11984 \, a_{0} + 499 \, a_{1} + 970 \, a_{2} - 3088 \, a_{3} + 1468 \, a_{4} + 5516 \, a_{5}\right)} X^{2} \\
 + 8 \, {\left(772 \, a_{0} + 29 \, a_{1} + 68 \, a_{2} - 136 \, a_{3} + 102 \, a_{4} + 212 \, a_{5}\right)} X - 1088 \, a_{0} - 12 \, a_{1} - 96 \, a_{2},
\end{multline*}
\begin{multline*}
P_2=-27392 \, a_{3} X^{15} - 64 \, {\left(428 \, a_{0} - 623 \, a_{3}\right)} X^{14} + 32 \, {\left(1246 \, a_{0} + 2546 \, a_{3} - 1889 \, a_{4}\right)} X^{13} \\
+ 16 \, {\left(5092 \, a_{0} - 8569 \, a_{3} + 8149 \, a_{4} - 2494 \, a_{5}\right)} X^{12} \\
- 8 \, {\left(17138 \, a_{0} + 6778 \, a_{3} - 2257 \, a_{4} - 12560 \, a_{5}\right)} X^{11} \\
- 4 \, {\left(13556 \, a_{0} + 217 \, a_{1} + 97 \, a_{2} - 38960 \, a_{3} + 49453 \, a_{4} + 5626 \, a_{5}\right)} X^{10} \\
+ 4 \, {\left(38960 \, a_{0} + 445 \, a_{1} + 331 \, a_{2} - 11640 \, a_{3} + 31708 \, a_{4} - 31222 \, a_{5}\right)} X^{9} \\
- {\left(46560 \, a_{0} - 363 \, a_{1} + 873 \, a_{2} + 46752 \, a_{3} - 4664 \, a_{4} - 126940 \, a_{5}\right)} X^{8} \\
- 4 \, {\left(11688 \, a_{0} + 683 \, a_{1} + 386 \, a_{2} - 15000 \, a_{3} + 16585 \, a_{4} + 9823 \, a_{5}\right)} X^{7} \\
+ {\left(60000 \, a_{0} + 1611 \, a_{1} + 2367 \, a_{2} - 23360 \, a_{3} + 59236 \, a_{4} - 31312 \, a_{5}\right)} X^{6} \\
- 2 \, {\left(11680 \, a_{0} - 93 \, a_{1} + 369 \, a_{2} + 3392 \, a_{3} + 10088 \, a_{4} - 25172 \, a_{5}\right)} X^{5} \\
- {\left(6784 \, a_{0} + 845 \, a_{1} + 593 \, a_{2} - 9904 \, a_{3} - 4256 \, a_{4} + 30148 \, a_{5}\right)} X^{4} \\
+ 4 \, {\left(2476 \, a_{0} + 176 \, a_{1} + 203 \, a_{2} - 1540 \, a_{3} + 469 \, a_{4} + 3221 \, a_{5}\right)} X^{3} \\
- {\left(6160 \, a_{0} + 255 \, a_{1} + 483 \, a_{2} - 1600 \, a_{3} + 664 \, a_{4} + 3008 \, a_{5}\right)} X^{2} \\
+ 2 \, {\left(800 \, a_{0} + 31 \, a_{1} + 73 \, a_{2} - 176 \, a_{3} + 132 \, a_{4} + 268 \, a_{5}\right)} X - 352 \, a_{0} - 6 \, a_{1} - 30 \, a_{2},
\end{multline*}
\[P_3=8 \, X^{2} - 18 \, X + 10.\]

Ainsi, on a

\[P_1(\zeta_3)+P_2(\zeta_3)\sqrt{P_3(\zeta_3)}=0,\]

et le nombre $\zeta_3$ est racine du polynôme

\[P=P_1^2-P_2^2P_3\in\Q[X]\]

de degré au plus $32$. \\

Or $\zeta_3$ a été défini comme un nombre réel algébrique tel que

\[[\Q(\zeta_3):\Q]\ge 33.\]
 
Ainsi, comme les coefficients de $P$ sont rationnels et que $\zeta_3$ est racine de $P$, on a 

\[P=0.\]

Les $4$ équations données par les monômes de degrés $32$, $30$, $28$ et $26$ forment alors un système échelonné qui donne

\[a_0=a_3=a_4=a_5=0.\]

Considérons enfin l'équation donnée par le monôme de degré $22$ :

\[14a_1^2+4a_1a_2-a_2^2=0,\]

ce qui donne

\[a_1=a_2\frac{-2\pm 3\sqrt 2}{14}.\]

Comme $(a_1,a_2)\in\Q^2$, on a $a_1=a_2=0$, et finalement on a bien obtenu

\[\forall i\in\{0,\ldots,5\},\quad a_i=0,\]

soit

\[\dim_{\Q}\vect_{\Q}(1,\zeta_1,\zeta_2,\zeta_3,\zeta_4,\zeta_5)=6.\]
\end{preuve}

On peut désormais démontrer le lemme \ref{lemme_crucial_R5_moqefivpzdin}.

\begin{preuve}
Soient $B\in\mathfrak R_5(2)$ et $(Y_1,Y_2)$ une base de $B$ donnée par le lemme \ref{lemme_pour_coordPluck_premsentreelles_egouzbelvc}. Comme $\binom 52=10$, notons $(\eta_1,\ldots,\eta_{10})$ les coordonnées de Plücker de $B$ associées à la base $(Y_1,Y_2)$. D'après le lemme \ref{lemme_pour_coordPluck_premsentreelles_egouzbelvc}, on a $(\eta_1,\ldots,\eta_{10})\in\Z^{10}$ et 

\begin{equation}\label{coord_Pluck_premsentreelles_R5_zoigbemoubv}
\pgcd(\eta_1,\ldots,\eta_{10})=1.
\vspace{2.6mm}
\end{equation}

De plus, ce vecteur vérifie les relations de Plücker données par le théorème \ref{relation_de_Plucker_oerighbfduvb} (déjà explicitées page \pageref{explications_relations_plucker_R5_baonoevbovbn}, en dessous du système \eqref{relationPlucker5_mroegomgme}) :

\begin{equation}\label{relation_Pluck_eta_R5_aomefihv}
\begin{cases} \eta_2\eta_5=\eta_3\eta_4+\eta_1\eta_6 \\\eta_2\eta_8=\eta_3\eta_7+\eta_1\eta_9 \\ \eta_4\eta_8=\eta_5\eta_{7}+\eta_1\eta_{10} \\ \eta_4\eta_9=\eta_6\eta_{7}+\eta_2\eta_{10} \\ \eta_5\eta_9=\eta_6\eta_{8}+\eta_3\eta_{10}.\end{cases}
\vspace{2.6mm}
\end{equation}

Soit $(X_\xi^{(1)},X_\xi^{(2)},X_\xi^{(3)})$ une base de $A_\xi$ associée à $\xi$, le vecteur représentant de ses coordonnées de Plücker dont on s'est servi pour définir $A_\xi$. Notons $M_\xi$ la matrice de $\M_5(\R)$ dont les colonnes sont respectivement $X_\xi^{(1)},X_\xi^{(2)},X_\xi^{(3)},Y_1,Y_2$. \\

Remarquons que

\[A_\xi\cap B=\{0\}\iff \det M_\xi\ne \{0\}.\]

On calcule le déterminant de $M$ par un développement de Laplace par rapport aux trois premières colonnes (corollaire \ref{corollaire_dev_Laplace_amiovomvbnainv}), et on trouve

\[\det M_\xi=\xi_1\eta_{10}-\xi_2\eta_9+\xi_3\eta_8+\xi_4\eta_7-\xi_5\eta_6+\xi_6\eta_5-\xi_7\eta_4+\xi_8\eta_3-\xi_9\eta_2+\xi_{10}\eta_1.\]

Supposons par l'absurde que $\det M_\xi=0$, \emph{i.e.} supposons par l'absurde que \hbox{$A_\xi\cap B\ne \{0\}$}. On a alors en utilisant les relations \eqref{relation_xi_zeta_R5_apihvdiifebnv} entre les $\xi_i$ et les $\zeta_i$ que
\small\begin{align}
0
	&=\det (M_\xi)\nonumber \\
	&=\eta_{10}-(\zeta_2+\zeta_5)\eta_9-\zeta_1\eta_8+(1+\zeta_1+\zeta_5)\eta_7-\zeta_2\eta_6+(2\zeta_2-\zeta_5)\eta_5+\zeta_3\eta_4+\zeta_3\eta_3-\zeta_4\eta_2+\zeta_5\eta_1\nonumber \\
	&=\eta_{10}+\eta_7+(-\eta_8+\eta_7)\zeta_1+(-\eta_9-\eta_6+2\eta_5)\zeta_2+(\eta_4+\eta_3)\zeta_3-\eta_2\zeta_4+(-\eta_9+\eta_7-\eta_5+\eta_1)\zeta_5\label{det_Mxi_R5_zomibvdoibv}.
\end{align}\normalsize
Or d'après le lemme \ref{dimensionrationnelle_amorigomegin},

\[\dim_{\Q}\vect_{\Q}(1,\zeta_1,\zeta_2,\zeta_3,\zeta_4,\zeta_5)=6,\]

et les $\eta_i$ sont des nombres rationnels, donc

\begin{equation}\label{equationslineairesetai_morgnmeofnvef}
\begin{cases} \eta_{10}=-\eta_7\\\eta_8=\eta_7\\\eta_6=-\eta_9+2\eta_5\\\eta_4=-\eta_3\\\eta_2=0\\\eta_1=\eta_9-\eta_7+\eta_5.\end{cases}
\vspace{2.6mm}
\end{equation}

Avec les relations \eqref{equationslineairesetai_morgnmeofnvef} nouvellement établies, on réécrit les relations de Plücker pour les $\eta_i$, qui de \eqref{relation_Pluck_eta_R5_aomefihv} deviennent

\begin{equation}\label{eq_x2468_aqzsedr}
\begin{cases}
	\eta_3^2-2\eta_5^2+2\eta_5\eta_7-\eta_5\eta_9-\eta_7\eta_9+\eta_9^2=0 \\
	-\eta_3\eta_7-\eta_5\eta_9+\eta_7\eta_9-\eta_9^2=0 \\
	-\eta_3\eta_7-\eta_7^2+\eta_7\eta_9=0 \\
	-2\eta_5\eta_7-\eta_3\eta_9+\eta_7\eta_9=0 \\
	\eta_3\eta_7-2\eta_5\eta_7+\eta_5\eta_9+\eta_7\eta_9=0.
\end{cases}
\vspace{2.6mm}
\end{equation}

Il reste à montrer qu'il n'existe pas de solution rationnelle non nulle $(\eta_3,\eta_5,\eta_7,\eta_9)$ à ce système polynomial. \\

On appelle $I$ l'idéal de $\Q[\eta_3,\eta_5,\eta_7,\eta_9]$ engendré par les équations du système \eqref{eq_x2468_aqzsedr}.

Avec la théorie des idéaux d'élimination (à laquelle une introduction peut être trouvée page 392 de \cite{lang02}), on regarde l'idéal $J$ dans lequel on a éliminé $\eta_3$, puis $\eta_5$. Autrement dit,

\[J=I\cap \Q[\eta_3,\eta_5].\]

Pour cela on se sert encore une fois du logiciel de calcul formel SageMath, dans lequel les instructions suivantes sont saisies :

\begin{verbatim}
L_1=eta_3^2-2*eta_5^2+2*eta_5*eta_7-eta_5*eta_9-eta_7*eta_9+eta_9^2;
L_2=-eta_3*eta_7-eta_5*eta_9+eta_7*eta_9-eta_9^2;
L_3=-eta_3*eta_7-eta_7^2+eta_7*eta_9;
L_4=-2*eta_5*eta_7-eta_3*eta_9+eta_7*eta_9;
L_5=eta_3*eta_7-2*eta_5*eta_7+eta_5*eta_9+eta_7*eta_9;

R.<eta_3,eta_5,eta_7,eta_9>=QQ[];

I=R.ideal(L_1,L_2,L_3,L_4,L_5);

J=I.elimination_ideal(eta_3).elimination_ideal(eta_5).

J
\end{verbatim}

ce qui renvoie l'idéal

\[J=(2\eta_7^3-4\eta_7\eta_9^2+\eta_9^3).\]

Ainsi, toute solution de \eqref{eq_x2468_aqzsedr} se projette en une solution de 

\begin{equation}\label{equation_projetee_R5_apbionfobin}
2\eta_7^3-4\eta_7\eta_9^2+\eta_9^3=0.
\vspace{2.6mm}
\end{equation}

Si $\eta_9=0$, alors $\eta_7=0$. Sinon on peut poser $X=\eta_7/\eta_9$, et l'équation \eqref{equation_projetee_R5_apbionfobin} devient

\[P(X)=2X^3-4X+1=0.\]

Supposons qu'il existe $p/q\in\Q^*$ ($0$ n'est pas racine de $P$) une racine rationnelle non nulle de $P$ sous forme irréductible. On a alors

\[2p^3-4pq^2+q^3=0.\]

Comme $\pgcd(p,q)=1$, on a d'après le lemme de Gauss que $p\mid 1$ et $q\mid 2$. On vérifie alors que

\[P(-1),P(-1/2),P(1/2),P(1)\ne 0.\]

Ainsi, on a $\eta_7=\eta_9=0$. \\

La première équation de \eqref{eq_x2468_aqzsedr} devient alors

\[\eta_3^2=2\eta_5^2,\]

et donc $\eta_3=\eta_5=0$ car $\sqrt 2\notin\Q$. \\

Finalement, en utilisant les équations \eqref{equationslineairesetai_morgnmeofnvef}, on obtient

\[\forall i\in\{1,\ldots,10\},\quad \eta_i=0,\]

ce qui est absurde car les $\eta_i$ sont les coordonnées de Plücker d'un sous-espace vectoriel de $\R^5$ de dimension $2$. \\

Il n'existe donc pas de plan $B$ rationnel tel que $\det M_\xi\ne 0$, \emph{i.e.} tel que $A_\xi\cap B\ne \{0\}$. Cela montre que

\[A_\xi\in\mathfrak I_5(3,2)_1\]

et établit \eqref{condition_irr_R5_eofivbemob}, la première partie du lemme \ref{lemme_crucial_R5_moqefivpzdin}. \\ 

Comme la base $(Y_1,Y_2)$ de $B$ donnée par le lemme \ref{lemme_pour_coordPluck_premsentreelles_egouzbelvc} est aussi une $\Z$-base de $B\cap\Z^n$, on peut utiliser le lemme \ref{lien_proximite_hauteur_zozofvbcz} qui donne

\begin{equation}\label{egalite_phi_R5_erogmbodvb}
\varphi(A_\xi,B)=\abs{\det(M_\xi)}\frac{c_1}{H(B)}
\vspace{2.6mm}
\end{equation}

pour une certaine constante $c_1>0$ ne dépendant que de $A_\xi$. Or d'après l'équation \eqref{det_Mxi_R5_zomibvdoibv}, on a
\begin{align*}
\abs{\det (M_\xi)}\phantom{\eta_1-\eta_4+(-\eta_4-\eta_3)\zeta_1+(-2\eta_6+\eta_5-\eta_2)\zeta_2+(\eta_8-\eta_7)\zeta_3-\eta_9\zeta_4+(\eta_{10}+\eta_6-\eta_4-)} \\
	=\abs{\eta_1-\eta_4+(-\eta_4-\eta_3)\zeta_1+(-2\eta_6+\eta_5-\eta_2)\zeta_2+(\eta_8-\eta_7)\zeta_3-\eta_9\zeta_4+(\eta_{10}+\eta_6-\eta_4-\eta_2)\zeta_5},
\end{align*} 
et il reste donc à minorer cette quantité en fonction de la hauteur de $B$. \\

Comme les coordonnées de Plücker $\eta=(\eta_1,\ldots,\eta_{10})$ de $B$ sont entières et premières entre elles (équation \eqref{coord_Pluck_premsentreelles_R5_zoigbemoubv}), on peut utiliser la définition \ref{def_hauteur_coord_plucker_reuighgzmne} de la hauteur : 

\[H(B)^2=\sum_{i=1}^{10} \eta_i^2=\norme{\eta}^2,\]

où $\norme\cdot$ désigne la norme euclidienne canonique sur $\R^{10}$, ce qui va permettre de minorer $\abs{\det(M_\xi)}$ en fonction de la hauteur de $B$. \\

Soit $\epsilon>0$.

Comme d'après le lemme \ref{dimensionrationnelle_amorigomegin} on a

\[\dim_\Q\vect_\Q (1,\zeta_1,\zeta_2,\zeta_3,\zeta_4,\zeta_5)=6\]

et que ces $6$ nombres sont des réels algébriques, la proposition \ref{prop_approx_formelin_nombrealg_neorifmnbv} donne une constante $c_2>0$ ne dépendant que de $A_\xi$ et de $\epsilon$, telle que

\begin{equation}\label{minoration_fonction_hauteur_R5_apefoineovn}
\forall q=(a_0,\ldots,a_5)\in\Z^6\setminus\{(0,0,0,0,0,0)\},\quad \abs{\sum_{i=0}^5 a_i\zeta_i}\ge c_2\norme {q}^{-5-\epsilon}.
\vspace{2.6mm}
\end{equation}

Remarquons que pour 

\[q=(\eta_1-\eta_4,-\eta_4-\eta_3,-2\eta_6+\eta_5-\eta_2,\eta_8-\eta_7,-\eta_9,\eta_{10}+\eta_6-\eta_4-\eta_2),\]

on a $q\ne (0,0,0,0,0,0)$, sinon le système \eqref{equationslineairesetai_morgnmeofnvef} serait vérifié et on a déjà établi que cela était impossible. De plus,

\[\norme{q}^2\le 5(\eta_1^2+\cdots+\eta_{10}^2)+24\cdot\max_{1\le i\le 10}(\abs{\eta_i})^2 \le 29\norme{\eta}^2.\]

La minoration \eqref{minoration_fonction_hauteur_R5_apefoineovn} donne donc

\[\abs{\det(M_\xi)}\ge c_3 \norme{\eta}^{-5-\epsilon}\]

pour une certaine constante $c_3>0$ qui ne dépend que de $A_\xi$ et de $\epsilon$. \\

En combinant cette minoration avec \eqref{egalite_phi_R5_erogmbodvb}, on obtient une constante $c_4>0$ (qui ne dépend que de $A_\xi$ et de $\epsilon$) telle que 

\[\varphi(A_\xi,B)\ge \frac{c_4}{H(B)^{6+\epsilon}},\]

ce qui termine la preuve du lemme \ref{lemme_crucial_R5_moqefivpzdin}.
\end{preuve}

\begin{remarque}

En appliquant le lemme \ref{minoration_psiphi_eomivocvbn} avec $j=2$, le lemme \ref{lemme_crucial_R5_moqefivpzdin} permet de montrer de la même façon que

\[\muexp 6322\le \frac 62=3,\]

mais on sait déjà grâce au théorème \ref{borne_sup_de_Schmidt_roimghnev} que $\muexp 6322\le 5/4$.

\end{remarque}

\section{Commentaires généraux}\label{section_commentaires_généraux_amirbnamoi}

Il est probable que cette méthode s'adapte à d'autres cas particuliers et permette de diminuer certaines majorations de Schmidt. Cependant, dans chaque cas particulier des calculs de plus en plus lourds semblent nécessaires. Dans le cas de $\R^5$, toute la difficulté était de construire un sous-espace $A_\xi$ suffisamment compliqué pour que le système \eqref{eq_x2468_aqzsedr} n'ait pas de solution rationnelle non triviale -- ce qui implique que $A_\xi\in\mathfrak I_5(3,2)_1$ -- mais aussi suffisamment simple pour arriver à montrer effectivement que ce système n'a pas de solution rationnelle non triviale. \\

Résumons l'idée globale des preuves développées dans les sections \ref{section_R4_maoefvbn} et \ref{section_R5_gaomfnbamobo}. On cherche à construire un sous-espace $A$ qui vérifie les deux conditions suivantes qui s'opposent : 
\begin{enumerate}[$\bullet$]
	\item avoir \emph{beaucoup} de coordonnées de Plücker linéairement indépendantes sur $\Q$ pour que $A$ soit $(e,1)$-irrationnel ;
	\item avoir \emph{peu} de coordonnées de Plücker linéairement indépendantes sur $\Q$ pour obtenir le meilleur exposant possible grâce à la proposition \ref{prop_approx_formelin_nombrealg_neorifmnbv}.
\end{enumerate}
C'est tout l'équilibre entre ces deux conditions qui est au c\oe ur de cette méthode. \\

Pour montrer que le sous-espace $A$ est $(e,1)$-irrationnel, on suppose par l'absurde qu'il existe un sous-espace rationnel $B\in\mathfrak R_n(e)$ intersectant $A$ de façon non triviale. Le fait que les coordonnées de Plücker $\zeta$ de $B$ soient rationnelles, alors que les coordonnées de Plücker $\xi$ de $A$ ne le sont pas, permet d'obtenir des équations entre les $\zeta_i$. Plus $\dim_\Q\vect_\Q(\xi)$ est grand, plus on obtient un grand nombre d'équations entre les $\zeta_i$. On utilise alors les relations de Plücker entre les $\zeta_i$, qui forment un système $(S)$ d'équations non linéaires. Les $\zeta_i$ étant rationnels, des relations judicieuses entre eux peuvent permettre de montrer que $\zeta$ n'est pas solution de $(S)$. Ceci peut fonctionner même si le nombre de relations entre les $\zeta_i$ est petit devant le nombre d'équations de $(S)$, grâce à la non linéarité du système et au fait que les $\zeta_i$ soient rationnels.

\section{Application d'un résultat de Moshchevitin}\label{section_app_Moshchevitin_aoimbnoibav}

Dans cette section, on déduit le résultat suivant du théorème \ref{th_Moshchevitin_2020_general_ergoighvoifvn} de Moshchevitin, grâce au théorème du Going-up (théorème \ref{goingup_eorihgzmefvnfon}).

\begin{theoreme}\label{theoreme_avec_moshchevitin_goingup_aeronfvn}

Si $n=2d$ avec $d\ge 2$, alors

\[\muexp nd{d-1}1\le \frac{2d^2}{d+1}.\]

\end{theoreme}

\begin{preuve}
Posons $n=2d$ et soit $\epsilon>0$. D'après le théorème \ref{th_Moshchevitin_2020_general_ergoighvoifvn} de Moshchevitin appliqué avec $\omega(t)=t^{-n-\epsilon}$, il existe un sous-espace vectoriel $A_\epsilon\in\mathfrak I_n(d,d)_1$ et une constante $c_1>0$ dépendant uniquement de $\epsilon$, tels que

\begin{equation}\label{preuve6321_aoregubelufb}
\forall D\in\mathfrak R_n(d),\quad\psi_1(A_\epsilon,D)\ge \frac {c_1}{H(D)^{n+\epsilon}}.
\vspace{2.6mm}
\end{equation}

On a $A_\epsilon\in\mathfrak I_n(d,d)_1\subset \mathfrak I_n(d,d-1)_1$, notons $\alpha=\muexpA n{A_\epsilon}{d-1}1$. Il existe une infinité de sous-espaces rationnels \hbox{$B\in\mathfrak R_n(d-1)$} tels que 

\begin{equation}\label{ineg_psi_1_A_B_arpeibnapifvndpiva}
\psi_1(A_\epsilon,B)\le \frac 1{H(B)^{\alpha-\epsilon}}.
\vspace{2.6mm}
\end{equation}

Le théorème \ref{goingup_eorihgzmefvnfon} du Going-up appliqué à $A_\epsilon$ et $B$ avec

\[\begin{cases} H=H(B)\\x_1=\alpha-\epsilon\\y_1=0\\\forall i\in\{2,\ldots,d-1\},\quad x_i=y_i=0\end{cases}\]

montre l'existence de $C\in\mathfrak R_n(d)$ contenant $B$, tel que

\begin{equation}\label{preuve6321_erjgmoiefhgnoeinv}
\psi_1(A_\epsilon,C)\le \frac {c_2}{H(C)^{(\alpha-\epsilon)(n-(d-1))/(n-d)}}=\frac {c_2}{H(C)^{(\alpha-\epsilon)(d+1)/d}}
\vspace{2.6mm}
\end{equation}

avec $c_2>0$ dépendant uniquement de $\epsilon$ (on rappelle que $n=2d$). \\

On a $A_\epsilon\in\mathfrak I_n(d,d)_1$, donc pour tout $C\in\mathfrak R_n(d)$, on a $\psi_1(A_\epsilon,C)\ne 0$. Ainsi, s'il n'existait qu'un nombre fini de sous-espaces rationnels $C$ vérifiant l'inégalité \eqref{preuve6321_erjgmoiefhgnoeinv}, on aurait une constante $c_3>0$ telle que 

\begin{equation}\label{minoration_psi1_A_eps_C_apihnramoegbobb}
\forall C\in\mathfrak R_n(d),\quad \psi_1(A_\epsilon,C)>c_3.
\vspace{2.6mm}
\end{equation}

Or il existe une infinité de sous-espaces $B\in\mathfrak R_n(d-1)$ vérifiant l'inégalité \eqref{ineg_psi_1_A_B_arpeibnapifvndpiva}, donc il en existe de hauteur arbitrairement grande, si bien que $\psi_1(A_\epsilon,B)\le c_3$. D'après le corollaire \ref{proximite_et_inclusions_aeoibnfvoidsn}, le sous-espace $C$ obtenu à partir d'un tel $B$ vérifie \eqref{preuve6321_erjgmoiefhgnoeinv} avec

\[\psi_1(A_\epsilon,C)\le \psi_1(A_\epsilon, B)\le c_3,\]

ce qui contredit \eqref{minoration_psi1_A_eps_C_apihnramoegbobb}. Il existe donc une infinité de sous-espaces $C\in\mathfrak R_n(d)$ pour lesquels \eqref{preuve6321_erjgmoiefhgnoeinv} est vérifiée. \\

La hauteur de $C$ pouvant être arbitrairement grande, et les constantes $c_1$ et $c_2$ dépendant uniquement de $\epsilon$, les inégalités \eqref{preuve6321_aoregubelufb} et \eqref{preuve6321_erjgmoiefhgnoeinv} imposent

\[(\alpha-\epsilon)\frac {d+1}d\le n+\epsilon.\]

Finalement,

\[\muexp nd{d-1}1\le \muexpA n{A_\epsilon}{d-1}1=\alpha\le \frac{2d^2}{d+1}+\epsilon\left(\frac{d}{d+1}+1\right),\]

et ceci étant vrai pour tout $\epsilon>0$, on obtient

\[\muexp nd{d-1}1\le \frac{2d^2}{d+1}.\]
\end{preuve}

\begin{corollaire}\label{cor_6321_eorglhleozfubv}

On a 

\[\begin{cases}\muexp 6321\le 9/2=4.5\\\muexp 8431\le 32/5=6.4.\end{cases}\]

\end{corollaire}

\begin{remarque}

Le corollaire \ref{cor_6321_eorglhleozfubv} améliore les majorants de Schmidt (théorème \ref{borne_sup_de_Schmidt_roimghnev}) :

\[\begin{cases} \muexp 6321\le 5\\\muexp 8431\le 8.\end{cases}\]

La nouvelle majoration est même asymptotiquement meilleure que celle de Schmidt. En effet, on a obtenu un majorant équivalent quand $n$ tend vers l'infini à $2d$, alors que le théorème \ref{borne_sup_de_Schmidt_roimghnev} donne

\[\muexp nd{d-1}1\le \left\lfloor\frac{d^2}2\right\rfloor.\]

\end{remarque}

\chapter{Approximation de sommes directes de sous-espaces vectoriels}\label{chap_somme_directes_vaoinvoin}

Dans ce chapitre, l'objectif est de minorer $\muexp ndej$ pour tout $n\ge 4$. La stratégie consiste à décomposer le sous-espace vectoriel $A$ approché en une somme directe de sous-espaces de dimensions plus petites (en l'occurence des droites). Il est alors possible d'approcher simultanément ces droites (grâce au théorème d'approximation simultanée de Dirichlet), et d'en déduire une approximation de $A$. \\

Ceci permet d'améliorer la minoration connue de $\muexp ndej$ : c'est le théorème \ref{premiereborneobtenue_amorimeovbn}, qui donne 

\[\muexp ndej\ge \frac{(n-j)(jn-jd+j^2/2+j/2+1)}{j^2(n-e)(n-d+j/2+1/2)}.\]

\section{Les résultats obtenus}\label{section_resultats_chap_4_aoibnbaionvaoin}

\subsection{Un résultat d'approximation}

Le résultat principal de ce chapitre est le théorème suivant. Sa preuve consiste à construire un sous-espace vectoriel de dimension $j$ de l'espace $A$ qu'on cherche à approcher. On décompose ce sous-espace en somme directe de droites. On approche ensuite simultanément ces droites par des droites rationnelles de hauteurs comparables entre elles, grâce au théorème d'approximation simultanée de Dirichlet. On conclut enfin grâce au théorème du Going-up.

\begin{theoreme}\label{premiereborneobtenue_amorimeovbn}

Soient $n\ge 4$ et $d,e\in\{1,\ldots,n-1\}$ tels que $d+e\le n$ ; soit $j\in\{1,\ldots,\min(d,e)\}$. On a

\[\muexp ndej\ge \frac{(n-j)(jn-jd+j^2/2+j/2+1)}{j^2(n-e)(n-d+j/2+1/2)}.\]

\end{theoreme}

\begin{remarque}

La minoration de $\muexp ndej$ obtenue dans le théorème \ref{premiereborneobtenue_amorimeovbn} est équivalente quand $n$ tend vers l'infini (avec $j$, $d$ et $e$ fixés) à $1/j$. Celle-ci est asymptotiquement meilleure que celle de Schmidt (théorème \ref{premiere_borne_Schmidt_hgoembeer}) qui est équivalente quand $n$ tend vers l'infini à 

\[\frac d{jn}\xrightarrow[n\to+\infty]{} 0.\]

Le théorème \ref{premiere_borne_Schmidt_hgoembeer} donne par exemple que pour tout $n\ge 5$, 

\[\muexp n222\ge \frac 1{n-2},\]

ce qui est moins bon que la minoration

\[\muexp n 222\ge \frac n{2n-1}\]

donnée par le théorème \ref{premiereborneobtenue_amorimeovbn}. 

\end{remarque}

On en déduit immédiatement du théorème \ref{premiereborneobtenue_amorimeovbn} le corollaire suivant.

\begin{corollaire}\label{corollaire_aeroibfnoiefbivo}

Soient $n\ge 4$ et $d\le n/2$ un entier. On a

\[\muexp nddd\ge\frac{2dn-d^2+d+2}{2d^2n-d^3+d^2}.\]

\end{corollaire}

De plus, en combinant le corollaire \ref{corollaire_aeroibfnoiefbivo} et le corollaire \ref{cor_maj_de_Saxce_amovrvrhbeofin} de Saxcé, on obtient pour $n\ge 2d$ l'encadrement

\[\frac{2dn-d^2+d+2}{2d^2n-d^3+d^2}\le \muexp nddd\le \frac n{d(n-d)},\]

ce qui donne le corollaire suivant.

\begin{corollaire}\label{cor_limite_avec_maj_desaxce_amoberubfvs}

On a 

\[\lim_{n\to+\infty} \muexp nddd = \frac 1d.\]

\end{corollaire}

\subsection{Comportement de la proximité par sommes directes}\label{ss_section_resultats_reconstruction_agmoibeamoiv}

L'idée principale de la preuve du théorème \ref{premiereborneobtenue_amorimeovbn} est de décomposer un sous-espace en somme directe de sous-espaces plus petits. Il est donc utile de savoir comment se comporte la proximité par l'opération de somme directe. 

\begin{proposition}\label{resultat_reconstruction_proximite_oaeirgbvodbv}

Soient $n\ge 4$, et $F_1,\ldots,F_\ell, B_1,\ldots,B_\ell$, $2\ell$ sous-espaces vectoriels de $\R^n$ tels que pour tout $i\in\{1,\ldots,\ell\}$, $\dim F_i=\dim B_i=d_i$. Supposons que les $F_i$ engendrent un sous-espace de dimension $k=d_1+\cdots+d_\ell$ et de même pour les $B_i$. Posons 

\[F=\bigoplus_{i=1}^\ell F_i \qquad\text{ et }\qquad B=\bigoplus_{i=1}^\ell B_i.\]

Alors on a

\[\psi_k(F,B)\le c_{F,n}\sum_{i=1}^\ell\psi_{d_i}(F_i,B_i)\]

où $c_{F,n}>0$ est une constante qui dépend uniquement de $F_1,\ldots,F_\ell$ et de $n$.

\end{proposition}

\begin{remarque}

Intuitivement, cela signifie que si les sous-espaces $F_i$ et $B_i$ sont presque confondus pour tout $i$, alors $F$ et $B$ sont presque confondus.

\end{remarque}

\section{Les preuves}

Démontrons ici les résultats de la section \ref{section_resultats_chap_4_aoibnbaionvaoin}.

\subsection{Proximité et sommes directes}

On démontre la proposition \ref{resultat_reconstruction_proximite_oaeirgbvodbv}.

\begin{preuve}
L'idée de la preuve est de décomposer chaque $F_i$ et chaque $B_i$ en somme directe de droites bien choisies. On commence par un bref lemme qui indiquera comment bien choisir ces droites. 

\begin{lemme}\label{minorationinverseavecdroitesbienchoisies_paeivmsoivn}

Soient $D$ et $E$ des sous-espaces vectoriels de $\R^n$ de même dimension $k$. Il existe $D_1,\ldots,D_k$ des droites de $D$ et $E_1,\ldots,E_k$ des droites de $E$, telles que 

\[D=\bigoplus_{i=1}^k D_i \quad \text { et }\quad E=\bigoplus_{i=1}^k E_i,\]

et qui vérifient de plus :

\[\psi_k(D,E)\le\sum_{i=1}^k \psi_1(D_i,E_i)\le k\psi_k(D,E).\]

\end{lemme}

\begin{preuve}[Lemme \ref{minorationinverseavecdroitesbienchoisies_paeivmsoivn}]
D'après le théorème \ref{la_raison_pour_laquelle_les_angles_sont_canoniques}, il existe des bases orthonormales $(X_1,\ldots,X_k)$ de $D$ et $(Y_1,\ldots,Y_k)$ de $E$ telles que

\[\forall i\in\{1,\ldots,k\},\quad \psi_i(D,E)=\psi(X_i,Y_i).\]

Or, toujours d'après ce théorème, on a 

\[\forall i\in\{1,\ldots,k\},\quad \psi_i(D,E)\le \psi_k(D,E).\]

Ainsi, on peut poser pour tout $i\in\{1,\ldots,k\}$

\[D_i=\vect(X_i)\quad\text { et }\quad E_i=\vect(Y_i).\]

On a alors

\[\sum_{i=1}^k\psi_1(D_i,E_i)=\sum_{i=1}^k \psi(X_i,Y_i)=\sum_{i=1}^k \psi_i(D,E)\le k\psi_k(D,E)\]

ce qui conclut la preuve de la majoration. \\

Pour la minoration, on remarque que tous les $\psi_1(D_i,E_i)$ sont positifs, et donc 

\[\psi_k(D,E)=\psi_1(D_k,E_k)\le \sum_{i=1}^k \psi_1(D_i,E_i).\]
\end{preuve}

Revenons à la preuve de la proposition \ref{resultat_reconstruction_proximite_oaeirgbvodbv}. \\

Soit $i\in\{1,\ldots,\ell\}$. D'après le lemme \ref{minorationinverseavecdroitesbienchoisies_paeivmsoivn}, il existe des droites $D_{i,1},\ldots,D_{i,d_i}$ de $F_i$ et des droites $E_{i,1},\ldots,E_{i,d_i}$ de $B_i$ telles que

\begin{equation}\label{majorationinvgraceaulemme}
\sum_{j=1}^{d_i} \psi_1(E_{i,j},D_{i,j})\le d_i\psi_{d_i}(F_i,B_i) \le n\psi_{d_i}(F_i,B_i).
\vspace{2.6mm}
\end{equation}

Soient $e_{i,1},\ldots,e_{i,d_i}$ des vecteurs directeurs unitaires de $D_{i,1},\ldots,D_{i,d_i}$ respectivement et $b_{i,1},\ldots,b_{i,d_i}$ des vecteurs directeurs unitaires de $E_{i,1},\ldots,E_{i,d_i}$ respectivement, tels que

\begin{equation}\label{eij_bij_prod_scal_pos_amoeriahgoudbv}
\forall j\in\{1,\ldots,d_i\},\quad e_{i,j}\cdot b_{i,j}\ge 0.
\vspace{2.6mm}
\end{equation}

D'après le théorème \ref{la_raison_pour_laquelle_les_angles_sont_canoniques}, il existe des vecteurs unitaires $X_k\in F$ et $Y_k\in B$ tels que 

\[\psi_k(F,B)=\psi(X_k,Y_k).\]

Soit $Z=\lambda_1 Y_1+\cdots+\lambda_k Y_k$ un vecteur unitaire de $B$. D'après le théorème \ref{la_raison_pour_laquelle_les_angles_sont_canoniques}, on a

\begin{equation}\label{majoration_k_eme_angle_ganoimfabhoivbsoi}
\abs{X_k\cdot Z}=\abs{\sum_{i=1}^k \lambda_i X_k\cdot Y_i}\le\sum_{i=1}^k \abs{\lambda_i \delta_{i,k} X_k\cdot Y_i} =\abs{\lambda_k}X_k\cdot Y_k\le X_k\cdot Y_k
\vspace{2.6mm}
\end{equation}

car $Z$ étant unitaire, $\abs{\lambda_k}\le 1$. L'inégalité \eqref{majoration_k_eme_angle_ganoimfabhoivbsoi} donne donc 
\[\psi_k(F,B)=\psi(X_k,Y_k)\le \min_{Z\in B\setminus\{0\}} \psi(X_k,Z)=\psi_1(\vect(X_k),B).\]

De plus, $\vect(Y_k)\subset B$, donc d'après le corollaire \ref{proximite_et_inclusions_aeoibnfvoidsn}, on a

\[\psi_k(F,B)=\psi(X_k,Y_k)\ge \psi_1(\vect(X_k),B).\]

Finalement,

\begin{equation}\label{egalite_psi_k_aroeihnveoivn}
\psi_k(F,B)=\psi_1(\vect(X_k),B).
\vspace{2.6mm}
\end{equation}

Décomposons $X_k$ dans la base $(e_{1,1},\ldots,e_{\ell,d_\ell})$ : 

\[X_k=\sum_{i=1}^\ell\sum_{j=1}^{d_i} x_{i,j} e_{i,j},\]

et posons

\[Y=\sum_{i=1}^\ell\sum_{j=1}^{d_i} x_{i,j}b_{i,j}\in B.\]

Comme $X_k$ est unitaire, on a d'après le lemme \ref{lemme_maj_psiXY_amoignaoivnovibs} :

\[\psi(X_k,Y)\le \norme{X_k-Y}=\norme{\sum_{i=1}^\ell\sum_{j=1}^{d_i} x_{i,j} (e_{i,j}-b_{i,j})}\le \sum_{i=1}^\ell\sum_{j=1}^{d_i} \abs{x_{i,j}} \norme{e_{i,j}-b_{i,j}},\]

où $\norme\cdot$ désigne -- ici et dans toute la suite -- la norme euclidienne canonique. \\

Considérons maintenant pour $i\in\{1,\ldots,\ell\}$ et $j\in\{1,\ldots,d_i\}$ les applications coordonnées : 

\[\fonction{p_{i,j}} F \R {\displaystyle\sum_{i=1}^\ell\sum_{j=1}^{d_i} x_{i,j}e_{i,j}}{x_{i,j}.}\]

Ces applications sont continues sur le compact

\[K=\{x\in F,\ \norme x=1\},\]

elles sont donc bornées sur cet ensemble. Ainsi, il existe $c_{F,n}^{(1)}$ une constante dépendant seulement de $e_{1,1},\ldots,e_{\ell,d_\ell}$ telle que

\[\forall x=\sum_{i=1}^\ell\sum_{j=1}^{d_i} x_{i,j} e_{i,j}\in K,\quad \abs{x_{i,j}}\le c_{F,n}^{(1)}.\]

Revenons au calcul principal : comme pour tous $i,j$ on a $e_{i,j}\cdot b_{i,j}\ge 0$ d'après \eqref{eij_bij_prod_scal_pos_amoeriahgoudbv}, on a en appliquant le lemme \ref{lemme_minoration_psiXY_amoifbvoisbdvoibs} que

\[\psi(X_k,Y)\le c_{F,n}^{(1)} \sum_{i=1}^\ell\sum_{j=1}^{d_i} \norme{e_{i,j}-b_{i,j}}\le c_{F,n}^{(2)} \sum_{i=1}^\ell\sum_{j=1}^{d_i} \psi_1(D_{i,j},E_{i,j})\]

car les $e_{i,j}$ et les $b_{i,j}$ sont des vecteurs unitaires, avec $c_{F,n}^{(2)}=\sqrt 2c_{F,n}^{(1)}$. \\

Finalement, avec l'inégalité \eqref{majorationinvgraceaulemme}, on a

\begin{equation}\label{majoration_psi_preuve_prox_amoefivbv}
\psi(X_k,Y)\le c_{F,n}^{(2)}n \sum_{i=1}^\ell \psi_{d_i}(F_i,B_i).
\vspace{2.6mm}
\end{equation}

Or, l'équation \eqref{egalite_psi_k_aroeihnveoivn} donne 

\[\psi_k(F,B)\le\psi_1(\vect(X_k),B)\le \psi(X_k,Y),\]

car $Y\in B$, et avec l'inégalité \eqref{majoration_psi_preuve_prox_amoefivbv} on obtient

\[\psi_k(F,B)\le c_{F,n} \sum_{i=1}^\ell \psi_{d_i}(F_i,B_i),\]

ce qui termine la preuve de la proposition \ref{resultat_reconstruction_proximite_oaeirgbvodbv}.
\end{preuve}

\subsection{Minoration de $\muexp ndej$}

Démontrons le théorème \ref{premiereborneobtenue_amorimeovbn}.

\begin{preuve}
Soit $F\in\mathfrak I_n(d,e)_j$ un sous-espace $(e,j)$-irrationnel de $\R^n$. \\

Montrons que $F$ possède une famille orthonormale $(f_1,\ldots,f_j)$ telle que pour tout $\ell\in\{1,\ldots,j\}$, $f_\ell$ ait au moins $d-\ell$ coordonnées nulles. \\

Pour cela on construit les $f_\ell$ par récurrence sur $\ell\in\{0,\ldots,j\}$. Pour $\ell=0$, il n'y a rien à construire et le résultat est évident. \\

On suppose désormais que les vecteurs $f_1,\ldots,f_\ell$ sont construits pour un certain \hbox{$\ell\in\{0,\ldots,j-1\}$}. On note $G$ le supplémentaire orthogonal de $\vect(f_1,\ldots,f_\ell)$ dans $F$. On a

\[\codim(\R^{n-d+\ell+1}\times\{0\}^{d-\ell-1})=d-\ell-1<d-\ell=\dim F-\ell=\dim G-1,\]

donc

\[G\cap(\R^{n-d+\ell+1}\times\{0\}^{d-\ell-1})\ne \{0\}.\]

On choisit $f_{\ell+1}\in G\cap(\R^{n-d+\ell+1}\times\{0\}^{d-\ell-1})$ unitaire. Le vecteur $f_{\ell+1}$ a bien au moins $d-(\ell+1)$ coordonnées nulles, et il est orthogonal à $f_1,\ldots,f_\ell$, ce qui termine la récurrence. \\

On fixe donc dans toute la suite une famille orthonormale $(f_1,\ldots,f_j)$ de $F$ telle que pour tout $\ell\in\{1,\ldots,j\}$, le vecteur $f_\ell$ ait au moins $d-\ell$ coordonnées nulles. \\

Notons $\underline x$ le vecteur formé de toutes les coordonnées non nulles des $f_\ell$ pour \hbox{$\ell\in\{1,\ldots,j\}$}. Comme chaque $f_\ell$ a au moins $d-\ell$ coordonnées nulles pour $\ell\in\{1,\ldots,j\}$, $\underline x$ a au plus 

\[\sum_{\ell=1}^j (n-(d-\ell))=jn-jd+\frac 12j^2+\frac 12 j\]

coordonnées. Notons $N\in\{1,\ldots,jn-jd+j^2/2+j/2\}$ le nombre de coordonnées de $\underline x$. \\

On a $\underline x\in\R^N\setminus\Q^N$, sinon $(f_1,\ldots,f_j)$ serait une base rationnelle d'un sous-espace vectoriel de dimension $j$ de $F$, ce qui contredirait l'hypothèse $F\in\mathfrak I_n(d,e)_j$. \\

Donc d'après le théorème d'approximation simultanée de Dirichlet (théorème \ref{approx_simul_dirichlet_aemofibbv}), il existe une infinité de couples $(p,q)\in\Z^N\times \N^*$ tels que $\pgcd(p_1,\ldots,p_N,q)=1$ et 

\begin{equation}\label{approxsimultaneedeDir}
\norme{\underline x-\frac pq}_{\infty}\le \frac 1{q^{1+1/N}}.
\vspace{2.6mm}
\end{equation}

On se donne un tel couple $(p,q)$. \\

Rappelons que $\underline x$ est le vecteur formé par les coordonnées non nulles de $f_1,\ldots,f_j$. Pour $i\in\{1,\ldots,j\}$, notons $p_i$ le sous-vecteur de $p$ correspondant aux coordonnées approchant celles de $f_i$, complété avec des zéros de sorte que $p_i\in\Z^n$ soit proche de $qf_i$. On a

\[\forall i\in\{1,\ldots,j\},\quad \norme{f_i-\frac{p_i}q}_{\infty}\le \frac 1{q^{1+1/N}}.\]

Posons

\[B=\vect(p_1,\ldots,p_j),\]

et notons $p_\perp(f_i)$ le projeté orthogonal de $f_i$ sur $\vect(p_i/q)$. On peut alors mener le calcul suivant, illustré sur la figure \ref{calculapproxsimultanee} (voir aussi le lemme \ref{lemme_proximite_proj_ortho_apamoebnamon}) : 

\begin{equation}\label{majoration_psi_moaeirbsovmin}
\psi(f_i,p_i/q)=\sin\widehat{\left(f_i, p_i/q\right)}=\frac{\norme{f_i-p_\perp(f_i)}}{\norme{f_i}}\le\norme{f_i-\frac{p_i}q}\le \frac{c_1}{q^{1+1/N}}
\vspace{2.6mm}
\end{equation}

car $\norme{f_i}= 1$, avec $c_1>0$ ne dépendant que de $n$. \\

\begin{figure}[H]
\begin{center}
\includegraphics[scale=1]{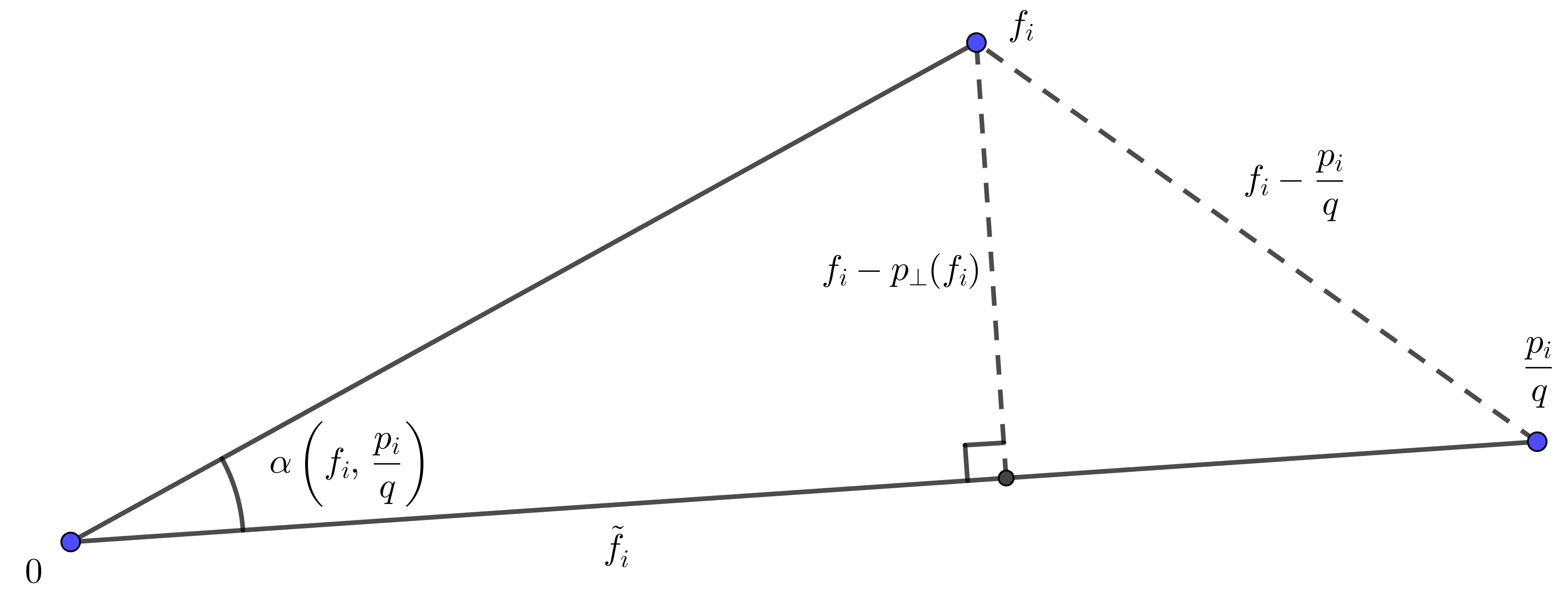}
\caption{Approximation de $f_i$ par $p_i/q$}
\label{calculapproxsimultanee}
\end{center}
\end{figure}

Il reste alors à majorer la hauteur de $B=\vect(p_1,\ldots,p_j)$. \\

L'inégalité \eqref{approxsimultaneedeDir} donne 
\[\norme{p}_\infty-\norme{q\underline x}_\infty\le \norme{q\underline x-p}_\infty\le q^{-1/N}\le 1,\]

donc pour tout $i\in\{1,\ldots,j\}$ :

\[\norme{p_i}_\infty\le \norme{p}_\infty\le 1+\norme{q\underline x}_\infty\le c_2q\]

avec $c_2>0$ qui ne dépend que de $F$. \\

Pour $E$ sous-espace vectoriel de $\R^n$ et $P$ une famille de vecteurs linéairement indépendants de $E$, notons $\vol_E(P)$ le volume du parallélotope engendré par les vecteurs de $P$ vus dans l'espace euclidien $E$. Comme $(p_1,\ldots,p_j)$ est un sous-réseau de $B\cap\Z^n$, on a d'après la proposition \ref{defequiv_hauteur_aveclecovol_eohuflbg} : 

\[H(B)\le\vol_B(p_1,\ldots,p_j)\le \prod_{i=1}^j \norme{p_i}\le c_3 q^j\]

avec $c_3>0$ qui ne dépend que de $F$. Ainsi, il existe une constante $c_4>0$ telle que 

\begin{equation}\label{majoration_1q_aprighmoidvn}
\frac 1q\le \frac{c_4}{H(B)^{1/j}}.
\vspace{2.6mm}
\end{equation}

Rappelons qu'en posant $B_i=\vect(p_i)$, la relation \eqref{majoration_psi_moaeirbsovmin} donne

\[\forall i\in\{1,\ldots,j\},\quad \psi_1(\vect(f_i),B_i)\le \frac{c_1}{q^{(N+1)/N}}.\]

Notons

\[\tilde F_j=\vect(f_1,\ldots,f_j)\]

qui est de dimension $j$ car la famille $(f_1,\ldots,f_j)$ est libre. Alors d'après la proposition \ref{resultat_reconstruction_proximite_oaeirgbvodbv}, on a

\begin{equation}\label{majoration_psi_j_aomerifnvmeoinv}
\psi_j(\tilde F_j,B)=\psi_j\left(\bigoplus_{i=1}^j\vect(f_i),\bigoplus_{i=1}^jB_i\right)\le c_5 \sum_{i=1}^j \psi_1(\vect(f_i),B_i)\le \frac{c_6}{q^{(N+1)/N}}
\vspace{2.6mm}
\end{equation}

avec $c_5,c_6>0$ ne dépendant que de $n$ et de $F$. \\

De plus, on a $F\supset \tilde F_j$, donc d'après le corollaire \ref{proximite_et_inclusions_aeoibnfvoidsn} :

\[\psi_j(F,B)\le \psi_j(\tilde F_j,B).\]

Ainsi, avec \eqref{majoration_1q_aprighmoidvn} et \eqref{majoration_psi_j_aomerifnvmeoinv}, il existe une constante $c_7>0$ ne dépendant que de $n$ et de $F$ telle que 

\begin{equation}\label{majoration_psij_aefombrriehv}
\psi_j(F,B)\le \frac{c_6}{q^{(N+1)/N}}\le \frac{c_7}{H(B)^{(N+1)/(jN)}}\le \frac{c_{7}}{H(B)^{(jn-jd+j^2/2+j/2+1)/(j(jn-jd+j^2/2+j/2))}},
\vspace{2.6mm}
\end{equation}

car $N\le jn-jd+j^2/2+j/2$. \\

Appliquons finalement le théorème du Going-up par récurrence sur $\ell\in\{j,\ldots,e\}$. Si $\ell=j$, le sous-espace $B$ construit convient. Soit $\ell\in\{j,\ldots,e-1\}$, supposons par hypothèse de récurrence qu'il existe un sous-espace rationnel $C_\ell\in\mathfrak R_n(\ell)$ tel que $B\subset C_\ell$ et 

\[\psi_j(F,C_\ell)\le \frac{c^{(\ell)}}{H(C_\ell)^{\alpha(n-j)/(n-\ell)}}\]

avec $c^{(\ell)}>0$ dépendant uniquement de $n$ et de $F$ et 

\[\alpha=\frac{(jn-jd+j^2/2+j/2+1)}{j^2(n-d+j/2+1/2)}.\]

Alors d'après le théorème \ref{goingup_eorihgzmefvnfon} du Going-up appliqué avec

\[\begin{cases} H=H(C_\ell)\\x_j=\alpha(n-j)/(n-\ell)\\y_j=0\\\forall i\in\{1,\ldots,\min(d,e)\}\setminus\{j\},\quad x_i=y_i=0,\end{cases}\]

il existe un sous-espace rationnel $C_{\ell+1}\in\mathfrak R_n(\ell+1)$ tel que $C_{\ell+1}\supset C_\ell\supset B$, et il existe $c^{(\ell+1)}>0$ dépendant uniquement de $n$ et $F$, tels que

\[\psi_j(F,C_{\ell+1})\le \frac{c^{(\ell+1)}}{H(C_{\ell+1})^{\alpha(n-\ell)(n-j)/((n-\ell)(n-\ell-1))}}= \frac{c^{(\ell+1)}}{H(C_{\ell+1})^{\alpha(n-j)/(n-\ell-1)}},\]

ce qui conclut la récurrence. \\

Posons $C=C_e$, on a ainsi construit un sous-espace rationnel $C\in\mathfrak R_n(e)$ tel que $C\supset B$ et

\begin{equation}\label{maj_psi_j_finale_zoihgfnmoibneo}
\psi_j(F,C)\le \frac{c_8}{H(C)^{\beta}}
\vspace{2.6mm}
\end{equation}

avec $c_8>0$ ne dépendant que de $n$ et de $F$, et 

\[\beta=\frac{(n-j)(jn-jd+j^2/2+j/2+1)}{j^2(n-e)(n-d+j/2+1/2)}.\]

De plus, on a $F\in\mathfrak I_n(d,e)_j$, donc pour tout $C\in\mathfrak R_n(e)$, on a $\psi_j(F,C)\ne 0$. Ainsi, s'il n'existait qu'un nombre fini de sous-espaces rationnels $C$ vérifiant \eqref{maj_psi_j_finale_zoihgfnmoibneo}, on aurait une constante $c>0$ telle que

\begin{equation}\label{minoration_psi_j_par_constante_aoeifnvoienv}
\forall C\in\mathfrak R_n(e),\quad \psi_j(F,C)>c.
\vspace{2.6mm}
\end{equation}

Or l'inégalité \eqref{majoration_psij_aefombrriehv} donne, compte tenu du corollaire \ref{proximite_et_inclusions_aeoibnfvoidsn} : 

\[\psi_j(F,C)\le\psi_j(F,B)\xrightarrow[q\to\infty]{}0,\]

ce qui contredit \eqref{minoration_psi_j_par_constante_aoeifnvoienv}. \\

Finalement, il existe une infinité de sous-espaces rationnels $C$ ainsi construits approchant $F$ à l'exposant $\beta$. On a donc $\muexpA nFej\ge \beta$, d'où 

\[\muexp ndej\ge \beta,\]

ce qui termine la démonstration du théorème \ref{premiereborneobtenue_amorimeovbn}.
\end{preuve}

\chapter{Inclusion dans un sous-espace vectoriel rationnel}\label{chap_inclusion_sev_rationnel_amoriegmoingoa}

On s'intéresse dans ce chapitre à un cas particulier : le cas où le sous-espace vectoriel $A$ que l'on souhaite approcher est inclus dans un sous-espace rationnel de $\R^n$. Cette étude sera notamment utile pour le chapitre \ref{chapitre_spectre_aoimhremvoin}. \\

Le résultat principal est le théorème \ref{th_inclusion_sev_rationnel_apeivpinpiaenv} ci-dessous. En le combinant notamment avec des résultats du chapitre \ref{chap_cas_particuliers_oamofbvoaube}, on obtient de nouvelles majorations de $\muexp ndd1$, notamment $\muexp 5221\le 3$. Ces corollaires sont énoncés et démontrés dans la section \ref{section_app_du_th_inclusiion_sev_rationnel_vaoihbmoifv}, juste après l'énoncé du théorème \ref{th_inclusion_sev_rationnel_apeivpinpiaenv} et l'esquisse de sa preuve (section \ref{section_inclusion_sev_rationnel_gaomribaoivsd}). Enfin, les lemmes qui forment cette esquisse sont démontrés dans la section \ref{sssection_preuves_lemmes_sevrationnelinclus_ameioboeivboibvz}.

\section{Le théorème principal}\label{section_inclusion_sev_rationnel_gaomribaoivsd}

L'énoncé principal de ce chapitre est le suivant.

\begin{theoreme}\label{th_inclusion_sev_rationnel_apeivpinpiaenv}

Soient $n\ge 2$ et $k\in\{2,\ldots,n\}$. Soient $d,e\in\{1,\ldots,k-1\}$ tels que $d+e\le k$, et $j\in\{1,\ldots,\min(d,e)\}$. Soit $A$ un sous-espace vectoriel de $\R^n$ de dimension $d$ tel qu'il existe un sous-espace vectoriel rationnel $F\in\mathfrak R_n(k)$ vérifiant $A\subset F$.

Notons $\varphi$ un isomorphisme rationnel de $F$ dans $\R^k$ et $\tilde A=\varphi(A)$, qui est un sous-espace vectoriel de dimension $d$ de $\R^k$.

Supposons que pour tout sous-espace rationnel $B'$ de $F$ de dimension $e$, on a 

\begin{equation}\label{condition_dirr_sur_A_ameirnmoafvonsd}
\dim (A\cap B')<j.
\vspace{2.6mm}
\end{equation}

Alors $A\in\mathfrak I_n(d,e)_j$, $\tilde A\in\mathfrak I_k(d,e)_j$ et 

\[\muexpA nAej=\muexpA k{\tilde A}ej.\]

\end{theoreme}

On peut déjà remarquer que l'hypothèse \eqref{condition_dirr_sur_A_ameirnmoafvonsd}, faite pour $B'\subset F$ rationnel de dimension $e$, est \emph{a priori} une version faible de l'hypothèse $A\in\mathfrak I_n(d,e)_j$ (\emph{i.e.} $\dim(A\cap B)<j$ pour tout $B\in\mathfrak R_n(e)$). Le théorème \ref{th_inclusion_sev_rationnel_apeivpinpiaenv} montre notamment que ces deux hypothèses sont équivalentes. \\

On se place dans toute la suite de cette sous-section dans le cadre des hypothèses du théorème \ref{th_inclusion_sev_rationnel_apeivpinpiaenv}. Des lemmes \ref{condition_dirr_sur_A_argpienvoinaofmbnvuoeabv}, \ref{premiere_inegalite_moangvomasvoausbv} et \ref{seconde_inegalite_aeoribnoibnoindfvdc} découle directement le théorème \ref{th_inclusion_sev_rationnel_apeivpinpiaenv}.

\begin{lemme}\label{condition_dirr_sur_A_argpienvoinaofmbnvuoeabv}

On a 

\[A\in\mathfrak I_n(d,e)_j\quad\text{ et }\quad\tilde A\in\mathfrak I_k(d,e)_j.\]

\end{lemme}

On démontre ensuite, pour commencer, l'inégalité la plus facile : elle vient du fait que tout sous-espace vectoriel $\tilde B$ de $\R^k$ s'écrit $\varphi(B)$, avec $B$ un sous-espace vectoriel de $\R^n$ (inclus dans $F$).

\begin{lemme}\label{premiere_inegalite_moangvomasvoausbv}

On a 

\[\muexpA nAej\ge\muexpA k{\tilde A}ej.\]

\end{lemme}

Pour démontrer l'inégalité inverse, on a besoin de savoir comment se comporte la proximité entre deux sous-espaces vis-à-vis d'une projection.

\begin{lemme}\label{inegalitesurladistance_ainreoaeonbe}

Soit $\mathcal R$ une partie non vide de $\R^n$ telle que $\mathcal R\cap F^\perp=\emptyset$ et telle qu'il existe une constante $c>0$ vérifiant 

\begin{equation}\label{hyp_lemme_proximite_projection_aoeibnoivnc}
\forall X\in\mathcal R,\quad \norme{p_F^\perp(X)}\ge c\norme X
\vspace{2.6mm}
\end{equation}

où $p_F^\perp$ est la projection orthogonale sur $F$.

Soit $D$ un sous-espace vectoriel de $\R^n$ tel que $\dim D\ge j$, et vérifiant $D\subset \mathcal R\cup\{0\}$.

Alors il existe une constante $c'>0$ dépendant seulement de $c$ telle que

\[\psi_j(A,D)\ge c'\psi_j(A,p_F(D)).\]

\end{lemme}

Il est alors possible de montrer l'inégalité inverse. L'idée est que si un sous-espace rationnel approche $A$, alors son projeté orthogonal sur le sous-espace rationnel $F$ approche bien $A$ lui-aussi ; via $\varphi$ on obtient alors un sous-espace vectoriel de $\R^k$ qui approche bien $\tilde A$.

\begin{lemme}\label{seconde_inegalite_aeoribnoibnoindfvdc}

On a 

\[\muexpA nAej\le\muexpA k{\tilde A}ej.\]

\end{lemme}

\section{Applications du théorème principal}\label{section_app_du_th_inclusiion_sev_rationnel_vaoihbmoifv}

Déduisons du théorème \ref{th_inclusion_sev_rationnel_apeivpinpiaenv} une amélioration de la meilleure majoration connue de $\muexp ndej$ dans quelques cas particuliers. \\

Pour commencer, la proposition \ref{amelioration_R5_grace_inclusion_R4_aoibvoeubv} améliore la majoration de Schmidt donnée par le théorème \ref{borne_sup_de_Schmidt_roimghnev} : $\muexp 5221\le 4$.

\begin{proposition}\label{amelioration_R5_grace_inclusion_R4_aoibvoeubv}

On a 

\[\muexp 5221\le 3.\]

\end{proposition}

\begin{preuve}
D'après la proposition \ref{proposition_R4_mu_Axi_zgiozovdib}, le sous-espace $A_{\sqrt 2}$ défini dans la sous-section \ref{resultats_dans_R4_jaomfibenovn} vérifie

\[\muexpA 4{A_{\sqrt 2}}21=3.\]

Soit $\rho$ un isomorphisme rationnel de $\R^4$ dans $\R^4\times\{0\}\subset \R^5$. Posons $A=\rho(A_{\sqrt 2})$. \\

Comme $A_{\sqrt 2}\in\mathfrak I_4(2,2)_1$, le sous-espace $A_{\sqrt 2}$ vérifie l'hypothèse \eqref{condition_dirr_sur_A_ameirnmoafvonsd} du théorème \ref{th_inclusion_sev_rationnel_apeivpinpiaenv}. Donc d'après le théorème \ref{th_inclusion_sev_rationnel_apeivpinpiaenv}, on a $A\in\mathfrak I_5(2,2)_1$ et 

\[\muexpA 5A21=\muexpA 4{A_{\sqrt 2}}21=3,\]

ce qui conclut la preuve de la proposition \ref{amelioration_R5_grace_inclusion_R4_aoibvoeubv}.
\end{preuve}

De façon similaire, on déduit la proposition suivante du résultat de Moshchevitin énoncé à la sous-section \ref{sssection_moshchevitin_gaoirhgpihvd}.

\begin{proposition}\label{amelioration_grace_inclusion_sev_rationnel_aoeifnboinvcc}

Soit $n\ge 6$ et $d\in\{3,\ldots,\lfloor n/2\rfloor\}$, on a

\[\muexp ndd1\le 2d.\]

\end{proposition}

\begin{preuve}
D'après le théorème \ref{th_Moshchevitin_2020_general_ergoighvoifvn}, il existe un sous-espace $\tilde A\in\mathfrak I_{2d}(d,d)_1$ tel que

\[\muexpA {2d}{\tilde A}d1\le 2d.\]

Soit $\rho$ un isomorphisme rationnel de $\R^{2d}$ dans $\R^{2d}\times\{0\}^{n-2d}\subset\R^n$. Posons $A=\rho(\tilde A)$. \\

Alors d'après le théorème \ref{th_inclusion_sev_rationnel_apeivpinpiaenv}, on a $A\in\mathfrak I_n(d,d)_1$ et 

\[\muexpA nAd1=\muexpA {2d}{\tilde A}d1\le 2d,\]

ce qui conclut la preuve de la proposition \ref{amelioration_grace_inclusion_sev_rationnel_aoeifnboinvcc}.
\end{preuve}

\begin{remarque}

D'après le théorème \ref{borne_sup_de_Schmidt_roimghnev}, on a

\[\muexp ndd1\le \left\lceil\frac{dn-d^2+1}{n-2d+1}\right\rceil.\]

La proposition \ref{amelioration_grace_inclusion_sev_rationnel_aoeifnboinvcc} améliore cette majoration dans de nombreux cas où $n$ est proche de $2d$. Par exemple les majorants $\muexp 7331\le 7$, $\muexp 9441\le 11$ et $\muexp {10}441\le 9$ sont améliorés respectivement par $6$, $8$ et $8$.

\end{remarque}

De façon identique à la proposition \ref{amelioration_grace_inclusion_sev_rationnel_aoeifnboinvcc}, le théorème \ref{theoreme_avec_moshchevitin_goingup_aeronfvn} permet d'obtenir grâce au théorème \ref{th_inclusion_sev_rationnel_apeivpinpiaenv} la proposition suivante.

\begin{proposition}\label{applicitation_deduction_Mosh_inclusion_sev_rationnel_naveoimfvn}

Soit $n\ge 6$ et $d\in\{3,\ldots,\lfloor n/2\rfloor\}$, on a

\[\muexp nd{d-1}1\le\frac{2d^2}{d+1}.\]

\end{proposition}

\begin{remarque}

La proposition \ref{applicitation_deduction_Mosh_inclusion_sev_rationnel_naveoimfvn} améliore plusieurs majorations connues, notamment quand $n$ est proche de $2d$. Par exemple, les majorants $\muexp 9431\le 7$, $\muexp{11}541\le 10$ et $\muexp{12}541\le 9$ sont améliorés respectivement par $32/5$, $25/3$ et $25/3$.

\end{remarque}

\section{Démonstration du théorème \ref{th_inclusion_sev_rationnel_apeivpinpiaenv}}\label{sssection_preuves_lemmes_sevrationnelinclus_ameioboeivboibvz}

Démontrons ici les lemmes énoncés dans la section \ref{section_inclusion_sev_rationnel_gaomribaoivsd}. \\ 

Commençons par la preuve du lemme \ref{condition_dirr_sur_A_argpienvoinaofmbnvuoeabv} : $A\in\mathfrak I_n(d,e)_j$ et $\tilde A\in\mathfrak I_k(d,e)_j$.

\begin{preuve}
Soit $B\in\mathfrak R_n(e)$. On remarque que $B'=B\cap F$ est un sous-espace rationnel de dimension $e'\le e\le \dim F$. Donc il existe un sous-espace rationnel $B''\subset F$ contenant $B'$, avec $\dim B''=e$. L'hypothèse \eqref{condition_dirr_sur_A_ameirnmoafvonsd} montre que $\dim(A\cap B'')<j$, d'où 

\[\dim(A\cap B')<j.\]

Or on a

\[A\cap B=A\cap F\cap B=A\cap B'\]

car $A\subset F$, donc $\dim(A\cap B)<j$ : on a $A\in\mathfrak I_n(d,e)_j$. \\

Soit $\tilde B\in\mathfrak R_k(e)$. Posons

\[B=\varphi^{-1}(\tilde B)\in\mathfrak R_n(e).\]

Or $\varphi$ est un isomorphisme, donc

\[\dim(\tilde A\cap \tilde B)=\dim(\varphi(A)\cap\varphi(B))=\dim(\varphi(A\cap B))=\dim(A\cap B)<j\]

car $B\in\mathfrak R_n(e)$ et $A\in\mathfrak I_n(d,e)_j$. Cela montre que $\tilde A\in\mathfrak I_k(d,e)_j$ et termine la preuve du lemme \ref{condition_dirr_sur_A_argpienvoinaofmbnvuoeabv}.
\end{preuve}

Démontrons maintenant l'inégalité du lemme \ref{premiere_inegalite_moangvomasvoausbv} : $\muexpA nAej\ge \muexpA k{\tilde A}ej$.

\begin{preuve}
Soit $\alpha<\muexpA k{\tilde A}ej$. Alors il existe une suite $(\tilde B_N)_{N\ge 0}$ de sous-espaces rationnels de $\R^k$ de dimension $e$, deux à deux distincts, tels que pour tout $N$ suffisamment grand :

\[\psi_j(\tilde A,\tilde B_N)\le \frac{1}{H(\tilde B_N)^\alpha}.\]

Pour tout $N\in\N$, on pose

\[B_N=\varphi^{-1}(\tilde B_N)\in\mathfrak R_n(e)\]

car $\varphi$ est un isomorphisme rationnel. \\

D'après la proposition \ref{inegalitesurlahauteur_aroibvaoivbnoi}, il existe $c_{\varphi^{-1}}$ telle que pour tout $N\in\N$, 

\[H(B_N)=H(\varphi^{-1}(\tilde B_N))\le c_{\varphi^{-1}}H(\tilde B_N).\]

De plus, d'après la proposition \ref{proximite_transfo_rationnelle_inv_moiefhvoizcnzou}, il existe une constante $c_{\varphi^{-1}}'>0$ telle que

\[\psi_j(A,B_N)=\psi_j(\varphi^{-1}(\tilde A),\varphi^{-1}(\tilde B_N))\le c_{\varphi^{-1}}'\psi_j(\tilde A,\tilde B_N).\]

Finalement, pour $N$ assez grand,

\[\psi_j(A,B_N)\le c_{\varphi^{-1}}'\psi_j(\tilde A,\tilde B_N)\le \frac{c_{\varphi^{-1}}'}{H(\tilde B_N)^\alpha}\le \frac {c_1}{H(B_N)^\alpha},\]

où $c_1>0$ ne dépend que de $\varphi$. Comme les $B_N$ sont encore deux à deux distincts, on en déduit

\[\muexpA nAej\ge\alpha.\]

Ceci étant vrai pour tout $\alpha<\muexpA k{\tilde A}ej$, on obtient

\[\muexpA nAej\ge \muexpA k{\tilde A}ej,\]

ce qui termine la preuve du lemme \ref{premiere_inegalite_moangvomasvoausbv}.
\end{preuve}

Démontrons le lemme \ref{inegalitesurladistance_ainreoaeonbe} donnant le comportement de la proximité vis-à-vis d'une projection. La preuve de ce lemme suit les idées de la preuve du lemme 13 page 446 de l'article \cite{schmidt67}. Ce dernier démontre, sous des hypothèses plus faibles, un résultat similaire avec une transformation inversible au lieu d'une projection.

\begin{preuve}
En utilisant l'hypothèse \eqref{hyp_lemme_proximite_projection_aoeibnoivnc} du lemme et le fait que $p_F^\perp$ soit une projection orthogonale, il existe une constante $c>0$ telle que

\begin{equation}\label{encadrementdelanormedelaprojection}
\forall X\in\mathcal R,\quad c\norme X\le \norme{p_F^\perp(X)}\le \norme X.
\vspace{2.6mm}
\end{equation}

En particulier on a $c\le 1$ car $\mathcal R\ne\emptyset$, donc on peut supposer que $F\setminus\{0\}\subset \mathcal R$ puisque $\norme{p_F^\perp(X)}=\norme X$ pour tout $X\in F$. \\

Montrons qu'il existe une constante $c_1>0$ (qui dépend seulement de $c>0$) telle que

\begin{equation}\label{inegalite_projection_aerinbaoeinvipnczd}
\forall X\in F\setminus\{0\},\quad \forall Y\in\mathcal R,\quad \psi(X,p_F^\perp(Y))\le c_1\psi(X,Y).
\vspace{2.6mm}
\end{equation}

Soient $X\in F\setminus\{0\}$ et $Y\in\mathcal R$. Comme 

\[\psi\left(\frac X{\norme X},\frac Y{\norme Y}\right)=\psi(X,Y)\]

et que la projection $p_F^\perp$ est linéaire, on peut supposer que $\norme X=\norme Y=1$. De plus, on peut supposer $X\cdot Y\ge 0$ quitte à remplacer $Y$ par $-Y$. \\

On a
\begin{align*}
\psi^2(X,p_F^\perp(Y))
	&=1-\cos^2\widehat{(X,p_F^\perp(Y))}  \\
	&=1-\frac{(X\cdot p_F^\perp(Y))^2}{\norme{p_F^\perp(Y)}^2} \\
	&=\frac {\norme{p_F^\perp(Y)}^2-(X\cdot p_F^\perp(Y))^2}{\norme{p_F^\perp(Y)}^2}.
\end{align*}
Posons $\lambda=\norme{p_F^\perp(Y)}$ ; on remarque alors que
\begin{align*}
0\le (X\cdot p_F^\perp(Y)-\lambda)^2 
	&\iff 0\le (X\cdot p_F^\perp(Y))^2-2\lambda(X\cdot p_F^\perp(Y))+\lambda^2 \\
	&\iff \lambda^2-(X\cdot p_F^\perp(Y))^2\le 2(\lambda^2-\lambda(X\cdot p_F^\perp(Y))).
\end{align*}
Ainsi,
\begin{align*}
\psi^2(X,p_F^\perp(Y))
	&\le \frac{2(\lambda^2-\lambda(X\cdot p_F^\perp(Y)))}{\norme{p_F^\perp(Y)}^2} \\
	&=2\frac {\norme{p_F^\perp(Y)}^2-\norme{p_F^\perp(Y)}(X\cdot p_F^\perp(Y))}{\norme{p_F^\perp(Y)}^2} \\
	&=2\frac {\norme{p_F^\perp(Y)}-X\cdot p_F^\perp(Y)}{\norme{p_F^\perp(Y)}}.
\end{align*}

Utilisons alors la première inégalité de \eqref{encadrementdelanormedelaprojection} qui donne

\[\psi^2(X,p_F^\perp(Y))\le\frac 2{c} \left(\norme{p_F^\perp(Y)}-X\cdot p_F^\perp(Y)\right).\]

Or

\[\frac 12\left(\norme{X-p_F^\perp(Y)}^2-1-\norme{p_F^\perp(Y)}^2\right)=-X\cdot p_F^\perp(Y),\]

donc car $X\in F$ et que $p_F^\perp$ est une application linéaire :
\begin{align*}
\psi^2(X,p_F^\perp(Y))
	&\le\frac 1{c}\left(2\norme{p_F^\perp(Y)}+\norme{X-p_F^\perp(Y)}^2-1-\norme{p_F^\perp(Y)}^2\right) \\
	&=\frac 1{c}\norme{p_F^\perp(X-Y)}^2-\frac1{c}\left(\norme{p_F^\perp(Y)}-1\right)^2 \\
	&\le \frac 1{c}\norme{p_F^\perp(X-Y)}^2 \\
	&\le \frac {1}{c}\norme{X-Y}^2.
\end{align*}
De plus, on a $X\cdot Y\ge 0$, donc d'après le lemme \ref{lemme_minoration_psiXY_amoifbvoisbdvoibs} : 

\[\psi(X,Y)\ge\frac {\sqrt 2}2\norme{X-Y},\]

donc

\[\psi^2(X,p_F^\perp(Y))\le \frac {1}{c}\norme{X-Y}^2\le \frac 2{c}\psi(X,Y)^2\]

soit

\[\psi(X,p_F^\perp(Y))\le c_1\psi(X,Y)\]

ce qui est bien l'inégalité \eqref{inegalite_projection_aerinbaoeinvipnczd} voulue, avec $c_1=\sqrt{2/c}$. \\

D'après la proposition \ref{prop_proximite_aveclesmin_maoierngmoinv}, il existe $A_j$ un sous-espace de $A$ de dimension $j$, tel que

\[\forall X\in A_j\setminus\{0\},\quad \exists Y\in D\setminus\{0\},\quad \psi(X,Y)\le \psi_j(A,D).\]

Soit $X\in A_j\setminus\{0\}$. Il existe donc $Y\in D\setminus\{0\}$ tel que 

\[\psi(X,Y)\le \psi_j(A,D).\]

Or $X\in A_j\subset A\subset F$ et $Y\in D\setminus\{0\}\subset \mathcal R$, on peut donc utiliser l'inégalité \eqref{inegalite_projection_aerinbaoeinvipnczd}, qui donne

\[\psi(X,p_F^\perp(Y))\le c_1\psi(X,Y)\le c_1\psi_j(A,D).\]

On a donc trouvé $Y'=p_F^\perp(Y)\in p_F^\perp(D)$ non nul (car $Y\in\mathcal R$ donc $Y\notin F^\perp)$ tel que

\[\psi(X,Y')\le c_1\psi_j(A,D).\]

Autrement dit, on a montré que :

\[\forall X\in A\setminus\{0\},\quad \exists Y'\in p_F^\perp(D)\setminus\{0\},\quad \psi(X,Y')\le c_1\psi_j(A,D).\]

D'après la proposition \ref{prop_proximite_aveclesmin_maoierngmoinv}, $\psi_j(A,p_F^\perp(D))$ est le plus petit réel ayant cette propriété, donc

\[\psi_j(A,p_F^\perp(D))\le c_1\psi_j(A,D)\]

qui est l'inégalité cherchée.
\end{preuve}

Démontrons maintenant le lemme \ref{seconde_inegalite_aeoribnoibnoindfvdc}, c'est-à-dire que $\muexpA nAej\le\muexpA k{\tilde A}ej$. Pour cela on commence par un lemme.

\begin{lemme}\label{lemme_intersection_orthogonal_est_vide_aorfbaovbd}

Sous les hypothèses du théorème \ref{th_inclusion_sev_rationnel_apeivpinpiaenv}, pour tout $\alpha<\muexpA nAej$, il existe une suite $(B_N)_{N\in\N}$ de sous-espaces rationnels de $\R^n$ de dimension $e$, deux à deux distincts, tels que pour tout $N$ suffisamment grand :

\[B_N\cap F^\perp=\{0\}\quad\text{ et }\quad \psi_j(A,B_N)\le \frac 1{H(B_N)^\alpha}.\]

\end{lemme}

\begin{preuve}[Lemme \ref{lemme_intersection_orthogonal_est_vide_aorfbaovbd}]
Soit $\alpha'$ tel que $\alpha<\alpha'<\muexpA nAej$. Par définition de $\muexpA nAej$, il existe une suite $(B_N)_{N\in\N}$ de sous-espaces rationnels de $\R^n$ de dimension $e$, deux à deux distincts, tels que pour tout $N$ suffisamment grand

\[\psi_j(A,B_N)\le \frac 1{H(B_N)^{\alpha'}}.\]

La difficulté est que $B_N\cap F^\perp$ peut être non nul. Notons $(g_1,\ldots,g_{n-k})$ une famille libre de vecteurs de $\Q^n$ telle que

\begin{equation}\label{initialisation_rec_inclusionsevrationnel_painerinv}
\R^n=F\oplus\vect(g_1,\ldots,g_{n-k}).
\vspace{2.6mm}
\end{equation}

On va construire des vecteurs $g_{n-k+1},\ldots,g_n\in\Q^n$ pour compléter la famille $(g_1,\ldots,g_{n-k})$. On pourra alors poser pour $I\subset \{1,\ldots,n\}$ :

\[G_I=\vect\{g_i,\ i\in I\}.\]

Pour $\ell\in\N^*$, notons $\mathcal P_m(\ell)$ l'ensemble des parties à $m$ éléments de $\{1,\ldots,\ell\}$. On montre par récurrence finie sur $\ell\in\{n-k,\ldots,n\}$ qu'il existe $g_{n-k+1},\ldots,g_\ell$ des vecteurs de $\Q^n$ tels que

\begin{equation}\label{rec_construction_gi_inclusionsevrationnel_aepirbniv}
\begin{cases}\dim\vect(g_1,\ldots,g_\ell)=\ell\\\forall I\in\mathcal P_{n-k}(\ell),\quad G_I\cap F=\{0\}.\end{cases}
\vspace{2.6mm}
\end{equation} 

Comme $\dim F=k$, la récurrence est initialisée pour $\ell=n-k$ avec \eqref{initialisation_rec_inclusionsevrationnel_painerinv}. \\

Soit $\ell\in\{n-k,\ldots,n-1\}$, supposons que les $g_i$ pour $i\in\{n-k+1,\ldots,\ell\}$ ont été construits et vérifient \eqref{rec_construction_gi_inclusionsevrationnel_aepirbniv}. Posons

\[G=\vect(g_1,\ldots,g_\ell)\cup\bigcup_{K\in\mathcal P_{n-k-1}(\ell)} \big(F\oplus G_K\big).\]

L'ensemble $G$ est une union d'un nombre fini de sous-espaces vectoriels de dimension $n-1$, et d'un sous-espace de dimension $\ell\le n-1$. On peut donc se donner un vecteur

\[g_{\ell+1}\in\Q^n\setminus G.\]

Montrons que ce vecteur vérifie bien \eqref{rec_construction_gi_inclusionsevrationnel_aepirbniv}. Comme $g_{\ell+1}\notin\vect(g_1,\ldots,g_\ell)$, on a bien 

\[\dim\vect(g_1,\ldots,g_{\ell+1})=\ell+1.\]

On suppose par l'absurde qu'il existe $I\in\mathcal P_{n-k}(\ell+1)$ tel que $G_I\cap F\ne\{0\}$. Soit $u\in G_I\cap F\setminus\{0\}$. Par hypothèse de récurrence, on a $I\notin\mathcal P_{n-k}(\ell)$, donc $\ell+1\in I$. On peut donc décomposer $I$ sous la forme $I=K\cup\{\ell+1\}$ avec $K\in\mathcal P_{n-k-1}(\ell)$. Comme $u\notin F\cap G_K=\{0\}$, il existe $\alpha_{\ell+1}\ne 0$ et des $\alpha_i\in\R$ (pour $i\in K$) tels que

\[u=\alpha_{\ell+1} g_{\ell+1}+\sum_{i\in K} \alpha_ig_i.\]

Donc

\[g_{\ell+1}=\frac 1{\alpha_{\ell+1}}\left(u-\sum_{i\in K}\alpha_ig_i\right)\in F\oplus G_K,\]

ce qui est absurde par définition de $G$ car $g_{\ell+1}\notin G$. Donc $g_{\ell+1}$ vérifie \eqref{rec_construction_gi_inclusionsevrationnel_aepirbniv}, ce qui termine la récurrence. \\

Finalement, on a construit des vecteurs $g_1,\ldots,g_n$ vérifiant \eqref{rec_construction_gi_inclusionsevrationnel_aepirbniv}. \\

Pour $I\in\mathcal P_{n-k}(n)$, on a

\[\dim G_I=n-k=\dim F^\perp.\]

Donc $G_I\oplus F=\R^n$ : il existe un isomorphisme rationnel $\rho_I\in\GL_n(\Q)$ tel que

\[\begin{cases}{\rho_I}_{\vert F}=\id_F\\\rho_I(G_I)=F^\perp.\end{cases}\]

Soit $N\in\N$. Supposons par l'absurde que

\[\forall I\in\mathcal P_{n-k}(n),\quad \rho_I(B_N)\cap F^\perp\ne\{0\},\]

soit

\[\forall I\in\mathcal P_{n-k}(n),\quad B_N\cap \rho_I^{-1}(F^\perp)\ne\{0\},\]

\emph{i.e.} par définition de $\rho_I$ :

\[\forall I\in\mathcal P_{n-k}(n),\quad B_N\cap G_I\ne\{0\}.\]

Notons

\[J=\left\{i\in\{1,\ldots,n\},\ \exists \alpha\ne 0,\quad\exists\lambda_1,\ldots,\lambda_{i-1}\in\R,\quad \alpha g_i +\sum_{\ell=1}^{i-1} \lambda_\ell g_\ell\in B_N\right\}.\]

Supposons dans un premier temps que $\card(J)\le k$. Il existe donc $I\in\mathcal P_{n-k}(n)$ tel que $I\cap J=\emptyset$. Or $B_N\cap G_I\ne\{0\}$, donc il existe une famille de réels $(\beta_i)_{i\in I}\in\R^{I}$ non tous nuls, tels que

\[\sum_{i\in I}\beta_i g_i\in B_N.\]

Notons $i_0$ le plus grand $i\in I$ tel que $\beta_i\ne 0$. Alors

\[\sum_{i\in I}\beta_i g_i=\alpha g_{i_0}+\sum_{i=1}^{i_0-1}\beta_i g_i,\]

en posant $\beta_i=0$ si $i\notin I$ et $\alpha=\beta_{i_0}\ne 0$. Donc $i_0\in J\cap I$, ce qui est absurde. \\

On a donc $\card(J)>k$. Les éléments de $J$ donnent au moins $k+1$ vecteurs linéairement indépendants de $B_N$, ce qui est absurde car $\dim B_N\le k$. \\

Finalement, on a montré que

\begin{equation}\label{rhoi_cap_Fperp_egal0_aberoinbf}
\exists I\in\mathcal P_{n-k}(n),\quad \rho_I(B_N)\cap F^\perp=\{0\},
\vspace{2.6mm}
\end{equation}

et on se donne un tel $I\in\mathcal P_{n-k}(n)$. \\

On a $\rho_I\in\GL_n(\R)$, donc d'après la proposition \ref{proximite_transfo_rationnelle_inv_moiefhvoizcnzou}, il existe une constante $c(\rho_I)>0$ telle que

\[\psi_j(A,\rho_I(B_N))=\psi_j(\rho_I(A),\rho_I(B_N))\le c(\rho_I)\psi_j(A,B_N)\,;\]

en effet $\rho_I(A)=A$ puisque ${\rho_I}_{\vert F}=\id_F$. \\

Alors en posant 

\[c_2=\max_{I\in\mathcal P_{n-k}(n)}c(\rho_I)>0,\]

on a une constante indépendante de $B_N$ telle que

\[\psi_j(A,\rho_I(B_N))\le c_2\psi_j(A,B_N).\]

De plus

\[\dim (\rho_I(B_N))=\dim(B_N)=e\]

car $\rho_I$ est un isomorphisme. On peut donc appliquer la proposition \ref{inegalitesurlahauteur_aroibvaoivbnoi}, qui donne une constante $c'(\rho_I)>0$ telle que

\[H(\rho_I(B_N))\le c'(\rho_I)H(B_N),\]

et en posant 

\[c_3=\max_{I\in\mathcal P_{n-k}(n)}c'(\rho_I)>0,\]

on a donc une constante indépendante de $B_N$ telle que

\[H(\rho_I(B_N))\le c_3H(B_N).\]

On a donc, si $N$ est suffisamment grand : 

\[\psi_j(A,\rho_I(B_N))\le c_2\psi_j(A,B_N)\le \frac{c_2}{H(B_N)^{\alpha'}}\le \frac{c_2c_3^{-\alpha'}}{H(\rho_I(B_N))^{\alpha'}}\le \frac 1{H(\rho_I(B_N))^\alpha}\]

ce qui termine la preuve du lemme \ref{lemme_intersection_orthogonal_est_vide_aorfbaovbd}.
\end{preuve}

\begin{preuve}[Lemme \ref{seconde_inegalite_aeoribnoibnoindfvdc}]
Soit $\alpha<\muexpA nAej$. Le lemme \ref{lemme_intersection_orthogonal_est_vide_aorfbaovbd} fournit une suite $(B_N)_{N\in\N}$ de sous-espaces rationnels de $\R^n$ de dimension $e$, deux à deux distincts, tels que pour tout $N$ suffisamment grand :

\begin{equation}\label{def_des_B_N_inclusionsevrationnel_aboemfbv}
B_N\cap F^\perp=\{0\}\quad\text{ et }\quad\psi_j(A,B_N)\le\frac{1}{H(B_N)^\alpha}.
\vspace{2.6mm}
\end{equation}

Soit $N\in\N$ suffisamment grand. On note $p_F^\perp$ la projection orthogonale sur $F$. Comme $F\in\mathfrak R_n(k)$, $p_F^\perp$ est un endomorphisme rationnel de $\R^n$. \\

Posons $B_N'=p_F^\perp(B_N)$ le projeté orthogonal de $B_N$ sur $F$. Comme $B_N$ est un sous-espace vectoriel rationnel, $B_N'$ est aussi un sous-espace vectoriel rationnel. \\

On cherche à appliquer le lemme \ref{inegalitesurladistance_ainreoaeonbe}. Pour cela, il faut choisir convenablement un ensemble $\mathcal R$ vérifiant l'hypothèse du lemme. \\

Posons $\mathcal R$ -- représenté sur la figure \ref{figure_cone_aeoinoinevne} -- l'ensemble des vecteurs non nuls de $\R^n$ formant un angle strictement inférieur à $\pi/4$ avec le sous-espace $F$, \emph{i.e.}

\begin{equation}\label{def_partie_R_aivofsnvoizdnbv}
\mathcal R=\left\{X\in\R^n\setminus\{0\},\ \psi_1(F,\vect(X))< \frac{\sqrt 2}{2}\right\}.
\vspace{2.6mm}
\end{equation}

On a bien $\mathcal R\cap F^\perp=\emptyset$.

\begin{figure}[H]
\begin{center}
\includegraphics[scale=0.73]{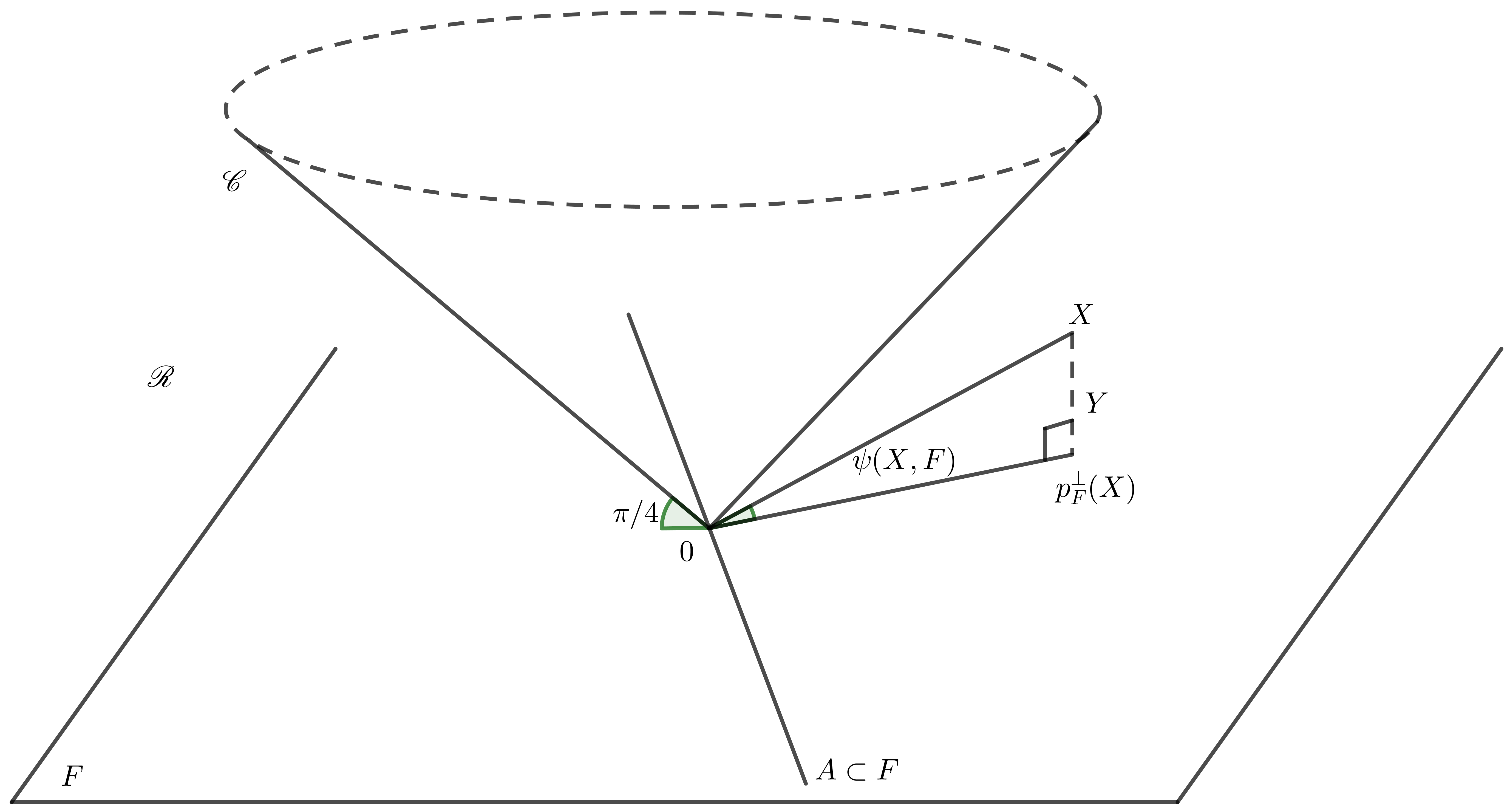}
\caption{Le cône $\mathcal C$ faisant un angle $\pi/4$ avec le sous-espace $F$}
\label{figure_cone_aeoinoinevne}
\end{center}
\end{figure}

Sur la figure \ref{figure_cone_aeoinoinevne}, on a représenté le cône 

\[\mathcal C=\{0\}\cup\left\{X\in\R^n\setminus\{0\},\ \psi_1(F,\vect(X))= \frac{\sqrt 2}{2}\right\},\]

l'ensemble $\mathcal R$ défini en \eqref{def_partie_R_aivofsnvoizdnbv} est alors la portion d'espace comprise entre $F$ et $\mathcal C$. \\

Soit $X\in\mathcal R$. En posant $Y=X-p_F^\perp(X)$ (représenté sur la figure \ref{figure_cone_aeoinoinevne}), on a

\[\norme{p_F^\perp(X)}^2+\norme Y^2=\norme X^2.\]

Or $X\in\mathcal R$, donc d'après les lemmes \ref{lemme_angle_projete_ortho_amotinfqovn} et \ref{lemme_proximite_proj_ortho_apamoebnamon} respectivement, on a

\[\psi(X,F)=\psi(X,p_F^\perp(X))=\frac{\norme Y}{\norme X} <\frac{\sqrt 2}2,\]

donc

\[\norme Y\le \frac{\sqrt 2}{2}\norme X.\]

Ainsi,

\[\norme{p_F^\perp(X)}^2=\norme X^2-\norme Y^2\ge \norme X^2-\frac 12 \norme X^2=\frac 12 \norme X^2,\]

et la partie $\mathcal R$ construite vérifie bien les hypothèses du lemme \ref{inegalitesurladistance_ainreoaeonbe} avec $c=\displaystyle\frac{\sqrt 2}{2}$. \\

Notons $B_{N,j}$ le sous-espace vectoriel de $B_N$ de dimension $j$ donné par la proposition \ref{prop_proximite_aveclesmin_maoierngmoinv}. Comme $N$ est supposé suffisamment grand, on peut supposer que 

\[\psi_j(A,B_N)\le \frac 12.\]

La proposition \ref{prop_proximite_aveclesmin_maoierngmoinv} donne que pour tout $Y\in B_{N,j}\setminus\{0\}$, il existe $X\in A\subset F$ non nul tel que 

\[\psi(X,Y)\le \psi_j(A,B_N)\le \frac 12.\]

On a donc

\[\psi_1(F,\vect(Y))\le \frac 12<\frac{\sqrt 2}2,\]

d'où $Y\in\mathcal R$. \\

Donc pour $N$ suffisamment grand,

\[B_{N,j}\setminus\{0\}\subset \mathcal R.\]

On peut ainsi appliquer le lemme \ref{inegalitesurladistance_ainreoaeonbe}, qui donne l'existence d'une constante $c_4>0$ ne dépendant ni de $A$ ni des $B_N$, telle que 

\begin{equation}\label{minoration_psijAB_N_gaepifbaiondv}
\psi_j(A,B_N)=\psi_j(A,B_{N,j})\ge c_4\psi_j(A,p_F^\perp(B_{N,j}))\ge c_4\psi_j(A,B_N')
\vspace{2.6mm}
\end{equation}

en appliquant le corollaire \ref{proximite_et_inclusions_aeoibnfvoidsn} car $B_N'=p_F^\perp(B_N)\supset p_F^\perp(B_{N,j})$. \\

Comme $B_N\cap F^\perp=\{0\}$ on a $\dim B_N'=e$ ; puisque $B_N'\subset F$, on peut poser

\[\tilde B_N=\varphi(B_N')\in\mathfrak R_k(e).\]

D'après la proposition \ref{proximite_transfo_rationnelle_inv_moiefhvoizcnzou}, il existe une constante $c_{\varphi}>0$ telle que

\[\psi_j(\varphi(A),\varphi(B_N'))\le c_\varphi\psi_j(A,B_N').\]

Soit $\beta>\muexpA k{\tilde A}ej$ ; on a ainsi en utilisant la minoration \eqref{minoration_psijAB_N_gaepifbaiondv} que pour $N$ suffisamment grand (en fonction de $\beta$) : 

\[\psi_j(A,B_N)\ge c_4\psi_j(A,B_N') \ge c_4c_{\varphi}^{-1}\psi_j(\varphi(A),\varphi(B_N'))=c_4c_{\varphi}^{-1}\psi_j(\tilde A,\tilde B_N)\ge \frac{c_5}{H(\tilde B_N)^\beta}\]

avec $c_5>0$. \\

Or d'après la proposition \ref{inegalitesurlahauteur_aroibvaoivbnoi}, il existe $c_6>0$ telle que

\[H(\tilde B_N)=H(\varphi(B_N'))\le c_6H(B_N'),\]

donc

\[\psi_j(A,B_N)\ge \frac{c_7}{H(B_N')^\beta}\]

avec $c_7>0$. \\

Comme $B_N\cap F^\perp=\{0\}$, on a

\[\dim (p_F^\perp(B_N))=\dim(B_N).\]

On peut donc appliquer la proposition \ref{inegalitesurlahauteur_aroibvaoivbnoi} qui fournit une constante $c_8>0$ telle que

\[H(B_N')=H(p_F^\perp(B_N))\le c_8 H(B_N),\]

donc, compte tenu de \eqref{def_des_B_N_inclusionsevrationnel_aboemfbv},

\[\frac 1{H(B_N)^\alpha}\ge \psi_j(A,B_N)\ge \frac{c_9}{H(B_N)^\beta}\]

avec $c_9>0$. Comme $H(B_N)$ tend vers l'infini quand $N\to+\infty$, on en déduit $\alpha\le\beta$. Ceci étant valable pour tout $\alpha<\muexpA nAej$ et pour tout $\beta>\muexpA k{\tilde A}ej$ on obtient

\[\muexpA nAej\le \muexpA k{\tilde A}ej.\]
\end{preuve}

\chapter{Le spectre de $\muexpA n{\bullet}{\ell}\ell$}\label{chapitre_spectre_aoimhremvoin}

Dans ce chapitre, on apporte une réponse partielle au problème \ref{probleme_spectre_oaerinemoivn}. \\

On s'intéresse ici au spectre de $\muexpA n{\bullet}{\ell}\ell$ sur $\mathfrak I_n(\ell,\ell)_\ell$, autrement dit à l'ensemble $\muexpA n{\mathfrak I_n(\ell,\ell)_\ell}{\ell}\ell$. Le résultat obtenu est le suivant :

\begin{theoreme}\label{theoreme_spectre_amoeribnefaomnv}

Soient $n\ge 2$ et $\ell\in\{1,\ldots,\lfloor n/2\rfloor\}$.

Alors 

\[\left[1+\frac 1{2\ell}+\sqrt{1+\frac{1}{4\ell^2}},+\infty\right]\subset\Big\{\muexpA nA\ell\ell,\ A\in\mathfrak I_n(\ell,\ell)_\ell\Big\}.\]

\end{theoreme}

Remarquons que

\[1+\frac 1{2\ell}+\sqrt{1+\frac{1}{4\ell^2}}\xrightarrow[\ell\to\infty]{} 2.\]

\begin{remarque}

Il serait intéressant de pouvoir remplacer la borne $1+1/(2\ell)+\sqrt{1+1/(4\ell^2})$ par l'exposant de Saxcé (voir sous-section \ref{sssection_de_Saxce_amoinbefoinvnie}), autrement dit de montrer que

\[\left[\frac{n}{\ell(n-\ell)},+\infty\right]\subset\Big\{\muexpA nA\ell\ell,\ A\in\mathfrak I_n(\ell,\ell)_\ell\Big\}.\]

\end{remarque}

Dans la section \ref{section_resultats_sev_exposant_donne_aefoubv} on esquisse la preuve du théorème \ref{theoreme_spectre_amoeribnefaomnv}. On commence par supposer que $n=2\ell$, et on fixe $\beta<+\infty$ dans l'intervalle du théorème \ref{theoreme_spectre_amoeribnefaomnv}. On cherche alors à construire un sous-espace vectoriel d'exposant $\beta$. Le lemme \ref{lemme_condition_dirr_sur_A_aomerbvozb} permet de construire un sous-espace $A\in\mathfrak I_{2\ell}(\ell,\ell)_\ell$. Le lemme \ref{majoration_psi2_A_BN_amroghozv} donne $\muexpA {2\ell}A\ell\ell\ge\beta$, tandis que les lemmes \ref{les_BN_sont_les_meilleurs_afnemobibvo} et \ref{minorationdupsielldeAetBN_eogihrgoeihgoi} donnent à eux deux $\muexpA{2\ell}A\ell\ell\le\beta$. Finalement, $\muexpA{2\ell}A\ell\ell=\beta$. Tous ces lemmes sont ensuite démontrés dans la section \ref{section_preuves_chapitre_exposant_donne_oaimmvobb}. \\

Le résultat est alors étendu au cas $n>2\ell$ grâce au théorème \ref{th_inclusion_sev_rationnel_apeivpinpiaenv}, et au cas $\beta=+\infty$ dans la preuve du théorème \ref{theoreme_spectre_amoeribnefaomnv} page \pageref{preuve_th_spectre_gamoerbgosn}. 

\section{Construction d'un sous-espace d'exposant prescrit}\label{section_resultats_sev_exposant_donne_aefoubv}

Soit $\ell\ge 1$ un entier ; posons $n=2\ell$. 

On se donne un réel $\beta$ tel que

\begin{equation}\label{condition_sur_beta_aomeibnoaidbnv}
\beta\ge1+\frac 1{2\ell}+\sqrt{1+\frac1{4\ell^2}},
\vspace{2.6mm}
\end{equation}

et on va construire $A\in\mathfrak I_n(\ell,\ell)_\ell$, un sous-espace vectoriel $(\ell,\ell)$-irrationnel de $\R^n$ tel que 

\[\muexpA nA\ell\ell=\beta.\]

Posons $\alpha=\ell\beta$, et pour tous $i,j\in\{1,\ldots,\ell\}$

\[\xi_{i,j}=\sum_{k=0}^\infty \frac{e^{(i,j)}_k}{\theta^{\lfloor\alpha^k\rfloor}}\]

où les $(e^{(i,j)}_k)_{k\in\N}$ sont des suites à déterminer, à valeurs dans $\{1,2\}$ si $i\ne j$ et à valeurs dans $\{2\ell,2\ell+1\}$ si $i=j$, et où $\theta$ est le plus petit nombre premier tel que

\begin{equation}\label{hyp_sur_theta_amorinvamnoavcbo}
\theta>(n+1)^{n/2}\left(\frac n2\right)!=\ell!\,(2\ell+1)^{\ell}.
\vspace{2.6mm}
\end{equation}

Le fait de choisir ici \emph{le plus petit} tel $\theta$ n'a pas d'autre intérêt que de permettre aux constantes de ne pas dépendre de $\theta$. En pratique, tout nombre premier $\theta$ vérifiant \eqref{hyp_sur_theta_amorinvamnoavcbo} conviendrait. \\

L'hypothèse \eqref{hyp_sur_theta_amorinvamnoavcbo} sur $\theta$ et le fait que les $e_k^{(i,j)}$ appartiennent à $\{1,2\}$ ou à $\{2\ell,2\ell+1\}$ seront utilisés dans la preuve du lemme \ref{lemmeminorationhauteurdesBN_aeoinmoaeivnv} pour minorer la hauteur des sous-espaces rationnels construits ci-dessous. \\

Notons $I_\ell$ la matrice identité de $\M_\ell(\R)$. Posons $M_\xi=(\xi_{i,j})_{(i,j)\in\{1,\ldots,\ell\}^2}\in\M_\ell(\R)$ la matrice des $\xi_{i,j}$ et $M_A$ la matrice définie par blocs :

\begin{equation}\label{def_MA_matricedeA_armoeigaoinv}
M_A=\begin{pmatrix} I_\ell \\ M_\xi \end{pmatrix}\in\M_{2\ell,\ell}(\R).
\vspace{2.6mm}
\end{equation}

On appelle $Y_1,\ldots,Y_{\ell}\in\R^{2\ell}$ les colonnes de $M_A$, et on note $A$ le sous-espace vectoriel de $\R^{2\ell}$ engendré par les $Y_i$ : 

\[A=\vect(Y_1,\ldots,Y_\ell).\]

Remarquons que $\rg(M_A)=\ell$, donc la famille $(Y_1,\ldots,Y_\ell)$ est libre et

\[\dim A=\ell.\]

Énonçons un premier lemme montrant que le sous-espace $A$ vérifie bien la condition de $(\ell,\ell)$-irrationalité.

\begin{lemme}\label{lemme_condition_dirr_sur_A_aomerbvozb}

Il existe des suites $(e^{(i,j)}_k)_{k\in\N}$ à valeurs dans $\{1,2\}$ si $i\ne j$ et à valeurs dans $\{2\ell,2\ell+1\}$ si $i=j$, telles que

\[A\in\mathfrak I_{n}(\ell,\ell)_1.\]

\end{lemme}

\emph{A fortiori}, on a $A\in\mathfrak I_{n}(\ell,\ell)_\ell$ puisque $\mathfrak I_n(\ell,\ell)_1\subset\mathfrak I_n(\ell,\ell)_\ell$. \\

Dans toute la suite, on fixe des suites $(e^{(i,j)}_k)_{k\in\N}$ vérifiant le lemme \ref{lemme_condition_dirr_sur_A_aomerbvozb}. \\

Le sous-espace $A$ étant alors construit, on va désormais construire des sous-espaces rationnels $B_N$ pour $N\ge 1$ approchant $A$ jusqu'au $\ell$-ième angle à l'exposant exactement $\beta$. Enfin, on montrera que ces sous-espaces $B_N$ sont ceux approchant le mieux $A$, ce qui donnera finalement $\muexpA nA\ell\ell=\beta$. \\

Posons pour $(i,j)\in\{1,\ldots,\ell\}^2$ et $N\ge 1$, 

\[f_N^{(i,j)}=\theta^{\lfloor\alpha^N\rfloor}\sum_{k=0}^N \frac{e^{(i,j)}_k}{\theta^{\lfloor\alpha^k\rfloor}}\in\Z\]

ainsi que $M_{B_N}$ la matrice par blocs

\[M_{B_N}=\begin{pmatrix} \theta^{\lfloor\alpha^N\rfloor}I_\ell \\ F_N\end{pmatrix}\in\M_{2\ell,\ell}(\Z)\]

où $F_N$ est la matrice $(f_N^{(i,j)})_{(i,j)\in\{1,\ldots,\ell\}^2}\in\M_\ell(\Z)$. \\

Notons de plus $X_N^{(1)},\ldots,X_N^{(\ell)}$ les colonnes de $M_{B_N}$, et posons

\[B_N=\vect(X_N^{(1)},\ldots,X_N^{(\ell)})\in\mathfrak R_{2\ell}(\ell).\] 

On peut remarquer que pour tous $i,j\in\{1,\ldots,\ell\}$, on a

\[\sum_{k=N+1}^\infty \frac{e_k^{(i,j)}}{\theta^{\lfloor\alpha^k\rfloor}}\le (2\ell+1)\sum_{j=\lfloor\alpha^{N+1}\rfloor}^\infty \frac 1{\theta^j}= \frac{2\ell+1}{\theta^{\lfloor\alpha^{N+1}\rfloor}}\cdot\frac 1{1-1/\theta}=\frac{2\ell+1}{\theta^{\lfloor\alpha^{N+1}\rfloor}}\cdot\frac{\theta}{\theta-1}< \frac{4\ell+2}{\theta^{\lfloor\alpha^{N+1}\rfloor}}\]

car $\theta> 2$, donc

\begin{equation}\label{le_sev_approche_bien_moaeirgnomdfvn}
0<\xi_{i,j}-\frac{f_N^{(i,j)}}{\theta^{\lfloor\alpha^N\rfloor}}=\sum_{k=N+1}^\infty \frac{e_k^{(i,j)}}{\theta^{\lfloor\alpha^k\rfloor}}<\frac {4\ell+2}{\theta^{\lfloor\alpha^{N+1}\rfloor}}.
\vspace{2.6mm}
\end{equation}

\begin{remarque}

Comme $\ell\ge 1$ et comme $\beta$ vérifie \eqref{condition_sur_beta_aomeibnoaidbnv}, on a 

\[\alpha=\ell\beta\ge\ell+\frac 12+\sqrt{\ell^2+\frac 14}\ge\frac{3+\sqrt 5}2.\]

Dans le cas $\ell=1$, on obtiendra donc comme cas particulier de la proposition \ref{exposant_donne_pour_A_aemofdunbvombc} le résultat connu sur l'exposant d'irrationalité de $\xi_{1,1}$ : 

\[\mu\left(\sum_{k=0}^\infty \frac{e^{(1,1)}_k}{\theta^{\lfloor\alpha^k\rfloor}}\right)=\alpha\]

avec $(e_k^{(1,1)})_{k\in\N}$ une suite à valeurs dans $\{2,3\}$, $\theta$ un nombre premier strictement supérieur à $3$ et $\alpha$ un réel supérieur à $(3+\sqrt 5)/2$. Les arguments de la section 8 de \cite{levesley06} permettent de montrer facilement ce résultat, mais la méthode développée ici est différente (et celle-ci ne traite pas le cas $\theta=3$). Si $2\le\alpha<(3+\sqrt 5)/2$, on a encore $\mu(\xi_{1,1})=\alpha$ grâce au théorème 2 de \cite{bugeaud08}. \\

\end{remarque}

\begin{remarque}

De façon similaire au cas $\ell=1$, on peut se demander si on a encore $\muexpA {2\ell}A\ell\ell=\beta$ lorsque $\beta$ est plus petit que $1+1/(2\ell)+\sqrt{1+1/(4\ell^2)}$.

\end{remarque}

On montre dans le lemme suivant que les sous-espaces rationnels $B_N$ approchent $A$ à l'exposant $\alpha/\ell$ : 

\begin{lemme}\label{majoration_psi2_A_BN_amroghozv}

Il existe une constante $c_1>0$ dépendant uniquement de $A$ telle que

\[\forall N\ge 1,\quad \psi_\ell(A,B_N)\le \frac{c_1}{H(B_N)^{\alpha/\ell}}.\]

\end{lemme}

Ce lemme donne

\[\muexpA {2\ell}A{\ell}{\ell}\ge \alpha/\ell.\]

Le lemme suivant montre que pour tout $\epsilon>0$, les $B_N$ sont les seuls sous-espaces rationnels de dimension $\ell$ approchant $A$ à l'exposant $\alpha/\ell+\epsilon$, ce qui implique que le meilleur exposant possible est déterminé par les $B_N$.

\begin{lemme}\label{les_BN_sont_les_meilleurs_afnemobibvo}

Supposons qu'il existe $\epsilon>0$ et $C\in\mathfrak R_{2\ell}(\ell)$ tels que

\begin{equation}\label{hyp_sur_psi2_A_C}
\psi_\ell(A,C)\le \frac 1{H(C)^{\alpha/\ell+\epsilon}}.
\vspace{2.6mm}
\end{equation}

Alors si $H(C)$ est suffisamment grand (en fonction de $\ell$ et de $\epsilon$), il existe $N\in\N^*$ tel que $C=B_N$.

\end{lemme}

Maintenant qu'on sait que si $N\in\N^*$ est suffisamment grand, les sous-espaces rationnels $B_N$ sont ceux réalisant la meilleure approximation de $A$, on va montrer qu'ils approchent le sous-espace $A$ à l'exposant au plus $\alpha/\ell$.

\begin{lemme}\label{minorationdupsielldeAetBN_eogihrgoeihgoi}

Pour $N$ suffisamment grand, on a

\[\psi_\ell(A,B_N)\ge \frac{c}{H(B_N)^{\alpha/\ell}}\]

avec $c>0$ dépendant uniquement de $A$. 
\end{lemme}

Les lemmes \ref{les_BN_sont_les_meilleurs_afnemobibvo} et \ref{minorationdupsielldeAetBN_eogihrgoeihgoi} combinés donnent

\[\muexpA {2\ell}{A}{\ell}\ell\le \alpha/\ell.\]

En comparant avec le lemme \ref{majoration_psi2_A_BN_amroghozv}, on a ainsi obtenu la proposition suivante :

\begin{proposition}\label{exposant_donne_pour_A_aemofdunbvombc}

Soit $\ell\ge 2$. Il existe $A\in\mathfrak I_{2\ell}(\ell,\ell)_\ell$ un sous-espace vectoriel $(\ell,\ell)$-irrationnel de $\R^{2\ell}$ tel que

\[\muexpA {2\ell}A\ell\ell=\beta.\]

\end{proposition}

\section{Les preuves}\label{section_preuves_chapitre_exposant_donne_oaimmvobb}

Dans cette section, on prouve les résultats énoncés sans démonstration dans la section \ref{section_resultats_sev_exposant_donne_aefoubv} précédente. On démontre aussi le théorème \ref{theoreme_spectre_amoeribnefaomnv}. \\

Soit 

\[\beta\in\left[1+\frac 1{2\ell}+\sqrt{1+\frac{1}{4\ell^2}},+\infty\right[.\]

On reprend les notations de la section \ref{section_resultats_sev_exposant_donne_aefoubv}. \\

Commençons par démontrer le lemme \ref{lemme_condition_dirr_sur_A_aomerbvozb}, qui affirme que le sous-espace $A$ construit dans la section \ref{section_resultats_sev_exposant_donne_aefoubv} est un sous-espace vectoriel $(\ell,1)$-irrationnel de $\R^n$.

\begin{preuve}
On réindexe pour plus de clarté les $\xi_{i,j}$ pour $i,j\in\{1,\ldots,\ell\}$, en $\xi_1,\ldots,\xi_{\ell^2}$ par ordre lexicographique. On réindexe de la même façon les suites $(e_k^{(i,j)})_{k\in\N}$ en $(e_k^{(1)})_{k\in\N},\ldots,(e_k^{(\ell^2)})_{k\in\N}$. \\

Montrons qu'on peut choisir des suites $(e_k^{(1)})_{k\in\N},\ldots,(e_k^{(\ell^2)})_{k\in\N}$ telles que $\xi_1,\ldots,\xi_{\ell^2}$ soient algébriquement indépendants sur $\Q$. On raisonne pour cela par récurrence finie sur $t\in\{1,\ldots,\ell^2\}$. \\

L'exposant d'irrationalité de $\xi_1$ est minoré par $\alpha>2$ (et même égal à $\alpha$ d'après la section 8 de \cite{levesley06}), donc d'après le théorème de Roth (voir \cite{roth55} page 2), $\xi_1$ est transcendant. \\ 

Soit $t\in\{1,\ldots,\ell^2-1\}$ ; supposons que les réels $\xi_1,\ldots,\xi_t$ sont algébriquement indépendants sur $\Q$. L'ensemble des réels algébriques sur $\Q(\xi_1,\ldots,\xi_t)$ est dénombrable, tandis que l'ensemble des suites $(e_k^{(t+1)})_{k\in\N}$ ne l'est pas. On peut donc choisir une suite $(e_k^{(t+1)})_{k\in\N}$ à valeurs dans $\{1,2\}$ ou $\{2\ell,2\ell+1\}$ (selon que $t$ correspond à un couple $(i,j)$ avec $i\ne j$ ou $i=j$) telle que $\xi_1,\ldots,\xi_{t+1}$ soient algébriquement indépendants sur $\Q$, ce qui conclut la récurrence. \\

Vérifions que le sous-espace $A$ ainsi construit appartient à $\mathfrak I_n(\ell,\ell)_1$. Supposons par l'absurde qu'il existe un sous-espace rationnel $B$ de $\R^{2\ell}$ de dimension $\ell$ intersectant non trivialement $A$. Soit $M_B$ une matrice dont les colonnes forment une base rationnelle de $B$. Comme $A\cap B\ne\{0\}$, on a 

\[\det\begin{pmatrix} M_A & M_B\end{pmatrix}=0\]

où $M_A$ est la matrice définie à l'équation \eqref{def_MA_matricedeA_armoeigaoinv}. Comme $M_B\in\M_{2\ell,\ell}(\Q)$, en développant ce déterminant par rapport aux $\ell$ premières colonnes grâce à un développement de Laplace (corollaire \ref{corollaire_dev_Laplace_amiovomvbnainv}), on obtient un polynôme $P\in\Q[X_1,\ldots,X_{\ell^2}]$ tel que 

\[\det\begin{pmatrix} M_A & M_B\end{pmatrix}=P(\xi_1,\ldots,\xi_{\ell^2})=0.\]

Or $\xi_1,\ldots,\xi_{\ell^2}$ sont algébriquement indépendants sur $\Q$, donc $P=0$. \\

Comme

\[M_A=\begin{pmatrix} I_\ell \\ M_\xi\end{pmatrix},\]

en décomposant $M_B\in\M_{2\ell,\ell}(\Q)$ sous la forme

\[M_B=\begin{pmatrix} B_1\\ B_2\end{pmatrix}\]

avec $B_1,B_2\in\M_{\ell,\ell}(\R)$, l'égalité $P=0$ montre que 

\begin{equation}\label{detegal0pourtoutQ_agroianoivnoizc}
\forall Q\in\M_\ell(\R),\quad \Delta_Q=\det\begin{pmatrix} I_\ell & B_1 \\ Q & B_2\end{pmatrix}=0.
\vspace{2.6mm}
\end{equation}

Énonçons un lemme sur les déterminants qui sera utile dans la suite de cette preuve. 

\begin{lemme}\label{lemme_det_par_blocs_moaiboevodizvb}

Soient $A_1,A_2,A_3,A_4\in\M_\ell(\R)$ telles que $A_1A_2=A_2A_1$. Alors 

\[\det\begin{pmatrix} A_1 & A_2\\A_3&A_4\end{pmatrix}=\det(A_4A_1-A_3A_2).\]

\end{lemme}

\begin{preuve}[Lemme \ref{lemme_det_par_blocs_moaiboevodizvb}]
Supposons que $A_1$ est inversible. On peut alors remarquer que

\[\begin{pmatrix} A_1 & A_2\\A_3&A_4\end{pmatrix}\begin{pmatrix} I_\ell&-A_2A_1^{-1}\\0&I_\ell\end{pmatrix}=\begin{pmatrix}A_1&-A_1A_2A_1^{-1}+A_2\\A_3&-A_3A_2A_1^{-1}+A_4\end{pmatrix}=\begin{pmatrix} A_1&0\\A_3&-A_3A_2A_1^{-1}+A_4\end{pmatrix}\]

car $A_1$ et $A_2$ commutent, donc

\begin{equation}\label{equa_lemme_matrices_carres_abmoifnbaoinv}
\det\begin{pmatrix} A_1 & A_2\\A_3&A_4\end{pmatrix}=\det(-A_3A_2A_1^{-1}+A_4)\det(A_1)=\det(A_4A_1-A_3A_2).
\vspace{2.6mm}
\end{equation}

Si on ne suppose plus $A_1$ inversible, par densité de $\GL_{\ell}(\R)$ dans $\M_\ell(\R)$, on peut trouver une suite de matrices $(\tilde A_N)_{N\in\N}$ de $\GL_\ell(\R)$ qui converge vers $A_1$ et commutant avec $A_2$. Alors par continuité du déterminant, l'égalité \eqref{equa_lemme_matrices_carres_abmoifnbaoinv} est encore vraie.
\end{preuve}

Comme $I_\ell$ commute avec $B_1$, on peut appliquer le lemme \ref{lemme_det_par_blocs_moaiboevodizvb} ci-dessus, pour obtenir

\begin{equation}\label{Delta_Q_apibpinfvpinpvinedc}
\forall Q\in\M_\ell(\R),\quad \Delta_Q=\det(B_2-QB_1)=0.
\vspace{2.6mm}
\end{equation}

Soit $\lambda\in\R$. En prenant $Q=\lambda I_\ell$, on a $\det(B_2-\lambda B_1)=0$. Supposons par l'absurde que $B_1$ est inversible, on a alors

\[0=\Delta_Q=\det((B_2B_1^{-1}-\lambda I_\ell)B_1)=\det(B_2B_1^{-1}-\lambda I_\ell)\det(B_1)\]

donc car $\det(B_1)\ne0$,

\[\forall \lambda\in\R,\quad \det(B_2B_1^{-1}-\lambda I_\ell)=0.\]

On a montré que pour tout $\lambda\in\R$, $\lambda$ est valeur propre de $B_2B_1^{-1}$ ce qui est absurde. \\

Ainsi, $\det(B_1)=0$. Notons

\[r=\rg(B_1)<\ell.\]

Soient $U,V\in\GL_\ell(\R)$ deux matrices inversibles telles que

\[UB_1V=\begin{pmatrix} I_r & 0 \\ 0 & 0\end{pmatrix}=\begin{pmatrix} J_r& 0\end{pmatrix}\in\M_\ell(\R),\]

où on a noté

\[J_r=\begin{pmatrix} I_r\\0\end{pmatrix}\in\M_{\ell,r}(\R).\]

On décompose

\[UB_2V=\begin{pmatrix} C_1& C_2\end{pmatrix}\in\M_\ell(\R)\]

où la matrice $C_1\in\M_{\ell,r}(\R)$ est constituée des $r$ premières colonnes de $UB_2V$, et la matrice $C_2\in\M_{\ell,\ell-r}(\R)$ est constituée des $\ell-r$ dernières colonnes de $UB_2V$. Ainsi, on a des matrices équivalentes :

\begin{equation}\label{equivdeuxmatrices_moareinmoaunvomfn}
\begin{pmatrix} J_r & 0 \\ C_1 & C_2 \end{pmatrix}=\begin{pmatrix} U & 0 \\ 0 & U\end{pmatrix}
\begin{pmatrix} B_1\\B_2\end{pmatrix}V\in\M_{2\ell,\ell}(\R).
\vspace{2.6mm}
\end{equation}

Comme $I_\ell$ commute avec $UB_1V$, on a d'après le lemme \ref{lemme_det_par_blocs_moaiboevodizvb} que pour tout \hbox{$Q\in\M_\ell(\R)$} :
\begin{align*}
\det\begin{pmatrix} I_\ell & UB_1V \\ UQU^{-1} & UB_2V\end{pmatrix}
	&=\det(UB_2V-UQU^{-1}UB_1V) \\
	&=\det(U(B_2-QB_1)V) \\
	&=\det(U)\Delta_Q\det (V) \\
	&=0
\end{align*}
d'après \eqref{Delta_Q_apibpinfvpinpvinedc}. Ceci étant vrai pour tout $Q\in\M_\ell(\R)$, quitte à poser $Q'=UQU^{-1}$, on a aussi

\begin{equation}\label{detQprimeegal0_aomibneofinv}
\forall Q'\in\M_\ell(\R),\quad \Delta'_{Q'}=\det\begin{pmatrix} I_\ell & UB_1V \\ Q' & UB_2V\end{pmatrix}=0.
\vspace{2.6mm}
\end{equation}

Soit $R\in\M_{\ell,r}(\R)$. Définissons une matrice $Q'$ par blocs :

\[Q'=\begin{pmatrix} C_1-R & 0\end{pmatrix}\in\M_\ell(\R).\]

On obtient, encore d'après le lemme \ref{lemme_det_par_blocs_moaiboevodizvb} car $I_\ell$ et $UB_1V$ commutent, que 
\begin{align*}
0
	&=\Delta'_{Q'} \\
	&=\det(UB_2V-Q'UB_1V) \\
	&=\det\left(\begin{pmatrix} C_1 & C_2\end{pmatrix}-\begin{pmatrix} C_1-R & 0\end{pmatrix}\begin{pmatrix} I_r & 0\\0&0\end{pmatrix}\right) \\
	&=\det\begin{pmatrix} C_1-C_1+R & C_2\end{pmatrix} \\
	&=\det\begin{pmatrix} R & C_2\end{pmatrix}.
\end{align*}
Si on avait $\rg(C_2)=\ell-r$, alors d'après le théorème de la base incomplète, on pourrait trouver $R$ telle que $\rg\begin{pmatrix} R & C_2\end{pmatrix}=\ell$, ce qui est absurde car son déterminant $\Delta'_{Q'}$ serait alors non nul. Donc $\rg(C_2)<\ell-r$. Or avec \eqref{equivdeuxmatrices_moareinmoaunvomfn}, on en déduit

\[\rg(M_B)=\rg\begin{pmatrix} B_1 \\ B_2\end{pmatrix}=\rg\begin{pmatrix} J_r & 0 \\ C_1 & C_2 \end{pmatrix}=r+\rg(C_2)<r+\ell-r=\ell,\]

ce qui est absurde car $\dim B=\ell=\rg(M_B)$. \\

Finalement, $A$ intersecte trivialement tous les sous-espaces rationnels de $\R^{2\ell}$ de dimension $\ell$, c'est-à-dire $A\in\mathfrak I_{2\ell}(\ell,\ell)_1$. Cela termine la preuve du lemme \ref{lemme_condition_dirr_sur_A_aomerbvozb}. 
\end{preuve}

On fixe désormais dans toute la suite, des suites $(e_k^{(i,j)})_{k\in\N}$ vérifiant le lemme \ref{lemme_condition_dirr_sur_A_aomerbvozb}. \\

Avant de démontrer le lemme \ref{majoration_psi2_A_BN_amroghozv}, montrons un lemme qui majore la hauteur des $B_N$. On verra plus loin (grâce au lemme \ref{lemmeminorationhauteurdesBN_aeoinmoaeivnv}) que cette majoration est optimale, à constante multiplicative près. 

\begin{lemme}\label{majoration_hauteur_B_N_painoaeoefainveo}

Pour tout $N\ge 1$, on a 

\[H(B_N)\le c_2(\theta^{\lfloor\alpha^N\rfloor})^\ell\]

où $c_2>0$ dépend uniquement de $\ell$.

\end{lemme}

\begin{preuve}
On a

\begin{equation}\label{maj_f_N_ij_aouomfaeuboueabv}
\abs{f_N^{(i,j)}}\le (2\ell+1)\theta^{\lfloor\alpha^N\rfloor}\sum_{k=0}^N\frac 1{\theta^{\lfloor\alpha^k\rfloor}}\le 2(2\ell+1)\theta^{\lfloor\alpha^N\rfloor}
\vspace{2.6mm}
\end{equation}

car la minoration $\theta\ge 2$ donne

\[\sum_{k=0}^N\frac 1{\theta^{\lfloor\alpha^k\rfloor}}\le \sum_{k=0}^\infty \frac 1{\theta^k}=\frac{\theta}{\theta-1}\le 2.\]

On a donc
\begin{align*}
H(B_N)
	&\le \norme{X_N^{(1)}\wedge\cdots\wedge X_N^{(\ell)}} \\
	&\le \prod_{j=1}^\ell\norme{X_N^{(j)}} \\
	&\le (2(2\ell+1)\cdot \sqrt{2\ell})^\ell (\theta^{\lfloor\alpha^N\rfloor})^\ell
\end{align*}
car les $2\ell$ coefficients de chaque $X_N^{(j)}$ sont inférieurs à $2(2\ell+1)\cdot \theta^{\lfloor\alpha^N\rfloor}$.
\end{preuve}

Démontrons maintenant le lemme \ref{majoration_psi2_A_BN_amroghozv}, selon lequel 

\[\psi_\ell(A,B_N)\le \frac{c_1}{H(B_N)^{\alpha/\ell}}\]

pour tout $N\ge 1$, avec une constante $c_1>0$ qui dépend uniquement de $A$.

\begin{preuve}
Pour $i\in\{1,\ldots,\ell\}$, notons $Z_N^{(i)}=\theta^{-\lfloor\alpha^N\rfloor}X_N^{(i)}$. On a $\norme{Y_i}\ge 1$, ce qui donne en appliquant le lemme \ref{lemme_maj_psiXY_amoignaoivnovibs} que

\begin{equation}\label{majoration_psi_XNYi_zmofihgmenv}
\psi(X_N^{(i)},Y_i)=\psi(Z_N^{(i)},Y_i) \le \frac{\norme{Z_N^{(i)}-Y_i}}{\norme{Y_i}}\le c_3(\theta^{\lfloor\alpha^N\rfloor})^{-\alpha}
\vspace{2.6mm}
\end{equation}

avec $c_3>0$ ne dépendant que de $\ell$, en utilisant \eqref{le_sev_approche_bien_moaeirgnomdfvn}. D'après la proposition \ref{resultat_reconstruction_proximite_oaeirgbvodbv}, il existe des constantes $c_4,c_5>0$ ne dépendant que de $A$, telles que

\begin{equation}\label{maj_psi_ell_A_B_N_aoeribqfoisbcovbz}
\psi_\ell(A,B_N)\le c_4\sum_{i=1}^\ell \psi(X_N^{(i)},Y_i)\le \frac{c_5}{(\theta^{\lfloor\alpha^N\rfloor})^{\alpha}}.
\vspace{2.6mm}
\end{equation}

D'après le lemme \ref{majoration_hauteur_B_N_painoaeoefainveo} qui donne $H(B_N)\le c_2(\theta^{\lfloor\alpha^N\rfloor})^\ell$, on a donc

\[\psi_\ell(A,B_N) \le \frac{c_1}{H(B_N)^{\alpha/\ell}}\]

avec $c_1>0$ ne dépendant que de $A$. Ceci termine la preuve du lemme \ref{majoration_psi2_A_BN_amroghozv}.
\end{preuve}

Montrons maintenant que les $B_N$ sont les seuls sous-espaces rationnels de dimension $\ell$ approchant aussi bien $A$ (\emph{i.e.} à l'exposant $\alpha/\ell$). Ceci montrera que le meilleur exposant possible est déterminé par les $B_N$. Précisément, on démontre ici le lemme \ref{les_BN_sont_les_meilleurs_afnemobibvo} : si $\epsilon>0$ et $C\in\mathfrak R_{2\ell}(\ell)$ sont tels que

\[\psi_\ell(A,C)\le \frac 1{{H(C)}^{\alpha/\ell+\epsilon}}\]

et si $H(C)$ est suffisamment grand (en fonction de $\ell$ et de $\epsilon$), alors il existe $N\ge 1$ tel que $C=B_N$.

\begin{preuve}
Soient $\epsilon>0$ et $C\in\mathfrak R_{2\ell}(\ell)$ tels que

\begin{equation}\label{hyp_sur_psi2_A_C}
\psi_\ell(A,C)\le \frac 1{H(C)^{\alpha/\ell+\epsilon}}.
\vspace{2.6mm}
\end{equation}

Montrons que si $H(C)$ est suffisamment grand, alors il existe un entier $N\ge 1$ tel que $C=B_N$. Commençons par travailler avec $N\in\N^*$ quelconque, on le fixera par la suite. \\

Le sous-espace $C$ est un sous-espace rationnel de dimension $\ell$ : il existe donc \hbox{$v_1,\ldots,v_\ell\in\Z^{2\ell}$} tels que $(v_1,\ldots,v_\ell)$ soit une $\Z$-base de $C\cap\Z^{2\ell}$. On a alors

\[H(C)=\norme{v_1\wedge\cdots\wedge v_\ell}.\]

Pour montrer que $C=B_N$ pour un certain entier $N$, montrons que tous les $X_N^{(i)}$ pour $i\in\{1,\ldots,\ell\}$ sont dans $C=\vect(v_1,\ldots,v_\ell)$. Par égalité des dimensions, on pourra conclure que $C=B_N$. Soit $i\in\{1,\ldots,\ell\}$ ; considérons les $\ell+1$ vecteurs $X_N^{(i)},v_1,\ldots,v_\ell$, et posons

\[Q=\begin{pmatrix} X_N^{(i)}&v_1&\cdots&v_\ell\end{pmatrix}\in\M_{2\ell,\ell+1}(\Z).\]

Comme la famille $(v_1,\ldots,v_\ell)$ est libre, pour montrer que $X_N^{(i)}\in\vect(v_1,\ldots,v_\ell)$, il suffit de montrer que $\rg(Q)<\ell+1$, \emph{i.e.} que tous les mineurs $(\ell+1)\times(\ell+1)$ de $Q$ sont nuls. Pour cela, on montre que 

\[D=\norme{X_N^{(i)}\wedge v_1\wedge\cdots\wedge v_\ell}=0.\]

Notons $p_C^\perp$ la projection orthogonale sur le sous-espace $C$ (les notations introduites ici sont illustrées sur la figure \ref{projection_sur_C}), et $h$ le vecteur

\[h=p_C^\perp(X_N^{(i)})-X_N^{(i)}.\]

Alors il existe $\lambda_1,\ldots,\lambda_\ell\in\R$ tels que le vecteur $X_N^{(i)}$ s'écrive sous la forme

\[X_N^{(i)}=\sum_{j=1}^\ell\lambda_jv_j-h.\]

\begin{figure}[H]
\begin{center}
\includegraphics[scale=1.1]{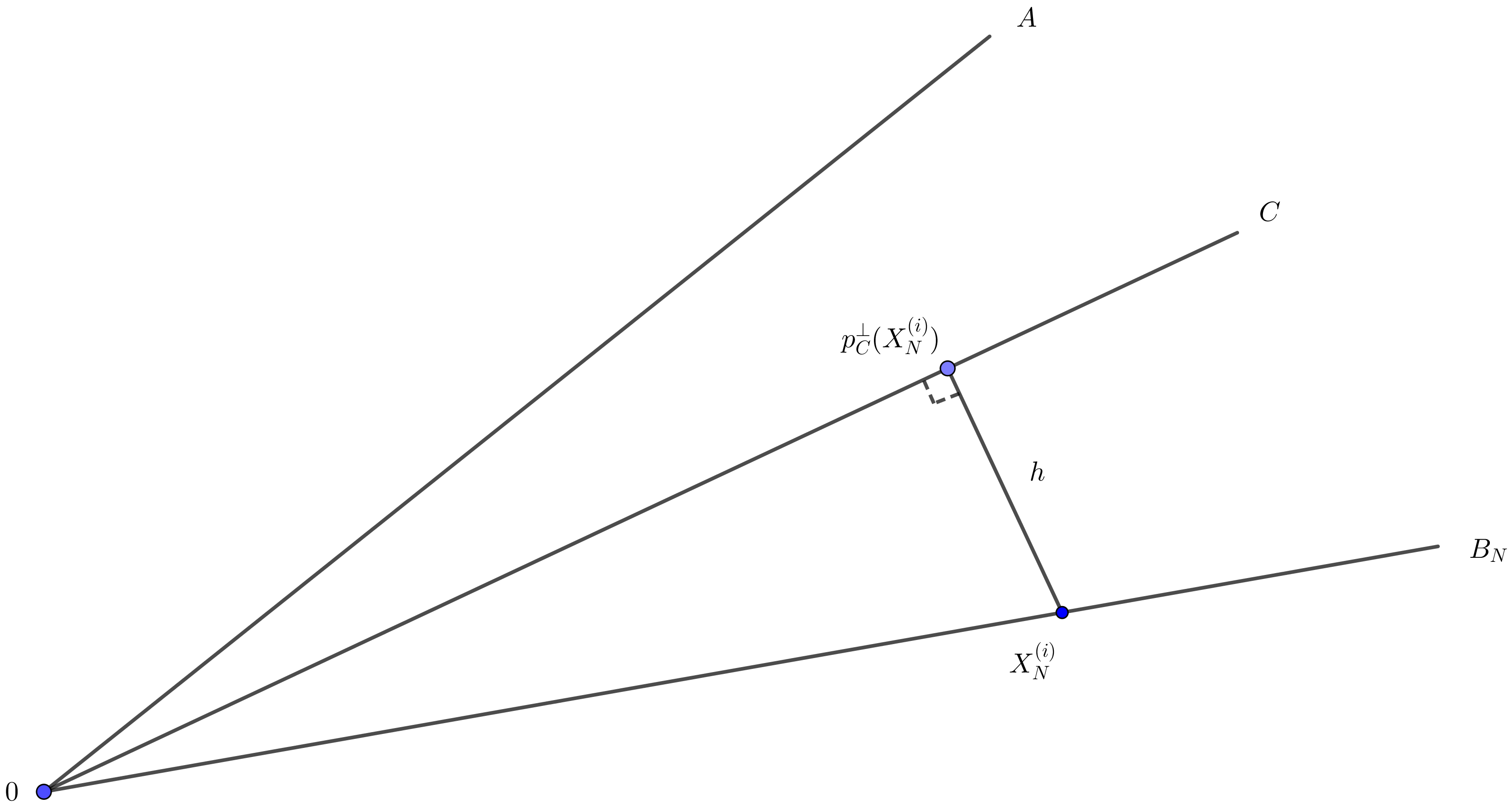}
\caption{Projection de $X_N^{(i)}$ sur $C$}\label{projection_sur_C}
\end{center}
\end{figure}

Or $h\in C^\perp$, donc
\begin{align*}
D
	&=\norme{X_N^{(i)}\wedge v_1\wedge\cdots\wedge v_\ell} \\
	&=\norme{\left(\sum_{j=1}^\ell\lambda_jv_j-h\right)\wedge v_1\wedge\cdots\wedge v_\ell} \\
	&=\norme h \cdot \norme{v_1\wedge\cdots\wedge v_\ell}\\
	&=\norme h H(C).
\end{align*}
De plus d'après le lemme \ref{lemme_proximite_proj_ortho_apamoebnamon}, on a

\[\norme h=\norme{X_N^{(i)}}\psi_1(X_N^{(i)},C),\]

donc 

\[D= \norme{X_N^{(i)}}\psi_1(X_N^{(i)},C)H(C)\le c_6 \theta^{\lfloor\alpha^N\rfloor} (\psi(X_N^{(i)},Y_i)+\psi_1(\vect(Y_i),C))H(C)\]

avec $c_6>0$ ne dépendant que de $\ell$, d'après l'équation \eqref{maj_f_N_ij_aouomfaeuboueabv} et en utilisant la proposition \ref{inegalite_triang_angles_aeoifgnivn}. \\ 

Or d'après le lemme \ref{lemme_transfert_psi1_psiell_baoeoearv}, on a

\[\psi_1(\vect(Y_i),C)\le \psi_\ell(A,C).\]

En utilisant alors la majoration \eqref{majoration_psi_XNYi_zmofihgmenv} de la preuve du lemme \ref{majoration_psi2_A_BN_amroghozv} pour majorer $\psi(Y_i,X_N^{(i)})$, ainsi que la majoration \eqref{hyp_sur_psi2_A_C} pour majorer $\psi_\ell(A,C)$, on trouve 

\begin{equation}\label{majoration_de_D_aobvobobdvipsn}
D\le c_6 \theta^{\lfloor\alpha^N\rfloor}H(C)\left(\frac {c_3}{(\theta^{\lfloor\alpha^N\rfloor})^{\alpha}}+\frac 1{H(C)^{\alpha/\ell+\epsilon}}\right)\le c_7 \left(\frac {H(C)}{\theta^{\lfloor\alpha^N\rfloor(\alpha-1)}}+\frac{\theta^{\lfloor\alpha^N\rfloor}}{H(C)^{\alpha/\ell-1+\epsilon}}\right)
\vspace{2.6mm}
\end{equation}

où $c_7>0$ ne dépend que de $\ell$. \\

On va désormais choisir $N$ : notons $N$ le plus grand entier tel que

\[\theta^{\alpha^N}\le H(C)^{\alpha/\ell-1+\epsilon/2}.\]

Alors

\begin{equation}\label{majoration_de_10puissancealphaN}
\theta^{\lfloor\alpha^N\rfloor}\le H(C)^{\alpha/\ell-1+\epsilon/2}.
\vspace{2.6mm}
\end{equation}

De plus, par maximalité de $N$, on a

\[\left(\theta^{\alpha^N}\right)^\alpha=\theta^{\alpha^{N+1}}> H(C)^{\alpha/\ell-1+\epsilon/2},\]

donc

\[\theta^{\alpha^N}>H(C)^{(\alpha/\ell-1+\epsilon/2)/\alpha}.\]

Or

\[\alpha=\ell\beta\ge\ell+\frac 12+\sqrt{\ell^2+\frac 14}=\frac{2\ell+1+\sqrt{1+4\ell^2}}2,\]

donc

\[\alpha^2-(2\ell+1)\alpha+\ell\ge0,\]

soit

\[(\alpha-\ell)(\alpha-1)\ge\ell\alpha,\]

donc

\[\frac{\alpha/\ell-1}{\alpha}\ge\frac 1{\alpha-1}.\]

Ainsi, 

\[\theta^{\alpha^N}>H(C)^{1/(\alpha-1)+\epsilon/(2\alpha)},\]

d'où

\begin{equation}\label{minoration_de_10puissancealphaN}
\theta^{\lfloor\alpha^N\rfloor}>\frac 1{\theta} H(C)^{1/(\alpha-1)+\epsilon/(2\alpha)}.
\vspace{2.6mm}
\end{equation}

En utilisant la majoration \eqref{majoration_de_10puissancealphaN} et la minoration \eqref{minoration_de_10puissancealphaN}, revenons à la majoration \eqref{majoration_de_D_aobvobobdvipsn} de $D$ : 
\begin{align*}
D
	&\le c_7 \left(\frac {H(C)}{\theta^{\lfloor\alpha^N\rfloor(\alpha-1)}}+\frac{\theta^{\lfloor\alpha^N\rfloor}}{H(C)^{\alpha/\ell-1+\epsilon}}\right) \\
	&\le c_7 \left(\frac{c_8 H(C)}{H(C)^{(1/(\alpha-1)+\epsilon/(2\alpha))(\alpha-1)}}+\frac{H(C)^{\alpha/\ell-1+\epsilon/2}}{H(C)^{\alpha/\ell-1+\epsilon}}\right) \\
	&\le c_9\left(\frac 1{H(C)^{(\alpha-1)\epsilon/(2\alpha)}}+\frac 1{H(C)^{\epsilon/2}}\right) \\
	&\xrightarrow[H(C)\to+\infty]{} 0
\end{align*}
avec $c_8,c_9>0$ dépendant uniquement de $\ell$. \\

Si $H(C)$ est suffisamment grand (en fonction de $\ell$ et de $\epsilon$, puisque $(\alpha-1)/(2\alpha)\ge(\ell-1)/(2\ell)$), on a $D<1$. Or $E=\norme{X_N^{(i)}\wedge v_1\wedge\cdots\wedge v_\ell}_\infty$ est un entier naturel tel que $E\le D$, donc $E=0$, c'est-à-dire

\[X_N^{(i)}\wedge v_1\wedge\cdots\wedge v_\ell=0.\]

Finalement, si $H(C)$ est suffisamment grand, on a montré que $C=B_N$ pour $N$ le plus grand entier tel que $\theta^{\alpha^N}\le H(C)^{\alpha/\ell-1+\epsilon/2}$. Ceci termine la preuve du lemme \ref{majoration_psi2_A_BN_amroghozv}. 
\end{preuve}

Avant de démontrer le lemme \ref{minorationdupsielldeAetBN_eogihrgoeihgoi}, commençons par minorer la hauteur des $B_N$ (ce qui montre que le lemme \ref{majoration_hauteur_B_N_painoaeoefainveo} est optimal à constante multiplicative près).

\begin{lemme}\label{lemmeminorationhauteurdesBN_aeoinmoaeivnv}

Pour tout $N$ suffisamment grand, on a

\[H(B_N)\ge \tilde c(\theta^{\lfloor\alpha^N\rfloor})^\ell\]

avec $\tilde c>0$ ne dépendant que de $A$.
\end{lemme}

\begin{preuve}[Lemme \ref{lemmeminorationhauteurdesBN_aeoinmoaeivnv}]
Soit $N\ge 1$. Pour montrer la minoration voulue sur la hauteur de $B_N$, on va d'abord montrer que la famille $(X_N^{(1)},\ldots,X_N^{(\ell)})$ est une $\Z$-base de $B_N\cap\Z^{2\ell}$. Pour cela, notons $P$ le parallélotope engendré par $X_N^{(1)},\ldots,X_N^{(\ell)}$, \emph{i.e.}

\[P=\left\{\sum_{i=1}^\ell \lambda_iX_N^{(i)},\ (\lambda_1,\ldots,\lambda_\ell)\in[0,1]^\ell\right\},\]

et montrons que les $2^\ell$ sommets de $P$ sont ses seuls points entiers. Notons $\mathcal S$ l'ensemble des $2^\ell$ sommets de $P$, \emph{i.e.}

\[\mathcal S=\left\{\sum_{i=1}^\ell \delta_iX_N^{(i)},\ (\delta_1,\ldots,\delta_\ell)\in\{0,1\}^\ell\right\}.\]

Supposons par l'absurde qu'il existe 

\[X\in (P\setminus\mathcal S)\cap\Z^{2\ell}.\]

On a alors $(\lambda_1,\ldots,\lambda_\ell)\in[0,1]^\ell\setminus\{0,1\}^\ell$ tel que 

\[X=\sum_{i=1}^\ell \lambda_i X_N^{(i)}\in\Z^{2\ell}.\]

Les $\ell$ premières coordonnées de $X$ donnent

\[\forall i\in\{1,\ldots,\ell\},\quad \lambda_i \theta^{\lfloor\alpha^N\rfloor}\in\Z.\]

Il existe donc des entiers $\gamma_1,\ldots,\gamma_\ell\in\{0,\ldots,\theta^{\lfloor\alpha^N\rfloor}\}$ tels que

\[\forall i\in\{1,\ldots,\ell\},\quad \lambda_i=\frac{\gamma_i}{\theta^{\lfloor\alpha^N\rfloor}}\]

car les $\lambda_i$ sont dans $[0,1]$ pour tout $i\in\{1,\ldots,\ell\}$. De plus, les $\ell$ dernières coordonnées de $X$ donnent

\[\forall i\in\{1,\ldots,\ell\},\quad \sum_{j=1}^\ell \lambda_j f_N^{(i,j)}\in\Z,\]

soit

\begin{equation}\label{random_riozagoeigbzouefm}\forall i\in\{1,\ldots,\ell\},\quad \sum_{j=1}^\ell \frac{\gamma_j}{\theta^{\lfloor\alpha^N\rfloor}}\cdot\theta^{\lfloor\alpha^N\rfloor} \sum_{k=0}^N\frac{e_k^{(i,j)}}{\theta^{\lfloor\alpha^k\rfloor}}=\sum_{k=0}^N\frac 1{\theta^{\lfloor\alpha^k\rfloor}}\sum_{j=1}^\ell\gamma_j e_k^{(i,j)}\in\Z.
\vspace{2.6mm}
\end{equation}

Pour $k\in\{0,\ldots,N\}$, notons $E_k$ la matrice $(e_k^{(i,j)})_{(i,j)\in\{1,\ldots,\ell\}^2}\in\M_{\ell}(\Z)$, et $\Gamma$ le vecteur colonne $\transp(\gamma_1,\ldots,\gamma_\ell)$. Ainsi, on peut réécrire sous forme matricielle les $\ell$ équations données par \eqref{random_riozagoeigbzouefm} :

\[\begin{pmatrix} \displaystyle\sum_{k=0}^N \frac{e_k^{(1,1)}}{\theta^{\lfloor\alpha^k\rfloor}} & \cdots & \displaystyle\sum_{k=0}^N \frac{e_k^{(1,\ell)}}{\theta^{\lfloor\alpha^k\rfloor}} \\ \vdots & & \vdots \\ \displaystyle\sum_{k=0}^N \frac{e_k^{(\ell,1)}}{\theta^{\lfloor\alpha^k\rfloor}} & \cdots & \displaystyle\sum_{k=0}^N \frac{e_k^{(\ell,\ell)}}{\theta^{\lfloor\alpha^k\rfloor}}\end{pmatrix}\begin{pmatrix} \gamma_1\\\vdots\\\vdots\\\gamma_\ell\end{pmatrix}\in\Z^\ell,\]

soit

\[\sum_{k=0}^N \frac{1}{\theta^{\lfloor\alpha^k\rfloor}} E_k\Gamma\in\Z^\ell,\]

ou encore, car on va s'intéresser plus particulièrement au dernier terme de la somme,

\[\sum_{k=0}^{N-1} \theta^{\lfloor\alpha^N\rfloor-\lfloor\alpha^k\rfloor} E_k\Gamma+E_N\Gamma\in\theta^{\lfloor\alpha^N\rfloor}\Z^\ell.\]

Comme $E_N\in\M_\ell(\Z)$, la transposée de sa comatrice appartient aussi à $\M_\ell(\Z)$, et donc

\[\sum_{k=0}^{N-1} \theta^{\lfloor\alpha^N\rfloor-\lfloor\alpha^k\rfloor} \transp\com(E_N)E_k\Gamma+\transp\com(E_N)E_N\Gamma\in\theta^{\lfloor\alpha^N\rfloor}\Z^\ell,\]

soit

\[\sum_{k=0}^{N-1} \theta^{\lfloor\alpha^N\rfloor-\lfloor\alpha^k\rfloor} \transp\com(E_N)E_k\Gamma+\det(E_N)\Gamma\in\theta^{\lfloor\alpha^N\rfloor}\Z^\ell.\]

Pour $i\in\{1,\ldots\ell\}$ et $k\in\{0,\ldots,N-1\}$, notons $L_{k,i}\in\M_{1,\ell}(\Z)$ la $i$-ème ligne du produit $\transp\com(E_N)E_k$. Ainsi,

\begin{equation}\label{random_omerganmoernbme}
\forall i\in\{1,\ldots,\ell\},\quad \sum_{k=0}^{N-1} (L_{k,i}\Gamma)\theta^{\lfloor\alpha^N\rfloor-\lfloor\alpha^k\rfloor}+\det(E_N)\gamma_i\in\theta^{\lfloor\alpha^N\rfloor}\Z.
\vspace{2.6mm}
\end{equation}

Soit $j\in\{1,\ldots,\ell\}$. On peut commencer par remarquer que 

\[e_N^{(j,j)}\ge2\ell>2(\ell-1)\ge\sum_{i\ne j} e_N^{(i,j)},\]

donc $E_N$ est une matrice à diagonale strictement dominante (le cas $\ell=1$ étant trivial car la somme de droite est alors vide -- la matrice $E_N$ est une matrice de $\M_1(\Z)$ dans ce cas). Elle est donc inversible, \emph{i.e.} 

\begin{equation}\label{random_iorembgzugmbez}
\det(E_N)\ne 0.
\vspace{2.6mm}
\end{equation}

De plus, on a $\abs{e_N^{(i,j)}}\le 2$ si $i\ne j$ et $\abs{e_N^{(i,i)}}\le 2\ell+1$. Donc

\[\abs{\det(E_N)}=\abs{\sum_{\sigma\in\mathfrak S_\ell}\epsilon(\sigma)\prod_{i=1}^\ell e_N^{(i,\sigma(i))}} \le\sum_{\sigma\in\mathfrak S_\ell}\prod_{i=1}^\ell \abs{e_N^{(i,\sigma(i))}} \le\ell!\,(2\ell+1)^\ell <\theta\]

par définition de $\theta$ (qui est le plus petit nombre premier vérifiant la minoration \eqref{hyp_sur_theta_amorinvamnoavcbo}). Ainsi, comme $0<\abs{\det(E_N)}<\theta$, on a $v_\theta(\det(E_N))=0$, donc

\[v_\theta(\det(E_N)\gamma_i)=v_\theta(\det(E_N))+v_\theta(\gamma_i)=v_\theta(\gamma_i).\]

On définit $u\ge 0$ et $i_0\in\{1,\ldots,\ell\}$ tels que

\[u=\min(v_\theta(\gamma_1),\ldots,v_\theta(\gamma_\ell))=v_\theta(\gamma_{i_0}).\]

Revenons alors à l'équation \eqref{random_omerganmoernbme}, qui s'écrit

\[\forall i\in\{1,\ldots,\ell\},\quad v_\theta\left(\sum_{k=0}^{N-1} (L_{k,i}\Gamma)\theta^{\lfloor\alpha^N\rfloor-\lfloor\alpha^k\rfloor}+\det(E_N)\gamma_i\right)\ge\lfloor\alpha^N\rfloor,\]

avec par convention $v_\theta(0)=+\infty$. \\

On a pour tout $k\in\{0,\ldots,N-1\}$,

\[v_\theta(\theta^{\lfloor\alpha^N\rfloor-\lfloor\alpha^k\rfloor})\ge \lfloor\alpha^N\rfloor-\lfloor\alpha^{N-1}\rfloor>0,\]

et $L_{k,i}\Gamma$ est une combinaison $\Z$-linéaire des $\gamma_i$, donc

\[v_\theta(L_{k,i}\Gamma)\ge \min(v_\theta(\gamma_1),\ldots,v_\theta(\gamma_\ell))=u.\]

Or si $v_\theta(a)\ne v_\theta(b)$, alors $v_\theta(a+b)=\min(v_\theta(a),v_\theta(b))$. Donc pour $i=i_0$, on a

\[v_\theta\left(\sum_{k=0}^{N-1} (L_{k,i_0}\Gamma)\theta^{\lfloor\alpha^N\rfloor-\lfloor\alpha^k\rfloor}+\det(E_N)\gamma_{i_0}\right)=u\ge\lfloor\alpha^N\rfloor\]

car la valuation $\theta$-adique du premier terme est strictement plus grande que $u$ et que celle du second est égale à $u$. Finalement, par définition de $u$, on a

\[\forall i\in\{1,\ldots,\ell\},\quad v_\theta(\gamma_i)\ge \lfloor\alpha^N\rfloor.\]

Or tous les $\gamma_i$ sont dans $\{0,\ldots,\theta^{\lfloor\alpha^N\rfloor}\}$, donc 

\[\forall i\in\{1,\ldots,\ell\},\quad \gamma_i\in\big\{0,\theta^{\lfloor\alpha^N\rfloor}\big\},\]

soit $X\in\mathcal S$ ce qui est absurde. \\

On a donc montré que les seuls points entiers de $P$ sont les points de $\mathcal S$, ce qui montre que la famille $(X_N^{(1)},\ldots,X_N^{(\ell)})$ est une $\Z$-base de $B_N\cap\Z^{2\ell}$. \\

Ainsi, d'après la proposition \ref{defequiv_hauteur_aveclecovol_eohuflbg}, on a

\[H(B_N)=\norme{X_N^{(1)}\wedge\cdots\wedge X_N^{(\ell)}}.\]

Or

\[\forall i\in\{1,\ldots,\ell\},\quad \theta^{-\lfloor\alpha^N\rfloor} X_N^{(i)}\xrightarrow[N\to\infty]{}Y_i,\]

donc

\[\norme{\theta^{-\lfloor\alpha^N\rfloor} X_N^{(1)}\wedge\cdots\wedge\theta^{-\lfloor\alpha^N\rfloor}X_N^{(\ell)}}\xrightarrow[N\to\infty]{}\norme{Y_1\wedge\cdots\wedge Y_\ell},\]

donc pour $N$ suffisamment grand, 

\[\norme{\theta^{-\lfloor\alpha^N\rfloor} X_N^{(1)}\wedge\cdots\wedge\theta^{-\lfloor\alpha^N\rfloor}X_N^{(\ell)}}\ge \frac 12\norme{Y_1\wedge\cdots\wedge Y_\ell}.\]

Finalement, pour $N$ suffisamment grand, on a

\[H(B_N)=\left(\theta^{\lfloor\alpha^N\rfloor}\right)^\ell\norme{\theta^{-\lfloor\alpha^N\rfloor} X_N^{(1)}\wedge\cdots\wedge\theta^{-\lfloor\alpha^N\rfloor}X_N^{(\ell)}}\ge \tilde c\left(\theta^{\lfloor\alpha^N\rfloor}\right)^\ell\]

avec $\tilde c>0$ ne dépendant que de $A$, ce qui termine la preuve du lemme \ref{lemmeminorationhauteurdesBN_aeoinmoaeivnv}. 
\end{preuve}

On peut désormais démontrer le lemme \ref{minorationdupsielldeAetBN_eogihrgoeihgoi}, selon lequel 

\[\psi_\ell(A,B_N)\ge \frac{c}{H(B_N)^{\alpha/\ell}}\]

pour tout $N$ suffisamment grand, où la constante $c$ dépend seulement de $A$.

\begin{preuve}[Lemme \ref{minorationdupsielldeAetBN_eogihrgoeihgoi}]
Posons $Z_N^{(1)}=\theta^{-\lfloor\alpha^N\rfloor}X_N^{(1)}$ et notons $p_A^\perp$ la projection orthogonale sur $A$. D'après le lemme \ref{lemme_angle_projete_ortho_amotinfqovn}, on a

\[\psi_1(\vect(Z_N^{(1)}),A)=\psi(Z_N^{(1)},p_A^\perp(Z_N^{(1)})).\]

Avec le lemme \ref{lemme_transfert_psi1_psiell_baoeoearv}, on a donc

\begin{equation}\label{minoration_psiell_AB_N_aroighmeobvo}
\psi_\ell(A,B_N)\ge \psi_1(\vect(Z_N^{(1)}),A)= \psi(Z_N^{(1)},p_A^\perp(Z_N^{(1)})).
\vspace{2.6mm}
\end{equation}

Posons

\[\Delta=p_A^\perp(Z_N^{(1)})-Y_1,\]

ainsi que

\[\omega=\norme{p_A^\perp(Z_N^{(1)})-Z_N^{(1)}}.\]

On décompose $p_A^\perp(Z_N^{(1)})$ dans la base $(Y_1,\ldots,Y_\ell)$ :

\[p_A^\perp(Z_N^{(1)})=\sum_{i=1}^\ell \lambda_i Y_i=\transp \begin{pmatrix}\lambda_1&\cdots&\lambda_\ell&\star&\cdots&\star\end{pmatrix}\]

où les $\star$ sont des coefficients non précisés, car pour tout $i\in\{1,\ldots,\ell\}$, 

\[Y_i=\transp \begin{pmatrix} \delta_{i,1}&\cdots&\delta_{i,\ell}&\star&\cdots&\star\end{pmatrix}\]

en notant $\delta$ le symbole de Kronecker. On a aussi pour tout $i\in\{1,\ldots,\ell\}$, 

\[Z_N^{(i)}=\transp \begin{pmatrix} \delta_{i,1}&\cdots&\delta_{i,\ell}&\star&\cdots&\star\end{pmatrix},\]

donc

\[\omega^2=\norme{p_A^\perp(Z_N^{(1)})-Z_N^{(1)}}^2\ge (\lambda_1-1)^2+\sum_{i=2}^\ell \lambda_i^2,\]

donc

\begin{equation}\label{majdetoutpar_delta_moairnbmeafkbvon}
\begin{cases}\abs{\lambda_1-1}\le\omega\\\forall i\in\{2,\ldots,\ell\},\quad \abs{\lambda_i}\le\omega.\end{cases}
\vspace{2.6mm}
\end{equation}

Soit $j\in\{1,\ldots,\ell\}$. On a

\[\norme{Y_j}^2=1+\sum_{i=1}^\ell\left(\sum_{k=0}^\infty \frac{e_k^{(i,j)}}{\theta^{\lfloor \alpha^k\rfloor}}\right)^2.\]

Or $\alpha>2$, $\theta\ge 2$ et $2\ell+1\le \theta$ d'après la condition \eqref{hyp_sur_theta_amorinvamnoavcbo} sur $\theta$, donc

\[\sum_{k=0}^\infty \frac{e_k^{(i,j)}}{\theta^{\lfloor \alpha^k\rfloor}}\le \sum_{k=1}^\infty \frac{e_k^{(i,j)}}{\theta^k}\le \frac{2\ell+1}{\theta}\cdot\frac 1{1-1/\theta}\le \frac{\theta}{\theta-1}\le 2.\]

Donc

\begin{equation}\label{maj_norme_Yj_arogmoeubv}
\norme{Y_j}\le \sqrt{1+4\ell}=c_{10}.
\vspace{2.6mm}
\end{equation}

Or 

\[\Delta=p_A^\perp(Z_N^{(1)})-Y_1=(\lambda_1-1)Y_1+\sum_{i=2}^\ell\lambda_i Y_i,\]

donc avec les deux majorations de \eqref{majdetoutpar_delta_moairnbmeafkbvon}, on a

\begin{equation}\label{maj_norme_Delta_aoerihgofbv}
\norme{\Delta}\le c_{11}\omega
\vspace{2.6mm}
\end{equation}

avec $c_{11}>0$ ne dépendant que de $\ell$. \\

On a
\begin{align}
\norme{Z_N^{(1)}\wedge p_A^\perp(Z_N^{(1)})}\nonumber
	&=\norme{Z_N^{(1)}\wedge(Y_1+p_A^\perp(Z_N^{(1)})-Y_1)}\nonumber \\
	&=\norme{Z_N^{(1)}\wedge Y_1+Z_N^{(1)}\wedge \Delta}\nonumber \\
	&\ge \norme{Z_N^{(1)}\wedge Y_1}-\norme{Z_N^{(1)}\wedge \Delta}\label{minoration_norme_Z_N1_proj_maoeboefb}.
\end{align}

Or

\[\begin{pmatrix} Z_N^{(1)} &Y_1\end{pmatrix}\in\M_{2\ell,2}(\R),\]

donc en notant $\eta_{i,j}$ le mineur $2\times 2$ correspondant aux lignes $i$ et $j$ avec $i<j$, on a

\[\norme{Z_N^{(1)}\wedge Y_1}=\sqrt{\sum_{1\le i<j\le 2\ell} \eta_{i,j}^2}\ge \abs{\eta_{1,\ell+1}},\]

avec 

\[\eta_{1,\ell+1}=\begin{vmatrix} 1 & 1 \\\displaystyle\sum_{k=0}^N\frac{e_k^{(1,1)}}{\theta^{\lfloor\alpha^k\rfloor}} &\displaystyle\sum_{k=0}^\infty\frac{e_k^{(1,1)}}{\theta^{\lfloor\alpha^k\rfloor}}\end{vmatrix}=\sum_{k=0}^\infty\frac{e_k^{(1,1)}}{\theta^{\lfloor\alpha^k\rfloor}}-\sum_{k=0}^N\frac{e_k^{(1,1)}}{\theta^{\lfloor\alpha^k\rfloor}}=\sum_{k=N+1}^\infty\frac{e_k^{(1,1)}}{\theta^{\lfloor\alpha^k\rfloor}}\ge \frac 1{\theta^{\lfloor\alpha^{N+1}\rfloor}}.\]

Donc

\[\norme{Z_N^{(1)}\wedge Y_1}\ge \frac 1{\theta^{\lfloor\alpha^{N+1}\rfloor}}.\]

De plus, en suivant la preuve de \eqref{maj_norme_Yj_arogmoeubv}, on démontre de même que $\norme{Z_N^{(1)}}\le c_{10}$, donc avec la majoration \eqref{maj_norme_Delta_aoerihgofbv} de la norme de $\Delta$, on a

\[\norme{Z_N^{(1)}\wedge \Delta}\le \norme{Z_N^{(1)}}\cdot\norme{\Delta}\le c_{12}\omega\]

avec $c_{12}>0$ ne dépendant que de $\ell$. Donc d'après la minoration \eqref{minoration_norme_Z_N1_proj_maoeboefb}, on obtient

\[\norme{Z_N^{(1)}\wedge p_A^\perp(Z_N^{(1)})}\ge \frac {1}{\theta^{\lfloor\alpha^{N+1}\rfloor}}-c_{12}\omega.\]

Or $\norme{p_A^\perp(Z_N^{(1)})}\le \norme{Z_N^{(1)}}\le c_{10}$, donc en utilisant le lemme \ref{lemme_proximite_proj_ortho_apamoebnamon}, on a
\begin{align*}
\omega
	&=\norme{p_A^\perp(Z_N^{(1)})-Z_N^{(1)}} \\
	&=\norme{Z_N^{(1)}}\psi(p_A^\perp(Z_N^{(1)}),Z_N^{(1)}) \\
	&=\norme{Z_N^{(1)}}\frac{\norme{p_A^\perp(Z_N^{(1)})\wedge Z_N^{(1)}}}{\norme{Z_N^{(1)}}\cdot\norme{p_A^\perp(Z_N^{(1)})}} \\
	&\ge c_{13}\norme{p_A^\perp(Z_N^{(1)})\wedge Z_N^{(1)}} \\
	&\ge \frac {c_{13}}{\theta^{\lfloor\alpha^{N+1}\rfloor}}-c_{14}\omega
\end{align*}
avec $c_{13},c_{14}>0$ ne dépendant que de $\ell$. Finalement,

\[\omega\ge \frac{c_{13}}{1+c_{14}}\cdot\frac 1{\theta^{\lfloor\alpha^{N+1}\rfloor}}=\frac {c_{15}}{\theta^{\lfloor\alpha^{N+1}\rfloor}}\]

avec $c_{15}>0$ ne dépendant que de $\ell$. Reprenons la minoration \eqref{minoration_psiell_AB_N_aroighmeobvo}, qui permet d'obtenir avec le lemme \ref{lemme_proximite_proj_ortho_apamoebnamon} :

\begin{equation}\label{minoration_psiell_oaioeanvoienvoeinvoevni}
\psi_\ell(A,B_N)\ge \psi(Z_N^{(1)},p_A^\perp(Z_N^{(1)}))=\frac{\omega}{\norme{Z_N^{(1)}}}\ge \frac {c_{16}}{\theta^{\lfloor\alpha^{N+1}\rfloor}}
\vspace{2.6mm}
\end{equation}

où $c_{16}>0$ ne dépend que de $\ell$. Or

\[\lfloor\alpha^{N+1}\rfloor\le \alpha^{N+1}\le \lfloor\alpha^N\rfloor\alpha+\alpha,\]

donc

\[\frac 1{\theta^{\lfloor\alpha^{N+1}\rfloor}}\ge \frac 1{\theta^{\lfloor\alpha^{N}\rfloor\alpha+\alpha}}=\frac 1{\theta^\alpha}\cdot\frac 1{\left(\theta^{\lfloor\alpha^{N}\rfloor}\right)^\alpha}.\]

De plus, d'après le lemme \ref{lemmeminorationhauteurdesBN_aeoinmoaeivnv} il existe une constante $c_{17}>0$ ne dépendant que de $A$ telle que

\[H(B_N)\ge c_{17}\left(\theta^{\lfloor\alpha^N\rfloor}\right)^\ell,\]

donc en reprenant la minoration \eqref{minoration_psiell_oaioeanvoienvoeinvoevni} on a

\[\psi_\ell(A,B_N)\ge \frac {c_{16}}{\theta^{\lfloor\alpha^{N+1}\rfloor}}\ge \frac {c_{16}}{\theta^{\alpha}} \cdot\frac 1{\theta^{\alpha\lfloor\alpha^N\rfloor}}\ge \frac{c_{18}}{H(B_N)^{\alpha/\ell}}\]

avec $c_{18}>0$ ne dépendant que de $A$. Ceci termine la preuve du lemme \ref{minorationdupsielldeAetBN_eogihrgoeihgoi}.
\end{preuve}

Finalement, démontrons le théorème \ref{theoreme_spectre_amoeribnefaomnv}.

\begin{preuve}\label{preuve_th_spectre_gamoerbgosn}
Supposons dans un premier temps que $n=2\ell$. Si 

\[\beta\in\left[1+\frac 1{2\ell}+\sqrt{1+\frac{1}{4\ell^2}},+\infty\right[,\]

la proposition \ref{exposant_donne_pour_A_aemofdunbvombc} donne un sous-espace $A\in\mathfrak I_n(\ell,\ell)_\ell$ tel que $\muexpA nA\ell\ell=\beta$. Si maintenant $\beta=+\infty$, posons pour $i,j\in\{1,\ldots,\ell\}$

\[\xi_{i,j}=\sum_{k=0}^\infty \frac{e_k^{(i,j)}}{3^{k^k}}\]

où les $(e_k^{(i,j)})_{k\in\N}$ sont des suites à déterminer, à valeurs dans $\{1,2\}$. En reprenant les notations de la section \ref{section_resultats_sev_exposant_donne_aefoubv}, posons $M_\xi=(\xi_{i,j})_{(i,j)\in\{1,\ldots,\ell\}^2}\in\M_\ell(\R)$ et notons $A_\infty$ le sous-espace vectoriel engendré par les colonnes de la matrice

\[\begin{pmatrix} I_\ell \\ M_\xi\end{pmatrix}\in\M_{2\ell,\ell}(\R).\]

De façon similaire au lemme \ref{lemme_condition_dirr_sur_A_aomerbvozb}, on peut choisir des suites $(e_k^{(i,j)})_{k\in\N}$ pour \hbox{$i,j\in\{1,\ldots,\ell\}$} telles que

\[A_\infty\in\mathfrak I_n(\ell,\ell)_\ell.\]

Posons pour $i,j\in\{1,\ldots,\ell\}$ et $N\ge 1$,

\[f_N^{(i,j)}=3^{N^N}\sum_{k=0}^N \frac{e_k^{(i,j)}}{3^{k^k}},\]

et notons $B_N\in\mathfrak R_{2\ell}(\ell)$ le sous-espace vectoriel rationnel engendré par les colonnes de 

\[\begin{pmatrix} 3^{N^N}I_\ell \\ F_N\end{pmatrix}\in\M_{2\ell,\ell}(\R)\]

où $F_N=(f_N^{(i,j)})_{(i,j)\in\{1,\ldots,\ell\}^2}$. De façon en tout point similaire au lemme \ref{majoration_psi2_A_BN_amroghozv}, on peut montrer qu'il existe une constante $c>0$ dépendant uniquement de $A$ telle que

\[\forall N\ge 1,\quad \psi_\ell(A_\infty,B_N)\le \frac{c}{H(B_N)^{N/\ell}}.\]

Ainsi,

\[\forall \kappa>0,\quad \forall N\ge \kappa\ell,\quad \psi_\ell(A_\infty,B_N)\le \frac{c}{H(B_N)^{\kappa}},\]

donc

\[\forall \kappa>0,\quad \muexpA {n}{A_\infty}\ell\ell\ge \kappa,\]

car en suivant la preuve du lemme \ref{les_BN_sont_les_meilleurs_afnemobibvo} (plus précisément l'inégalité \eqref{maj_psi_ell_A_B_N_aoeribqfoisbcovbz} page \pageref{maj_psi_ell_A_B_N_aoeribqfoisbcovbz}) on voit que $\psi_\ell(A_\infty,B_N)$ tend vers $0$ quand $N\to\infty$, donc il y a une infinité de sous-espaces $B_N$ deux à deux distincts. \\

Finalement,

\[\muexpA {n}{A_\infty}\ell\ell=+\infty.\]

Considérons enfin le cas où $n>2\ell$. Notons $\phi$ un isomorphisme rationnel de $\R^{2\ell}$ dans $\R^{2\ell}\times\{0\}^{n-2\ell}$. Posons $A'=\phi(A)$. D'après le théorème \ref{th_inclusion_sev_rationnel_apeivpinpiaenv}, on a $A'\in\mathfrak I_n(\ell,\ell)_\ell$ et 

\[\muexpA n{A'}\ell\ell=\muexpA {2\ell}A\ell\ell=\beta,\]

ce qui permet d'étendre le résultat aux entiers $n>2\ell$.
\end{preuve}

\newpage
\bibliographystyle{alpha-fr}
\bibliography{biblio_these}

\includepdf[pages=-]{quatrieme_de_couv_these_EJ.pdf}

\end{document}